\documentclass[msom,nonblindrev]{informs3} 

\OneAndAHalfSpacedXI

\usepackage{natbib}
 \bibpunct[, ]{(}{)}{,}{a}{}{,}%
 \def\bibfont{\small}%

\usepackage{pgfplots}
 \usepackage{comment}

\TheoremsNumberedThrough     

\EquationsNumberedThrough    

\usepackage{xcolor}

\newcommand{\umin}{U_{\textrm{min}}}
\newcommand{\umax}{U_{\textrm{max}}}

\newcommand{\C}{\vec{\mathcal{C}}}
\newcommand{\E}{\mathbb{E}}
\newcommand{\p}{\mathbb{P}}
\newcommand{\F}{\mathcal{F}}
\DeclareMathOperator{\poi}{Poi}
\DeclareMathOperator{\Mod}{\,mod\,}

\newcommand{\um}[1]{\textcolor{black}{#1}}
\newcommand{\md}[1]{\textcolor{black}{#1}}
\usepackage{amssymb}
\usepackage{graphicx}
\usepackage{booktabs}
\usepackage{tabularx}
\usepackage{threeparttable}
\usepackage{subfig}
\usepackage{tikz}

\newcolumntype{L}[1]{>{\raggedright\arraybackslash}p{#1}}
\newcolumntype{C}[1]{>{\centering\arraybackslash}p{#1}}
\newcolumntype{R}[1]{>{\raggedleft\arraybackslash}p{#1}}
\newcolumntype{J}[1]{>{\justifying\arraybackslash}p{#1}}

\begin{document}

\RUNAUTHOR{Mohring et al.} 

\RUNTITLE{Transparent time-dependent shipment pricing for e-fulfillment}

\TITLE{Same-day or next-day? Transparent time-dependent shipment pricing for e-fulfillment}

\ARTICLEAUTHORS{%
\AUTHOR{Uta Mohring}
\AFF{Department of Business Administration, University of Zurich, Switzerland, \EMAIL{uta.mohring@business.uzh.ch}}
\AUTHOR{Melvin Drent}
\AFF{Department of Information Systems and Operations Management, Tilburg University, the Netherlands, \EMAIL{m.drent@tilburguniversity.edu}}
\AUTHOR{Ivo Adan, Willem van Jaarsveld}
\AFF{School of Industrial Engineering, Eindhoven University of Technology, the Netherlands, \EMAIL{\{i.adan, w.l.v.jaarsveld\}@tue.nl}}
} 

\ABSTRACT{\textbf{Problem definition:}
Fulfilling online orders quickly is increasingly vital in e-commerce, with many companies offering same-day shipment to meet rising customer expectations.
However, firms risk overpromising when shipment options and fees are not aligned with their fulfillment capacity, leading to delayed orders and customer dissatisfaction.
Fulfillment centers may be unable to hand over all orders to parcel delivery companies by the agreed deadlines.
To avoid this, shipment options and fees should be adapted based on information that is both operationally relevant and marketable to customers—namely, the time remaining for the fulfillment center to collect and ship these orders.
\\
\textbf{Methodology/results:} 
We develop a parsimonious model of a fulfillment center with Poisson customer arrivals, stochastic processing capacity, and the objective of maximizing long-run average profit, and provide an exact steady-state analysis based on a periodic Markov chain.
As shipment fees affect both the volume and timing of same-day demand, direct analysis in the price domain is analytically intractable.
We therefore devise a problem-domain transformation that maps each shipment policy to its induced cumulative demand profile.
This transformation enables both structural insights—the optimal shipment policy includes a cutoff point and monotonically increasing fees—and tractable optimization, as the reformulated profit function is supermodular, allowing computation of the optimal shipment policy in polynomial time. We also propose a simple shipment policy with a two-level fee structure that is easy to communicate to customers and achieves near-optimal performance.\\
\textbf{Managerial implications:} 
Transparent time-dependent shipment policies provide a customer-friendly alternative to opaque dynamic pricing, allowing e-commerce companies to align same-day demand with fulfillment capacity while maintaining fairness and simplicity.
}%


\maketitle
\newpage
\section{Introduction}
Given the fast growth of e-commerce, fulfilling online orders faster and cheaper becomes more relevant than ever. 
Offering quick and convenient shipment options is a powerful marketing tool to differentiate from competitors and attract and retain customers. Empirical studies indicate that faster shipment boosts sales and profitability. For example, when promising one day faster shipment, sales for customized photo products was found to increase by 0.7\% \citep{Cui2022} while sales for apparel products increase by 1.5\% \citep{Fisher2019}. Conversely, sales decrease by 14.6\% when removing the highest-quality shipment option 
\citep{Cui2020}. Furthermore, the growing willingness of customers to pay for fast shipment, which has seen a notable rise from 41\% in 2023 to 70\% in 2024 \citep{Capgemini2025}, underscores the customers' impatience and their clear preference for quicker shipment options.   
In response, e-commerce companies increasingly promise same-day shipment, meaning orders placed until a certain cutoff point in the late afternoon or evening are promised to be shipped the same day to ensure next-day arrival. In dense urban areas, some companies even promise same-day delivery \citep{McKinsey2022}. 
From an operational perspective, however, fulfilling same-day shipment promises is challenging, especially as the end of the day approaches and the window for same-day shipment becomes tighter. 
Same-day shipment therefore requires a careful coordination between marketing and operations to ensure a smooth and efficient order fulfillment process.

Online retailers 
typically fulfill orders from so-called fulfillment centers.
In these centers, orders are picked from shelves and then prepared and consolidated into batches for shipment. 
While orders are collected and consolidated continuously throughout the day, the batches are dispatched for shipment only at predefined deadlines -- typically at the end of the day -- to meet agreed fixed handover times with external parcel delivery companies (e.g., FedEx or UPS) that handle the actual delivery of the orders \citep{Doerr2013}. Meeting these strict deadlines is crucial for the service quality perceived by customers as missing a deadline by only a few minutes may result in a delay for the customer of at least one day \citep{Ceven2017}.
As e-commerce companies offer increasingly faster service, they risk overpromising their shipment options: fulfillment centers may be unable to process the inflow of orders that require same-day shipment and need to be ready before the next deadline.
This risk is particularly high when different shipment options are priced using flat fees, giving customers no incentive to choose next-day over same-day shipment, even when fulfillment center operations would benefit from it. 

To mitigate this risk and manage limited capacity effectively, companies need to steer demand so that same-day shipment is selected primarily by customers who value it most and are willing to pay for it.
A natural approach is to introduce time- and fulfillment center load-based shipment options and shipment fees, which may provide operational advantages by helping to balance order inflow with available capacity and remaining time until handover.
While from an operational perspective appealing, such dynamic shipment policies are difficult to communicate and market effectively in e-commerce \citep{Agatz2013}. Frequent and unexpected price changes are unappealing and appear unfair to customers as they see different prices for the same service at different times without comprehensible reasons for these price differences \citep{Xia2004}. Hence, from the customer’s point of view, shipment policies that depend on the load of the fulfillment center are often perceived as opaque and unintuitive. Moreover, load-dependent pricing increases operational complexity by requiring frequent updates on the ordering platform \citep{Bergquist2025}. 

To address these concerns, we focus on a more practical and customer-focused alternative that has already seen some adoption: transparent time-dependent shipment options and corresponding fees that depend only on the remaining time until the next shipment deadline. These shipment options and fees are consistent each day and can be easily communicated on the ordering platform. Hence, any change in shipment options or fees is immediately visible and predictable to customers. Such a transparent time-dependent shipment policy is operationally relevant, simple to communicate, and perceived as predictable and fair. It thus provides a transparent alternative to opaque dynamic pricing, while still offering a lever to align customer demand with operational capacity.


In this paper, we study this transparent time-dependent shipment policy for online retailers. We do so in the context of products for which the shipment options and shipment fees mostly affect the customer's decision to choose same-day or next-day shipment, and not so much the actual decision to buy the specific product. This makes sense when companies expect little competition or when customers are loyal to the specific products and/or companies. For these settings, we explore how a transparent time-dependent shipment policy can effectively distribute customer demand between same-day and next-day shipment. Our study focuses on cutoff-based shipment options, where all orders placed before a specific cutoff point are eligible for same-day shipment, while orders placed after this point are only eligible for next-day shipment \citep{Mohring2022}.
Given the cutoff point decision, fees for same-day and next-day shipment need to be set for each moment in time. Consequently, a transparent time-dependent shipment policy raises the following interrelated questions: 
\begin{enumerate}
    \item[(1)] Cutoff point: Until what moment in time should same-day shipment be offered to customers?
    \item[(2)] \md{Shipment fees: How much should the shipment fees be?} 
    \item[(3)] \md{Differentiated fees: When and how should the shipment fees be adapted over time?} 
\end{enumerate}

We build a parsimonious model that provides answers to the questions posed above. We consider multiple operating cycles, where each cycle consists of a fixed number of time periods and ends with a deadline upon which orders that are due for shipment by this deadline need to be handed over to the parcel delivery company. Throughout each operating cycle, online customers arrive according to a Poisson process. Customers that complete their online transaction before the cutoff point are offered the choice between same-day or next-day shipment together with their corresponding fees (or, equivalently, express shipment or regular shipment). We assume that customers make this choice by trading off the utility that they derive from both shipment options. Customers that arrive after the cutoff point cannot choose same-day shipment, and hence their products are handed over to the parcel company upon the deadline at the end of the following operating cycle. The processing capacity of the fulfillment center in each time period is randomly distributed. If the available capacity is insufficient to process all orders due for shipment in a given operating cycle, then these orders are said to be late and carried over for shipment in the subsequent operating cycle. 
We are interested in characterizing transparent time-dependent shipment policies that maximize the long-run average profit, i.e., revenue minus cost for late orders. 
By developing a discrete-time periodic Markov chain model for the steady-state analysis of fulfillment centers described above, we are able to provide closed-form expressions for the relevant performance measures as well as structural properties of the transparent time-dependent shipment policy. 

The long-run average profit in our model depends on the entire temporal distribution of same-day orders, as each time-dependent shipment fee influences not only the total demand for same-day shipment but also its timing throughout the day. These interdependencies make the problem analytically challenging and render direct analysis of transparent time-dependent shipment policies in the price domain intractable. Our analysis therefore hinges on a novel problem-domain transformation that represents each policy through the cumulative demand profile it induces for same-day shipment. This transformation forms the basis for our main results: it (i) yields a supermodular reformulation that allows the optimal transparent shipment policy to be computed in polynomial time, and (ii) enables structural characterization of optimal transparent shipment policies, showing that same-day fees increase monotonically over time and that introducing a cutoff point and replacing a constant fee with a time-dependent fee improve profit.

Beyond these analytical contributions, we compare specific classes of transparent shipment policies to derive simple, implementable guidelines for e-commerce companies. Our numerical analysis demonstrates that a policy with a two-level fee structure achieves near-optimal performance while remaining fully transparent and easy to communicate to customers. These results show that e-commerce companies can derive substantial benefits from introducing transparent time-dependent shipment options and fees, balancing operational efficiency with customer trust.

The remainder of this paper is structured as follows: \S\ref{sec:Literature} positions our work within the related literature. In \S\ref{sec:ModelAnalysis}, we formalize the problem and introduce the Markov chain model. In \S\ref{sec:analyticalResults}, we establish structural properties of the transparent time-dependent shipment policy and propose the two-level time-dependent shipment policy. In \S\ref{sec:numericalAnalysis}, we demonstrate the performance of the two-level time-dependent shipment policy in a numerical study. \S\ref{sec:Conclusions} concludes the paper.

\section{Literature review}
\label{sec:Literature}
In terms of modeling, our work is related to the literature on pricing in queueing systems. In terms of application, it is related to the literature on order fulfillment operations, in particular, the streams on pricing decisions and deadline-oriented operations therein. In the following, we review the relevant papers and position our work in each of these literature streams.

\subsection{Pricing in queueing systems}
The literature on pricing in single-server queues studies dynamic and static pricing policies to regulate the total demand rate and customers' service class choices. 
We refer the reader to \cite{HassinHaviv2003} and \cite{Hassin2016} for comprehensive surveys. 
The stream on dynamic pricing, which builds on the seminal work of \cite{Low1974}, models prices as load-dependent variables that are continuously updated depending on the current queue length. Even if the queue length is observable for customers, it is intransparent to them how prices change as the system workload changes. As dynamic pricing policies are opaque and load-dependent, they differ fundamentally from the transparent time-dependent pricing policy studied in this paper.


More related to our work is the stream on static priority pricing, which includes \cite{MendelsonWhang1990}, \cite{Afeche2013}, \cite{AfecheBaronKerner2013}, \cite{AfechePalvin2016}, \cite{GavirneniKulkarni2016}, \cite{Sainathan2018}, \cite{NazerzadehRandhawa2018}, and \cite{CaoWangXie2019}. These papers study the design of the price/lead-time menu, i.e., the menu of service classes, each characterized by a price and an expected lead time. As prices and lead times do not depend on the system workload, the price/lead-time menu is transparent for customers, similar to our shipment policy. Customers are price- and delay-sensitive and select their preferred service class based on the posted menu, their service valuation, and their delay cost. If customers also have the option to balk, the service service is called non-captive. 
Our shipment policy can be seen as a price/lead-time menu with two service classes -- same-day and next-day shipment -- for captive shipment service in online retail. However, unlike the common notion of lead times as the time period between arrival and service completion of a customer order, lead times in e-fulfillment centers are deadline-oriented. That is, service classes promise service completion (shipment) by certain predefined and periodically recurring deadlines rather than an expected lead time until service completion. As a consequence, a time-dependent price/lead-time menu that depends on the remaining time until the next deadline appears to be an effective tool to regulate customers' service class choices. Our work is the first to study priority pricing for deadline-oriented service classes.

\subsection{Pricing as part of order fulfillment operations}
The order fulfillment literature discusses pricing in the context of order fulfillment optimization, last-mile delivery, and shipment pricing. 
The stream on order fulfillment optimization studies the optimal assignment of customer orders to fulfillment centers and includes papers \cite[e.g.,][]{Lei2018,Harsha2019,Lei2022} that use inventory level-based pricing -- often opaque from the customer’s perspective -- to mitigate demand-inventory-imbalances. 
We refer the reader to \cite{Acimovic2019} for a comprehensive review.

The stream on last-mile delivery includes papers that integrate vehicle routing problems with pricing decisions and focuses on workload-dependent pricing policies that are opaque from the customer's perspective. We refer the reader to \cite{Wasmuth2022} and \cite{Fleckenstein2023} for comprehensive reviews. The closest to our work is \cite{Banerjee2025}, who study a transparent pricing policy with price-differentiated same-day and next-day delivery options. 
While this paper focuses on two specific pricing schemes with a single and a two-level price for same-day delivery, which are consistent to our constant cutoff point and two-level time-dependent shipment policies, we provide a comprehensive performance analysis of general time-dependent pricing policies as well as an efficient algorithm to determine the optimal policy. Moreover, we study the benefit of time-dependent pricing for operations in the fulfillment center rather than vehicle routing and last-mile delivery.

Most related to our work is the stream on shipment pricing, which studies price partitioning strategies, contingent free shipment policies, and subscription plans for free shipment. 
Papers on price partitioning \citep[e.g.,][]{Gumus2013,Jiang2013} study under what conditions to bundle product price and shipment fee, and if not to bundle, how to set prices for product and shipment separately. 
Contingent free shipment policies offer free shipment to customers if their order exceeds a predefined threshold, otherwise customers pay for shipment. Related papers include \cite{Leng2005}, \cite{Leng2010}, and \cite{LiShengZhan2023}.
Shipment subscription plans charge an upfront fee and offer free shipment for a limited/unlimited number of orders throughout the subscription period. Related papers are \cite{Fang2021}, \cite{Sun2022}, and \cite{Balakrishnan2024}. 
All of the aforementioned papers analyze the benefit of shipment pricing from a marketing perspective. However, they ignore the operational challenges at the fulfillment center resulting from customers' short and heterogeneous delivery time expectations.  

\cite{Yao2012} and \cite{Sainathan2018} are the only papers that  study shipment pricing from an operational perspective as we do.
\cite{Yao2012} characterize the optimal price partitioning strategy and shipment pricing scheme of an online retailer that differentiates between expedited and standard shipment. Both the game-theoretical model and the empirical analysis reveal that the shipment price premium of expedited shipment over standard shipment increases in the standard shipment time, but decreases in the expedited shipment time and the retailer's on-time probability. 
\cite{Sainathan2018} studies pricing of an ancillary service (e.g., shipment in e-commerce) where customers are heterogeneous in terms of service valuation and delay. The author shows how the optimal pricing strategy -- free service, single service, or differentiated service -- depends on capacity restrictions and product type. 
Our work goes beyond these two papers as we analyze the benefit of transparent time-dependent pricing for expedited shipment 
to better balance the incoming customer orders with the scarce processing capacity to fulfill these on time.  

\subsection{Deadline-oriented order fulfillment operations}
One of the key challenges of fulfillment centers is that orders arrive and are processed continuously throughout the day while shipments are dispatched at predefined deadlines at the end of the day \citep{Doerr2013}. Nevertheless, deadline-oriented order fulfillment operations have received limited attention in the literature. \cite{Lan2024} study a same-day delivery problem in which vehicles are dispatched at predefined deadlines (e.g., once per hour). The authors determine at each decision epoch, which orders to dispatch and which to postpone to consolidate with potential future orders. 

More related to our work are papers that focus on operations in the fulfillment center.  
\cite{Doerr2013} propose the performance metric \textit{Next Scheduled Deadline} (NSD), which  measures the fraction of orders targeted for a particular deadline (truck departure) that are ready by this deadline. An order is targeted for a deadline if it arrives before the associated cutoff point. The authors study how to steer the cutoff point and the percentage of NSD in order to improve worker motivation and accelerate the operating speed when it matters most. 
\cite{Ceven2017} and \cite{MacCarthy2019} study wave picking strategies and determine the optimal number and timing of picking waves. The authors apply NSD as the performance metric for the service quality perceived by customers. \cite{Ceven2017} focus on fulfillment centers with daily deadlines, while \cite{MacCarthy2019} study the fulfillment process of buy-online-pickup-in-store retail services, where online orders are fulfilled in retail stores while also serving walk-in customers. 
\cite{Mohring2022} explore differentiated time-dependent shipment options that depend on the remaining time until the deadline. Specifically, same-day shipment is promised to all customers that order until the cutoff point and next-day shipment to all customers that order thereafter. The authors discuss the benefit of such cutoff-based shipment options and how to set the cutoff point to balance customer expectations in terms of order response time and service level. 

Although these models address challenges of deadline-oriented order fulfillment operations and discuss time-dependent shipment options, little research considers shipment menus from which the customers select their preferred option, and more importantly, suitable shipment pricing schemes. Our paper addresses this research gap.

\section{Model}
\label{sec:ModelAnalysis}
In \S\ref{sec:ProblemDescription}, we formalize problem setting and shipment policy. In \S\ref{sec:Model}, we introduce a Markov chain model for steady-state performance analysis. 

\subsection{Problem description}
\label{sec:ProblemDescription}
We study a stylized stochastic e-commerce fulfillment center with a transparent time-dependent shipment policy. In this section, we formalize the order fulfillment process, the shipment policy, and the customer choice model.   

\subsubsection*{Order fulfillment process.}
We consider the order fulfillment process of a company operating its own warehouse but outsourcing the deliveries to a dedicated parcel delivery company. That is, the order fulfillment process starts with receiving customer orders and includes picking, packing products to parcels, preparing parcels for shipment, and consolidating parcels into batches for shipment depending on their destinations. The order fulfillment process ends with handing over the parcels to the parcel delivery company that operates with a fixed schedule of truck departure times at the fulfillment center. Meeting the truck departure time is crucial for the order fulfillment process as missing a truck departure by only a few minutes may result in a delay to the customer of at least one day \citep{Doerr2013,Ceven2017}.

These recurring deadlines subdivide the discrete time axis $t \in \mathbb{N}_0$ into buckets, e.g. days, which we call \emph{operating cycles}. Operating cycles consist of $T$ time periods each and are indexed by $l\in\mathbb{N}_0$. Operating cycle $l$ corresponds to the time periods $\{l T, lT+1,\ldots,l T +T-1 \big\}$; it starts at deadline $lT$ and ends immediately before reaching the subsequent deadline $(l+1)T$; see Figure~\ref{fig:OperationPeriod} for illustration. 
For every time $t \in \mathbb{N}_0$, the corresponding operating cycle is
\begin{equation}
    l = \left\lfloor \frac{t}{T} \right\rfloor, 
\end{equation}
and the age $\tau$ within operating cycle $l$ is obtained as:
\begin{equation}
    \tau = t\Mod T. 
\end{equation} 

\begin{figure}[t]
\centering
    \begin{tikzpicture}[scale = 0.82]\small
        \draw[thick, -stealth] (0.25,0) -- (18.5,0);
        \draw(18.5,-0.4) node{time};
        
        \draw (2.25, 2.2) arc[radius = 0.2, start angle = 180, end angle = 90];
        \draw (2.45, 2.4) -- (9.05, 2.4);

        \draw (9.05, 2.4) arc[radius = 0.2, start angle = 270, end angle = 360];
        \draw (9.45, 2.4) arc[radius = 0.2, start angle = 270, end angle = 180];
        
        \draw (9.45, 2.4) -- (16.05, 2.4);
        \draw (16.05, 2.4) arc[radius = 0.2, start angle = 90, end angle = 0];
        
        \draw (9.25, 3.0) -- (9.25, 3.0) node{Operating cycle $l$};
        \draw (1.025, 3.0) -- (1.025, 3.0) node{$l-1$};
        \draw (17.575, 3.0) -- (17.575, 3.0) node{$l+1$};
        
        \draw[dashed] (2.05, 2.4) arc[radius = 0.2, start angle = 90, end angle = 0];
        \draw[dashed, stealth-] (0.25, 2.4) -- (2.05, 2.4);
        
        \draw[dashed] (16.25, 2.2) arc[radius = 0.2, start angle = 180, end angle = 90];
        \draw[dashed, -stealth] (16.55, 2.4) -- (18.5, 2.4);

        \draw (0.5, -0.15) -- (0.5, 0.15);
        \draw[thick] (2.25, -0.3) -- (2.25, 0.3);
        \draw (4.0, -0.15) -- (4.0, 0.15);
        \draw (5.75, -0.15) -- (5.75, 0.15);
        \draw (7.5, -0.15) -- (7.5, 0.15);
        \draw (9.25, -0.15) -- (9.25, 0.15);
        \draw (11.0, -0.15) -- (11.0, 0.15);
        \draw (12.75, -0.15) -- (12.75, 0.15);
        \draw (14.5, -0.15) -- (14.5, 0.15);
        \draw[thick] (16.25, -0.3) -- (16.25, 0.3);
        \draw (18.0, -0.15) -- (18.0, 0.15);
        
        \draw (0.5, 0.7) -- (0.5, 0.7) node{$\tau:$};
        \draw (2.25, 0.7) -- (2.25, 0.7) node{$0$};
        \draw (4.0, 0.7) -- (4.0, 0.7) node{$1$};
        \draw (14.5, 0.7) -- (14.5, 0.7) node{$T-1$};
        \draw (16.25, 0.7) -- (16.25, 0.7) node{$0$};
        
        \draw (0.5, 1.15) -- (0.5, 1.15) node{$t:$};
        \draw (2.25, 1.15) -- (2.25, 1.15) node{$lT$};
        \draw (4.0, 1.15) -- (4.0, 1.15) node{$lT+1$};
        \draw (14.3, 1.15) -- (14.3, 1.15) node{$lT + T-1$};
        \draw (16.3, 1.15) -- (16.3, 1.15) node{$(l+1)T$};
        \draw (2.25, 1.8) -- (2.25, 1.8) node{Deadline};
        \draw (16.25, 1.8) -- (16.25, 1.8) node{Deadline};
        
        \draw[dashed] (2.25, -0.4) -- (2.25, -0.8);
        \draw[dashed] (16.25, -0.4) -- (16.25, -0.8);
        
        
        \draw[thick, -stealth] (2.75, -1.0) -- (16.25, -1.0);
        \draw[thick, stealth-] (2.25, -1.0) -- (2.75, -1.0);
        \draw (2.25, -0.9) -- (2.25, -1.1);
        \draw (16.25, -0.9) -- (16.25, -1.1);
        
        \draw(9.25, -0.7) -- (9.25, -0.7) node{$T$};

    \end{tikzpicture}
    
    \caption{Representation of the time line and corresponding notation.} 
    \label{fig:OperationPeriod}
\end{figure}

We denote by $D_t$ the random \emph{customer demand} at time $t\in\mathbb{N}_0$, i.e., the number of orders released for order fulfillment at time $t$. Random variables $D_t$, $t \in \mathbb{N}_0$, are i.i.d. drawn from a Poisson distribution with expected customer demand~$\lambda$:
\begin{equation}
    D_t \sim \poi(\lambda).
\end{equation}
Most of our analytical results continue to hold for time-dependent customer demand. Related model extensions are discussed in Appendix \ref{sec:modelExtension}.

We model the order fulfillment process as a single-stage system, which includes all substeps of the order fulfillment process. Let $K_t$ denote the random \emph{processing capacity} available at time $t\in\mathbb{N}_0$. That is, at most $K_t$ orders are processed at time $t$. 
Random variables $K_t$, $t \in \mathbb{N}_0$, are i.i.d., independent of $D_t$, and specified by the discrete generally-distributed random variable~$K$ with support $\mathbb{N}_0$:
\begin{equation}
    K_t \sim K.
\end{equation}
All our analytical results continue to hold for time-dependent processing capacity.

The \emph{utilization} of the fulfillment center, denoted by $\rho$, is the ratio of the expected customer demand $\lambda$ by the expected processing capacity $\E[K]$, i.e., 
\begin{equation}
    \rho=\frac{\lambda}{\E[K]}.
    \label{eq:U}
\end{equation}
We focus on the case $\rho < 1$, which ensures system stability.

For ease of reference, Table~\ref{tab:parameters} lists the notation introduced above, as well as the notation introduced in the remainder of this section. 
\begin{table}[h]
    \small
    \begin{center}
        \caption{Model notation.}
        \label{tab:parameters}
        \begin{tabular}{lL{7.1cm}l}
            \toprule
             & Name & Remark\\
            \midrule
            \multicolumn{3}{l}{\textbf{Time axis}} \\
            $t$ & Time & $t\in\mathbb{N}_0$ \\
             $T$ & Duration of each operating cycle & $T \in \mathbb{N}$ \\
             $l$ & Operating cycle & $l=\lfloor t/T \rfloor \in \mathbb{N}_0$ \\
             $\tau$ & Age of operating cycle  & $\tau= t \Mod T$\\
           \midrule
            \multicolumn{3}{l}{\textbf{Order fulfillment}} \\
            $\lambda$ & Expected customer demand per time interval & $\lambda >0$\\
            $D_{t}$ & Customer demand at time $t$ & $D_t\sim \poi(\lambda)$\\
            $K$ & Processing capacity per time interval & Generic random variable w. support $\mathbb{N}_0$\\
            $K_t$ & Processing capacity at time $t$ &$K_t \sim K$\\
            $\rho$ & Utilization & $\rho = \lambda/\E[K]$ \\ 
            \midrule
            \multicolumn{3}{l}{\textbf{Shipment policy}} \\
            $\bar{p}$ & Product price  & $\bar{p} \in \mathbb{R}_{\geq0}$\\
            $\bar{\tau}$ & Cutoff point & $\bar{\tau}\in \{0,1,\ldots,T-1\}$ \\
            $f_{\tau}$ & Express shipment fee at age $\tau$ & $f_\tau \in \mathbb{R}_{>0}$ \\  
            $f_{t}$ & Express shipment fee at time t & $f_t = f_\tau$ with $\tau=t\Mod T$ \\  
            \midrule
            \multicolumn{3}{l}{\textbf{Customer choice}} \\
            $v$ & Product valuation & $v\geq \bar{p}$ \\
            $\umin$ & Minimum express shipment valuation & $\umin \in \mathbb{R}_{\ge 0}$\\
            $\umax$ & Maximum express shipment valuation  & $\umax \in [\umin,\infty)$\\
            $U$ & Express shipment valuation & Uniform random variable w. support $[\umin,\umax]$\\
            $w(f)$ & Willingness proportion at shipment fee $f$ & $w: \mathbb{R}\rightarrow [0,1]$\\ 
            \midrule
            \multicolumn{3}{l}{\textbf{Express and regular shipment orders}} \\
            $E_\tau$ & Express shipment orders at age $\tau$ & $E_\tau \sim \poi(w(f_\tau)\lambda)$ if $\tau\leq \bar{\tau}$, $E_\tau\equiv0$ otherwise\\
            $E_t$ & Express shipment orders at time $t$ & $E_t\sim E_\tau$ with $\tau=t\Mod T$\\
            $N_\tau$ & Regular shipment orders at age $\tau$ & $N_\tau \sim \poi((1-w(f_\tau))\lambda)$ if $\tau\leq \bar{\tau}$, \\
            & & $N_\tau\sim\poi(\lambda)$ otherwise\\
            $N_t$ & Regular shipment orders at time $t$ & $N_t\sim N_\tau$ with $\tau=t\Mod T$\\
            \bottomrule
        \end{tabular}
    \end{center}
    \end{table}

\subsubsection*{Shipment policy.} 
The shipment policy specifies the shipment options offered to customers and the fees charged for these shipment options. There are potentially two shipment options available to customers: \emph{express shipment} and \emph{regular shipment}.
More specifically, consider a customer arriving at time $t$ in operating cycle $l=\lfloor t/T\rfloor$. Then, express shipment corresponds to an order dispatch at the next deadline after $t$, i.e., at time $(l+1)T$, while regular shipment coincides with an order dispatch at the second deadline after $t$, i.e., at time $(l+2)T$. Thus, with a one-day operating cycle ($T=24$ hours), express shipment means shipping today for delivery tomorrow, while regular shipment mean shipping tomorrow for delivery in two days. 

To manage the inflow of express shipment orders and guarantee a minimum time buffer between their arrival and promised shipment date, we introduce the \emph{cutoff point}, denoted by $\bar{\tau}$. Express shipment is available to customers at time $t$ only if the age $\tau = t \Mod T$ of the operating cycle is weakly smaller than the cutoff point, i.e., if $\tau \le \bar{\tau}$. Regular shipment is available at every time $t$. 

Customers pay the \emph{product price}, denoted by $\bar{p}$, which is independent of the selected shipment option. In addition, they pay a shipment fee that depends on the selected shipment option, motivated by the general willingness of e-commerce customers to pay for faster shipment \citep{Capgemini2025}. We model a zero fee for regular shipment and a time-dependent positive fee for express shipment. 
Mathematically, this is equivalent to providing a discount for regular shipment, which in practice is often easier to market (e.g., Amazon).
Let $f_t$ denote the \emph{express shipment fee} at time~$t$. The fee $f_t$ depends on time $t$ by the age $\tau$ of the operating cycle, i.e., $f_t = f_\tau$ with $\tau = t \Mod T$.

Thus, the shipment policy specifies the cutoff point $\bar{\tau}$ and $(\bar{\tau}+1)$ express shipment fees $f_\tau$, $\tau = 0, 1, \dots, \bar{\tau}$. We denote the shipment policy by $\pi: = (\bar{\tau},f_0,f_1,\dots,f_{\bar{\tau}})$. Note that $\pi$ is time-dependent within each operating cycle as the shipment options offered and the fees charged for express shipment depend on the age $\tau$ of the operating cycle. However, since shipment options and fees are consistent in each operating cycle, $\pi$ can be easily communicated on the ordering platform and is transparent to customers.

\subsubsection*{Customer choice.}
\label{sec:DescriptionCustomerChoice}
Let $v$ be a customer's \emph{product valuation}, i.e., the valuation of owning the product, and $U$ denotes a customer's \emph{express shipment valuation}, i.e., the willingness to pay for receiving express instead of regular shipment. Then, a customer's \emph{total service valuation} is $v+U$ for express shipment and $v$ for regular shipment. As focusing on the choice between express and regular shipment, we assume, for clarity of results, that every customer that arrives at any time $t$ also buys the product, i.e., $v\geq \bar{p}$. The express shipment valuations of the customers are i.i.d. drawn from a uniform distribution with support $[\umin,\umax]$. 

Consider a customer who arrives at time $t$, i.e., at age $\tau = t\Mod T$ of operating cycle $l=\lfloor t/T\rfloor$. If $\tau >\bar{\tau}$, regular shipment is the only available shipment option and the customer chooses regular shipment. If $\tau \le \bar{\tau}$, express shipment is available at a price of $\bar{p}+f_\tau$ and regular shipment at a price of $\bar{p}$. The customer derives a utility of $v+U-\bar{p}-f_\tau$ from express shipment and $v-\bar{p}$ from regular shipment and chooses express shipment if 
\begin{equation*}
    v+U-\bar{p}-f_\tau > v-\bar{p} \quad \Leftrightarrow\quad  U > f_\tau.
\end{equation*}

Thus, the \emph{willingness proportion}, denoted by $w$, that gives the proportion of customers choosing express shipment is a function of the express shipment fee $f$, i.e., 
\begin{equation} 
    w(f) = \begin{cases}
        \p(U>f) &\text{if } f \in (\umin,\umax) \\
        1   &\text{if } f \le \umin \\
        0   &\text{if } f \ge \umax,
        \end{cases} 
        \quad \quad \quad \quad \text{with } \p(U>f)= \frac{\umax - f}{\umax - \umin}.
        \label{eq:willingnessToPay}
\end{equation} 
Note that $w$ is linear in $f$ and strictly decreases for $f \in [\umin,\umax]$. All customers choose express shipment if the fee is sufficiently small, i.e., $f \leq \umin$. Conversely, no customer chooses express shipment if $f \geq \umax$.

\subsubsection*{Express and regular shipment orders.}
Depending on the posted shipment policy, the customer demand $D_t$ arriving at time $t$ divides into \emph{express shipment orders}, denoted by $E_t$, and \emph{regular shipment orders}, denoted by $N_t$. Recalling that express shipment is available at a price of $\bar{p}+f_\tau$ if $\tau=t\Mod T \le \bar{\tau}$, we find together with (\ref{eq:willingnessToPay}) that 
\begin{equation}
    \E[E_t] = \begin{cases}
       w(f_\tau)\lambda & \text{if } \tau = t\Mod T \leq \bar{\tau}, \\
       0                    & \text{otherwise}, \\
    \end{cases}
    \quad \quad  \quad 
    \E[N_t] = \begin{cases}
       \left(1-w(f_\tau)\right)\lambda & \text{if } \tau = t\Mod T \leq \bar{\tau}, \\
      \lambda                   & \text{otherwise}. \\
    \end{cases}
    \label{eq:expressRegularOrdersRealization}
\end{equation}
Note that random variables $E_t$, $t\in\mathbb{N}_0$, and $N_t$, $t\in\mathbb{N}_0$, are Poisson-distributed and independent of each other \citep[Corollary 5.9]{Last2017}. As $E_t$, $t\in\mathbb{N}_0$, and $N_t$, $t\in\mathbb{N}_0$, depend on time~$t$ only by the age $\tau=t \Mod T$ of the operating cycle, we have $E_t\sim E_\tau$ and $N_t\sim N_\tau$, where
\begin{align}
    E_\tau  &\sim \begin{cases}
        \poi\big( w(f_\tau)\lambda\big) & \text{if }\tau \leq \bar{\tau}, \\
        0  & \text{otherwise} \\
    \end{cases}
    \quad \quad \quad \quad &\forall &\tau \in\{0,1,\dots, T-1\}, \label{eq:expressOrdersDistribution} \\
    N_\tau &\sim 
    \begin{cases}
        \poi\big((1-w(f_\tau))\lambda\big) & \text{if } \tau \leq \bar{\tau}, \\
        \poi\big(\lambda\big)                 & \text{otherwise} \\
    \end{cases}
    \quad \quad \quad \quad &\forall &\tau \in\{0,1,\dots, T-1\}. \label{eq:regularOrdersDistribution}    
\end{align}

\subsection{Model}
\label{sec:Model}
In this section, we introduce a discrete-time Markov chain model for steady-state performance analysis of an e-commerce fulfillment center with a transparent time-dependent shipment policy. 

\subsubsection{Discrete-time Markov chain}
\label{sec:ModelMarkovChain}
Let $t \in \mathbb{N}_0$ and $l=\lfloor t/T\rfloor$. We denote by $X^h_t$ the \emph{high-urgency order backlog} at time $t$, i.e., the number of unprocessed orders at time $t$ that are due by the next deadline after $t$, i.e., at time $(l+1)T$, or that are already too late. Further, $X^\Sigma_t$ denotes the \emph{total order backlog} at time $t$, i.e., the total of unprocessed orders at time $t$. 

The dynamics of the order backlog $(X^h_t, X^\Sigma_t)$ at time $t$ depends on the age $\tau = t\Mod T$ of the operating cycle. Specifically, the customer demand $D_t$ arriving at time $t$ divides into the express shipment orders $E_t$ and regular shipment orders $N_t$ as specified in (\ref{eq:expressOrdersDistribution})-(\ref{eq:regularOrdersDistribution}). Noting that the order fulfillment process of express and regular shipment orders is identical, but the orders differ in their promised shipment date, the processing capacity $K_t$ at time $t$ gives strict priority to high-urgency order backlog and express shipment orders. Any remaining processing capacity is used to process regular shipment orders or remains unused in case the system runs idle. 

For the dynamics of the total order backlog $X^\Sigma_t$, note that the total of unprocessed orders at time $t$ equals $X^\Sigma_t+D_t$. Any order of $X^\Sigma_t+D_t$ that remains unprocessed in period $t$ is carried over to the next period $t+1$. Thus, we have $X^\Sigma_{t+1}= \left( X^\Sigma_t + D_t -K_t \right)^{+}$, where $(x)^+:=\max(x,0)$. 

For the dynamics of the high-urgency order backlog $X^h_{t}$, recall that the high-urgency unprocessed orders at time $t$ are due by the next deadline after time $t$, i.e., at time  $(l+1)T$, or already too late. Their number equals $X^h_t+E_t$ and any order of $X^h_t+E_t$ that remains unprocessed in period $t$ is carried over to the next period $t+1$. We have $X^h_{t+1}= \left( X^h_t + E_t -K_t \right)^{+}$. 
If $t$ is the last period of operating cycle $l$, i.e., $t \Mod T = T-1$, the next period $t+1$ coincides with deadline $(l+1)T$ and the next deadline after time $t+1$ is $(l+2)T$. Thus, in addition to $X^h_t+E_t$, any order of $N_t$ that remains unprocessed in period $t$ is due by this deadline $(l+2)T$ or already too late. That is, $X^h_{t+1}=\left( X^\Sigma_t + D_t -K_t \right)^{+} = X^{\Sigma}_{t+1}$ if $t \Mod T = T-1$. 

Following these considerations, the transition of the order backlog from $(X^h_t,X^\Sigma_t)$ at time $t$ to $(X^h_{t+1},X^\Sigma_{t+1})$ at time $t+1$ is given by
\begin{align}
     X_{t+1}^h &= \begin{cases}
            \left(X^\Sigma_t + D_t - K_t \right)^+    & \text{if }t\Mod T = T-1, \\
            \left(X^h_t + E_t - K_t \right)^+     &\text{otherwise},   \\
    \end{cases} 
    \label{eq:transitionHighUrgencyBacklog}\\
    X_{t+1}^\Sigma &= \left(X^\Sigma_t + D_t - K_t \right)^+. \label{eq:transitionTotalBacklog}
\end{align}
As these transitions depend on time $t$ only by the age $\tau=t\Mod T$ of the operating cycle, $(X^h_t,X^\Sigma_t,\tau)$ with $\tau = t \Mod T$ represents the state of the underlying periodic Markov chain completely.

\subsubsection{Performance measures}
\label{sec:ModelPerformanceMeasures}
The (periodic) Markov chain that arises from any given shipment policy $\pi = (\bar{\tau},f_0,f_1,\ldots,f_{\bar{\tau}})$ has a stationary probability distribution since $\rho<1$. Therefore, defining performance measures in terms of long-run averages is appropriate.

\subsubsection*{Backorders.}
Let $t \in \mathbb{N}_0$ and $l=\lfloor t/T\rfloor$. Orders that are due by the deadline immediately after operating cycle $l$, i.e., at time $(l+1)T$, become backorders if they remain unprocessed in the last period of operating cycle $l$. We denote by $B_l$ the random number of \emph{backorders} at the end of operating cycle $l$. Backorders occur at the end of operating cycle $l$ if at time $t=lT+(T-1)$, the number of high-urgency unprocessed orders $X^h_t+E_t$ exceeds the processing capacity $K_t$. More specifically, we have 
\begin{equation} \label{eq:backorders}
    B_l=\big(X^h_{lT+(T-1)} + E_{lT+(T-1)} - K_{lT+(T-1)} \big)^+.
\end{equation}
The long-run average number of backorders per operating cycle equals 
\begin{equation}
\E[B]=\lim_{L\rightarrow\infty}\frac{1}{L}\sum_{l=0}^{L-1}\E[B_l].
\end{equation}

To motivate our later choice to penalize backorders, we seek to relate the backorders to a more customer-oriented performance measure: the average delay experienced by customers relative to the promised shipment date. Consider a customer who arrives at time $t$. The promised shipment date is $(l+1)T$ if the customer chooses express shipment and $(l+2)T$ otherwise. Further, each backorder in $B_l$ was promised to be shipped at or before $(l+1)T$, but will be shipped at or after $(l+2)T$. Thus, by \citet{little1961proof}, the average delay that a customer experiences relative to the promised shipment date equals $T\E[B]/(T\lambda)=\E[B]/\lambda$.

\subsubsection*{Profit.}
The profit of the e-fulfillment center is determined by revenues from the product price and shipment fees, net of penalty costs incurred for backorders. Noting that the product price revenue is independent of the shipment policy, it is omitted in the following. Let $R_l:=\sum_{\tau=0}^{T-1} f_\tau E_{lT+\tau}$ denote the \emph{revenue} earned in operating cycle $l$. The expected revenue per operating cycle is 
\begin{align}
    \E[R] &= \sum_{\tau=0}^{{T-1}} f_\tau \E[E_\tau]. \label{eq:revenue} 
\end{align}

The penalty costs quantify the loss of goodwill caused by orders that are dispatched for shipment after their promised shipment date. Let $\beta\in\mathbb{R}_{>0}$ be the penalty cost rate per backorder and time period. Then, $G_l:= R_l - \beta B_l$ denotes the \emph{profit} earned in operating cycle $l$. 
The expected profit per operating cycle is
\begin{align}
    \E[G] &= \E[R]  - \beta \E[B]. \label{eq:ProfitVariable}
\end{align}

\section{Structural properties of the shipment policy}
\label{sec:analyticalResults}
In this section, we establish analytical properties of the shipment policy. We are particularly interested in the effects of the shipment policy on the expected profit to gain structural insights into superior shipment policies. All proofs are relegated to Appendix \ref{AppendixProofs}. 
We use superscript $\pi$ to denote the dependence of variables on the shipment policy.
\begin{definition}[Set of shipment policies] \label{def:setShipmentPolicies}
    The set of shipment policies is given by
    \begin{equation*} 
        \Pi := \left\{ (\bar{\tau},f_0, f_1, \dots, f_{\bar{\tau}}) \mid  \bar{\tau} \in \{0,1,\dots,T-1\}, f_\tau \in [\umin, \umax]\, \forall \tau = 0,1,\dots,\bar{\tau} \right\}.
    \end{equation*}
    \hfill $\Box$
\end{definition}

Recall from (\ref{eq:ProfitVariable}) the expected profit $\E[G] = \E[R]  - \beta \E[B]$ and note that the shipment policy $\pi$ affects both the revenue and penalty cost component of the profit.
Analyzing the effects of $\pi$ on the expected revenue $\E[R]$ is straightforward and connects to well-known results of the revenue management literature \citep[\S5.2.1.2]{Talluri2006}. As $\E[R] = \lambda \sum_{\tau=0}^{\bar{\tau}} f_\tau w(f_\tau)$ by (\ref{eq:revenue}) and (\ref{eq:expressRegularOrdersRealization}), a revenue-maximizing e-commerce company sets the cutoff point $\bar{\tau}$ to the last time period of the operating cycle and charges a constant express shipment fee $f^r$, which satisfies $f^r = -w(f^r)/w'(f^r)$. That is, it offers express shipment for a constant fee of $f^r$ at every age $\tau$ of the operating cycle. These observations are formalized as follows:
\begin{lemma}[Revenue-maximizing shipment policy] \label{lem:revenueMaximizingPolicy}
Let $\pi=(\bar{\tau},f_0, f_1, \dots, f_{\bar{\tau}})\in\Pi$.
    \item[(1)] The expected revenue $\E[R^{\pi}]$ is separable and concave in the shipment fees $f_\tau$, $\tau = 0,1,\dots, \bar{\tau}$.
    \item[(2)] The unique revenue-maximizing shipment policy is $\pi^r=(T-1,f^r,\dots, f^r)$ with $f^r = \max\{\umax/2; \umin\}$.
\end{lemma}

The effects of $\pi$ on the penalty costs are more subtle. The penalty costs are linear in the expected backorders $\E[B]$, see (\ref{eq:ProfitVariable}), and $\E[B]$ depends on $\pi$ through the express shipment orders $E_\tau$, $\tau=0,1,\dots,\bar{\tau}$, see (\ref{eq:backorders}). In addition to the total express shipment orders, their distribution within the operating cycle is critical. To see this, consider two scenarios with the same total express shipment orders, where all these orders arrive at $\tau=0$ in the first scenario and at $\tau=T-1$ in the latter one. Since the processing capacity in each time period is limited and unused capacity is lost, the processing capacity of $T$ time periods is available to process the express shipment orders in the first scenario, whereas in the second scenario, it is only the processing capacity of time period $T-1$. Thus, the processing capacity can be better utilized to process express shipment orders in the first scenario, resulting in fewer expected backorders. 

To formalize the temporal distribution of express shipment orders, we introduce:
\begin{definition}[Express shipment demand profile]\label{def:demandProfile}
    Let $\pi\in\Pi$. We denote by $\mathbf{C}^\pi:= (C_{-1}^\pi, C_{0}^\pi, C_{1}^\pi,\dots, C_{T-1}^\pi)$ the \textit{express shipment demand profile} of shipment policy $\pi$, where $C_{-1}^\pi := 0$ and $C_\tau^\pi := \sum_{j=0}^\tau\E[E_j^\pi]$ $\forall \tau\in\{0,1,\dots,T-1\}$. I.e., $C_\tau^\pi$ denotes the expected cumulative express shipment orders up to age $\tau$, and $C_{T-1}^\pi$ the expected express shipment orders per cycle.
    \hfill $\Box$
\end{definition}
Clearly, for any $\pi\in\Pi$, the demand profile $\mathbf{C}^\pi$ is non-decreasing in $\tau$ and satisfies
\begin{align}
    C_\tau^\pi - C_{\tau-1}^\pi = \E[E_\tau^\pi] &= \begin{cases}
        w(f_\tau)\lambda &\text{if }\tau \leq \bar{\tau}, \\
        0 & \text{otherwise}
    \end{cases}\quad \quad \quad \quad \forall \tau \in \{0,1,\dots,T-1\}. 
     \label{eq:demandProfileShipmentPolicy}
\end{align}
An example demand profile is shown in Figure~\ref{fig:demandProfile}. 
\begin{figure}[t]
    \centering

    \begin{tikzpicture}[scale = 0.7]\small
        \draw[-stealth] (-1,0) -- (8,0);
        \draw(8,-0.4) node{$\tau$};
        \draw[-stealth] (-1,0) -- (-1,9);
        \draw(-1.4,9) node{$C_\tau$};
        \draw (0, -0.15) -- (0, 0.15);
        \draw (1, -0.15) -- (1, 0.15);
        \draw (2, -0.15) -- (2, 0.15);
        \draw (3, -0.15) -- (3, 0.15);
        \draw (4, -0.15) -- (4, 0.15);
        \draw (5, -0.15) -- (5, 0.15);
        \draw (6, -0.15) -- (6, 0.15);
        \draw (7, -0.15) -- (7, 0.15);
        \draw(0,-0.4) node{\scriptsize 0};
        \draw(1,-0.4) node{\scriptsize1};
        \draw(2,-0.4) node{\scriptsize2};
        \draw(3,-0.4) node{\scriptsize3};
        \draw(4,-0.4) node{\scriptsize4};
        \draw(5,-0.4) node{\scriptsize5};
        \draw(6,-0.4) node{\scriptsize6};
        \draw(7,-0.4) node{\scriptsize7};

        \draw[<->] (0,1.5) -- (0,0);
        \draw (0.65,0.75) node{\tiny $w(f_0)\lambda$};
        \draw[<->] (1,1.5) -- (1,3);
        \draw (1.65,2.25) node{\tiny$w(f_1)\lambda$};
        \draw[<->] (2,3) -- (2,4.5);
        \draw (2.65,3.75) node{\tiny$w(f_2)\lambda$};
        \draw[<->] (4,4.5) -- (4,8.5);
        \draw (4.65,6.5) node{\tiny $w(f_4)\lambda$};
        \draw[dotted] (0,1.5) -- (1,1.5);
        \draw[dotted] (1,3) -- (2,3);
        \draw[dotted] (3,4.5) -- (4,4.5);

        \filldraw[red] (-1,0) circle (2pt) node[anchor=west]{};
        \filldraw[red] (0,1.5) circle (2pt) node[anchor=west]{};
        \filldraw[red] (1,3) circle (2pt) node[anchor=west]{};
        \filldraw[red] (2,4.5) circle (2pt) node[anchor=west]{};
        \filldraw[red] (3,4.5) circle (2pt) node[anchor=west]{};        
        \filldraw[red] (4,8.5) circle (2pt) node[anchor=west]{};
        \filldraw[red] (5,8.5) circle (2pt) node[anchor=west]{};
        \filldraw[red] (6,8.5) circle (2pt) node[anchor=west]{};
        \filldraw[red] (7,8.5) circle (2pt) node[anchor=west]{};
        \draw[dashed, thick, red] (-1,0) -- (0,1.5);
        \draw[dashed, thick, red] (0,1.5) -- (1,3);
        \draw[dashed, thick, red] (1,3) -- (2,4.5);
        \draw[dashed, thick, red] (2,4.5) -- (3,4.5);
        \draw[dashed, thick, red] (3,4.5) -- (4,8.5);
        \draw[dashed, thick, red] (4,8.5) -- (5,8.5);
        \draw[dashed, thick, red] (5,8.5) -- (6,8.5);        
        \draw[dashed, thick, red] (6,8.5) -- (7,8.5);

        \filldraw[red] (3.5,1.5) circle (2pt); 
        \draw (5.8,1.5) node {\scriptsize \textcolor{black}{$\pi=(\bar{\tau},f_0, f_1, \dots, f_{\bar{\tau}})$}}; 
        \draw[gray] (3.2,1.9) rectangle (7.9,1.1);
    \end{tikzpicture}
    \caption{Example express shipment demand profile (dashed line for guidance only).}
    \label{fig:demandProfile}
\end{figure}
We note that $\mathbf{C}^\pi$ is increasing between $\tau=-1$ and $\tau=2$, which indicates that in expectation, there are customers at ages $\tau=0$, $\tau=1$, and $\tau=2$ who choose express shipment. Together with (\ref{eq:willingnessToPay}), it follows that $f_0,f_1,f_2<\umax$. In particular, $f_0 = f_1 = f_2$ as the increments of the demand profile at ages $\tau=0$, $\tau=1$, and $\tau=2$ are equal. Comparing these increments with the one at $\tau=4$, we find that the expected number of express shipment orders at $\tau=4$ is higher than at earlier ages and thus $f_4 < f_0 = f_1 = f_2$. 
At age $\tau=3$, no customer chooses express shipment, which indicates that $f_3 \geq \umax$. Similarly, either $f_\tau\geq \umax$ for any age $\tau\ge 5$, or alternatively, the cutoff point is $\bar{\tau}=4$ such that express shipment is not available to customers at any age $\tau\ge 5$.

This notion of express shipment demand profiles will be instrumental in analyzing how shipment policies affect expected backorders and penalty costs. In particular, we will apply a problem-domain transformation to work directly in the space of express shipment demand profiles. Comparing shipment policies in terms of these demand profiles rather than express shipment fees allows for insightful analytical results on their effects on expected profit. We first start with some basic results on the simplest class of shipment policies.

\subsection{Constant shipment policy}
\label{sec:constantShipmentPolicy}
As a baseline, we consider the simplest class of shipment policies that offer express shipment at every age $\tau$ of the operating cycle at a constant fee:
\begin{definition}[Constant shipment policy]\label{def:constantShipmentPolicy}
    Let $f\in [\umin,\umax]$. We denote by $\pi^c := (T-1, f, \hdots, f)$ a \emph{\underline{c}onstant shipment policy}. The set of constant shipment policies is denoted by $\Pi^{c}$. 
    \hfill \Halmos \endproof
\end{definition}

We know by Lemma \ref{lem:revenueMaximizingPolicy} that the constant shipment policy $\pi^c = (T-1, f^r, \hdots, f^r)$ with $f^r=\umax/2$ maximizes the expected revenue. However, as seeking to maximize the expected profit, we are also interested in the penalty costs and expected backorders of the constant shipment policy:   

\begin{lemma}[Constant shipment policy]\label{lem:constantShipmentPolicy}
\phantom{texttextext}
    \begin{enumerate}
        \item[(1)] The expected backorders $\E[B^{\pi^c}]$ are convex in the express shipment fee $f$ with unique backorder-minimizing fee $f^b$. 
        \item[(2)] There is a unique optimal constant shipment policy $\pi^{c*}=(T-1, f^*, \hdots, f^*)$ with $f^* \in[f^r,f^b]$.
    \end{enumerate}
\end{lemma}

Recalling the expected profit $\E[G] = \E[R]  - \beta \E[B]$ from (\ref{eq:ProfitVariable}) and that $\E[R]$ is concave in $f$ by Lemma~\ref{lem:revenueMaximizingPolicy}, $\E[G^{\pi^c}]$ is concave in $f$ and thus $\pi^{c*}$ is unique. 
At the optimal express shipment fee $f^*$, the marginal revenue and the marginal penalty cost of the constant shipment policy are equal. 
As the penalty cost is non-increasing, the marginal penalty cost is increasing, and the marginal revenue is decreasing in the express shipment fee, $f^*$ is at most as high as the revenue-maximizing fee $f^r$ and at least as high as the backorder-minimizing fee $f^b$. 

Furthermore, $f^*$ increases in the penalty cost rate $\beta$ since a profit-maximizing e-commerce company may want to reduce the expected backorders when facing a higher penalty cost rate. To do so, the e-commerce company charges a higher express shipment fee, which incentivizes more customers to choose regular shipment instead of express shipment and the total expected express shipment orders per cycle decreases, resulting in fewer expected backorders.

\subsection{Constant cutoff point shipment policy}
\label{sec:constantShipmentPolicyCutoffPoint}
A cutoff point for express shipment orders is widely used in practice to manage the inflow of express shipment orders. That is, express shipment is only available to customers who order until the cutoff point $\bar{\tau}$, i.e., at times $t$ when $\tau = t \Mod T \leq \bar{\tau}$.
\begin{definition}[Constant cutoff point shipment policy]\label{def:constantShipmentPolicyCutoff}
    Let $\bar{\tau}\in\{0,1,\dots, T-1\}$ and $f\in [\umin,\umax]$. We denote by $\pi^{cc} := (\bar{\tau}, f, \dots, f)$ a \emph{\underline{c}onstant \underline{c}utoff point shipment policy}. The set of constant cutoff point shipment policies is denoted by $\Pi^{cc}$. 
    \hfill $\Box$
    \end{definition}
Note that any constant shipment policy $\pi^c\in\Pi^c$ as introduced in Definition~\ref{def:constantShipmentPolicy} is a constant cutoff point shipment policy with $\bar{\tau} = T-1$, i.e., $\Pi^c \subset \Pi^{cc}$. 

Introducing a cutoff point $\bar{\tau} < T-1$ has several merits from an operational perspective as it guarantees a minimum time buffer between the arrival and shipment date of an express shipment order. This buffer provides time to execute operations effectively. Moreover, the cutoff point is easy to communicate and intuitive to understand for customers and therefore marketable in practice. For example Coolblue, a large European electronics e-tailer, prominently features ``ordered before 23:59, delivered tomorrow'' on their website.

To characterize the benefits of a cutoff point for express shipment, we exploit the problem-domain transformation to the space of express shipment demand profiles. In particular, we compare the express shipment demand profiles of a constant cutoff point shipment policy $\pi^{cc}\in\Pi^{cc}\setminus \Pi^c$ and a constant shipment policy $\pi^c\in\Pi^c$ that yield the same total express shipment orders per cycle, i.e., $C_{T-1}^{\pi^{cc}}=C_{T-1}^{\pi^c}$, illustrated in Figure~\ref{fig:demandProfileCutoffConstantPolicy}.
\begin{figure}[t]
    \centering

    \begin{tikzpicture}[scale = 0.7]\small
        \draw[-stealth] (-1,0) -- (8,0);
        \draw(8,-0.4) node{$\tau$};
        \draw[-stealth] (-1,0) -- (-1,9);
        \draw(-1.4,9) node{$C_\tau$};
        \draw (0, -0.15) -- (0, 0.15);
        \draw (1, -0.15) -- (1, 0.15);
        \draw (2, -0.15) -- (2, 0.15);
        \draw (3, -0.15) -- (3, 0.15);
        \draw (4, -0.15) -- (4, 0.15);
        \draw (5, -0.15) -- (5, 0.15);
        \draw (6, -0.15) -- (6, 0.15);
        \draw (7, -0.15) -- (7, 0.15);
        \draw(0,-0.4) node{\scriptsize0};
        \draw(1,-0.4) node{\scriptsize1};
        \draw(2,-0.4) node{\scriptsize2};
        \draw(3,-0.4) node{\scriptsize3};
        \draw(4,-0.4) node{\scriptsize4};
        \draw(5,-0.4) node{\scriptsize5}; \draw(5, -0.85) node[red]{$\bar{\tau}$};
        \draw(6,-0.4) node{\scriptsize6};
        \draw(7,-0.4) node{\scriptsize7};

        \draw[<->] (5,8) -- (5,20/3);
        \draw (5.65,7.45) node{\tiny$w(f')\lambda$};
        \draw[dotted] (4,20/3) -- (5,20/3);    
        \draw[<->] (5,6) -- (5,5);
        \draw (5.65,5.5) node{\tiny$w(f)\lambda$};
        \draw[dotted] (4,5) -- (5,5);   

        \filldraw[blue] (-1.1cm,0) -- (-0.9cm,0) -- (-1cm,0.2cm);
        \filldraw[blue] (-0.1cm,1) -- (0.1cm,1) -- (0,1.2cm);
        \filldraw[blue] (0.9cm,2) -- (1.1cm,2) -- (1cm,2.2cm);
        \filldraw[blue] (1.9cm,3) -- (2.1cm,3) -- (2cm,3.2cm);
        \filldraw[blue] (2.9cm,4) -- (3.1cm,4) -- (3cm,4.2cm);
        \filldraw[blue] (3.9cm,5) -- (4.1cm,5) -- (4cm,5.2cm);
        \filldraw[blue] (4.9cm,6) -- (5.1cm,6) -- (5cm,6.2cm);
        \filldraw[blue] (5.9cm,7) -- (6.1cm,7) -- (6cm,7.2cm);
        \filldraw[blue] (6.9cm,8) -- (7.1cm,8) -- (7cm,8.2cm);
        \draw[dashed, thick, blue] (-1,0) -- (7,8);
        
        \filldraw[red] (-1,0) circle (2.5pt) node[anchor=west]{};
        \filldraw[red] (0,4/3) circle (2.5pt) node[anchor=west]{};
        \filldraw[red] (1,8/3) circle (2.5pt) node[anchor=west]{};
        \filldraw[red] (2,12/3) circle (2.5pt) node[anchor=west]{};
        \filldraw[red] (3,16/3) circle (2.5pt) node[anchor=west]{};
        \filldraw[red] (4,20/3) circle (2.5pt) node[anchor=west]{};
        \filldraw[red] (5,24/3) circle (2.5pt) node[anchor=west]{};        
        \filldraw[red] (6,8) circle (2.5pt) node[anchor=west]{};        
        \filldraw[red] (7,8) circle (2.5pt) node[anchor=west]{};
        \draw[dashed, thick, red] (-1,0) -- (5,8);
        \draw[dashed, thick, red] (5,8) -- (7,8);

        \filldraw[red] (3.5,2) circle (2.5pt); 
        \draw (5.6,2) node {\scriptsize \textcolor{black}{$\pi^{cc}=(\bar{\tau}, f',\dots, f')$}}; 
        \filldraw[blue] (3.4cm,1.4) -- (3.6cm,1.4) -- (3.5cm,1.6cm); 
        \draw (5.8,1.5) node {\scriptsize \textcolor{black}{$\pi^c=(T-1, f,\dots, f)$}};   
        \draw[gray] (3.2,2.35) rectangle (7.9,1.2);
        
    \end{tikzpicture}
    \caption{Express shipment demand profile of constant cutoff point shipment policy (dashed line for guidance only).} 
    \label{fig:demandProfileCutoffConstantPolicy}
\end{figure}
Note that $\pi^c$ and $\pi^{cc}$ distribute the same total express shipment orders differently within the operating cycle. While $\pi^c$ distributes the expected express shipment orders equally over all time periods, $\pi^{cc}$ distributes them equally over the time periods until the cutoff point. That is, $\pi^{cc}$ has a higher expected fraction of express shipment orders at early ages of the operating cycle than $\pi^{c}$. 
Having a higher fraction of express shipment orders at early ages ensures that the processing capacity in these time periods is better utilized to process express shipment orders. In particular in busy operating cycles, the scarce processing capacity is fully utilized for express shipment orders and none is wasted on orders for regular shipment, resulting in fewer expected backorders. Thus, $\pi^{cc}$ generates fewer expected backorders than $\pi^{c}$.

Noting that the demand increments of $\mathbf{C}^{\pi^{cc}}$ are larger than the ones of $\mathbf{C}^{\pi^{c}}$ and recalling that the expected express shipment orders are non-decreasing in the fee, $\pi^{cc}$ charges a lower express shipment fee than $\pi^{c}$, i.e., $f'<f$. 
Together with $\E[R^{\pi^{cc}}] = f'C^{\pi^{cc}}_{T-1}$ and $\E[R^{\pi^c}] = fC^{\pi^c}_{T-1}$, 
it follows that $\pi^{cc}$ yields a lower expected revenue than $\pi^{c}$. Consequently, $\pi^{cc}$ achieves a higher expected profit than $\pi^{c}$ only if the penalty cost rate is sufficiently high such that the backorder cost savings of $\pi^{cc}$ countervails its revenue loss compared to $\pi^{c}$. These observations are formalized as follows:
\begin{lemma}[Effect of cutoff point]
    Let $\pi^c = (T-1, f,\dots, f)\in \Pi^c$ and $\pi^{cc} = (\bar{\tau}, f',\dots, f')\in\Pi^{cc}\setminus\Pi^c$. Suppose that $C_{T-1}^{\pi^c} = C_{T-1}^{\pi^{cc}}$. 
    \begin{enumerate}
         \item[(1)] $\pi^{cc}$ charges a smaller express shipment fee than $\pi^{c}$, i.e., $f' < f$.       
        \item[(2)] $\pi^{cc}$ generates fewer expected backorders than $\pi^{c}$, i.e., $\E[B^{\pi^{cc}}] \leq \E[B^{\pi^{c}}]$.
        \item[(3)] $\pi^{cc}$ yields a higher expected profit than $\pi^{c}$, i.e., $\E[G^{\pi^{cc}}] \geq \E[G^{\pi^{c}}]$, if and only if 
        \begin{equation}
            \beta  \geq (f-f') \cdot \frac{C_{T-1}^{\pi^{cc}}}{\E[B^{\pi^{c}}] -\E[B^{\pi^{cc}}]}. 
        \end{equation}
    \end{enumerate}
    \label{lem:constantShipmentPolicyCutoffPoint}
\end{lemma}

This lemma provides e-commerce companies that operate a constant shipment policy at their fulfillment centers with a clearly actionable insight: They can reduce the expected backorders by introducing a cutoff point for express shipment and reducing the express shipment fee. This will yield a higher expected profit if the penalty cost rate is sufficiently high. 

\subsection{Two-level time-dependent shipment policy}
\label{sec:twoLevelTimeDependentShipmentPolicy}
We propose the two-level time-dependent shipment policy that differentiates two express shipment fees and switches between these once per operating cycle; see Figure~\ref{fig:ServicePhases} for an illustration. This policy has certain practical appeal as an adjustment of the express shipment fee allows for a more differentiated and targeted control of the inflow of express shipment orders than a cutoff point. Moreover, a shipment policy consisting of two fee values and one fee switching point is still intuitive for customers and thus marketable in practice. 
\begin{figure}[t]
    \centering

    \begin{tikzpicture}[scale = 0.9]\small
        \draw[thick, -stealth] (0.5,0) -- (18,0);
        \draw(18,-0.4) node{\scriptsize time};

        \draw (0.5, -0.15) -- (0.5, 0.15);
        \draw[thick] (2.25, -0.3) -- (2.25, 0.3);
        \draw (4.0, -0.15) -- (4.0, 0.15);
        \draw (5.75, -0.15) -- (5.75, 0.15);
        \draw[thick] (7.5, -0.3) -- (7.5, 0.3);
        \draw (9.25, -0.15) -- (9.25, 0.15);
        \draw (11.0, -0.15) -- (11.0, 0.15);
        \draw[thick] (12.75, -0.3) -- (12.75, 0.3);
        \draw (14.5, -0.15) -- (14.5, 0.15);
        \draw[thick] (16.25, -0.3) -- (16.25, 0.3);
        
        \draw (2.25, 0.7) -- (2.25, 0.7) node{$\tau = 0$};
        \draw (7.5, 0.7) -- (7.5, 0.7) node{$\tau = \hat{\tau}$};
        \draw (12.75, 0.7) -- (12.75, 0.7) node{$\tau = \bar{\tau}$};
        \draw (16.25, 0.7) -- (16.25, 0.7) node{$\tau = 0$};
        
        \draw (2.25, 1.2) -- (2.25, 1.2) node{\footnotesize Deadline};
        \draw (7.5, 1.2) -- (7.5, 1.2) node{\footnotesize Fee switching point};
        \draw (12.75, 1.2) -- (12.75, 1.2) node{\footnotesize Cutoff point};
        \draw (16.25, 1.2) -- (16.25, 1.2) node{\footnotesize Deadline};
        
        \draw[dashed] (2.25, -0.4) -- (2.25, -1.6);
        \draw[dashed] (7.5, -0.4) -- (7.5, -1.6);
        \draw[dashed] (12.75, -0.4) -- (12.75, -1.6);
        \draw[dashed] (16.25, -0.4) -- (16.25, -1.6);
        
        \draw[thick, -stealth] (3.5, -1.6) -- (7.5, -1.6);
        \draw[thick, stealth-] (2.25, -1.6) -- (3.5, -1.6);
        \draw (2.25, -1.5) -- (2.25, -1.7);
        \draw (7.5, -1.5) -- (7.5, -1.7);
        
        \draw[thick, -stealth] (8.5, -1.6) -- (12.75, -1.6);
        \draw[thick, stealth-] (7.5, -1.6) -- (8.5, -1.6);
        \draw (7.5, -1.5) -- (7.5, -1.7);
        \draw (12.75, -1.5) -- (12.75, -1.7);
        
        \draw[thick, -stealth] (13.75, -1.6) -- (16.25, -1.6);
        \draw[thick, stealth-] (12.75, -1.6) -- (13.75, -1.6);
        \draw (12.75, -1.5) -- (12.75, -1.7);
        \draw (16.25, -1.5) -- (16.25, -1.7);
        
        \draw(4.875, -2.2) -- (4.875, -2.2) node{\footnotesize Express shipment};
        \draw(4.875, -2.7) -- (4.875, -2.7) node{$f_\tau = f^e$}; 
        
       \draw(10, -2.2) -- (10, -2.2) node{\footnotesize Last-minute express shipment};
        \draw(10.125, -2.7) -- (10.125, -2.7) node{$f_\tau = f^l$}; 
        
        \draw(14.5, -2.2) -- (14.5, -2.2) node{\footnotesize no express shipment};
        
    \end{tikzpicture}
    
    \caption{Illustration of the two-level time-dependent shipment policy.}
    \label{fig:ServicePhases}
\end{figure}
\begin{definition}[Two-level time-dependent shipment policy] \label{def:twoLevelShipmentPolicy}
     Consider the following policy parameters:
     $\bar{\tau} \in\{0,1,\dots,T-1\}$ as the cutoff point; 
     $\hat{\tau}\in\mathbb{N}_0$, $\hat{\tau} < \bar{\tau}$, as the \emph{fee switching point};
     $f^e,f^l\in[\umin, \umax]$, $f^e\neq f^l$, as the \emph{express fee} and \emph{last-minute express fee}.
    The corresponding \emph{\underline{two}-level \underline{t}ime-dependent shipment policy} $\pi^{2t}=(\bar{\tau}, f_0, f_1, \dots, f_{\bar{\tau}})$ sets $f_\tau = f^e$ if $\tau\leq \hat{\tau}$ and $f_\tau = f^l$ if $\hat{\tau} < \tau \leq \bar{\tau}$. The set of two-level time-dependent shipment policies is denoted by $\Pi^{2t}$.   \hfill $\Box$ 
\end{definition}

The corresponding express shipment demand profile $\mathbf{C}^{\pi^{2t}}$, illustrated in Figure~\ref{fig:demandProfileTwoLevelPolicy}a, has a two-level structure similar to that of the express shipment fees. In particular, the increments of the demand profile at ages $\tau = 0,\dots,\hat{\tau}$ are constant and equal $w(f^e)\lambda$. At the fee switching point $\hat{\tau}$, the increment changes to $w(f^l)\lambda$ and all increments at ages $\tau=\hat{\tau}+1,\dots,\bar{\tau}$ equal $w(f^l)\lambda$, 

As longer time buffers between the arrival and shipment date of an express shipment order provide more flexibility to process this order on time and reduce the risk of backorders, an order fulfillment center prefers a high fraction of expected express shipment orders at early ages of the operating cycle that reduces to a lower fraction at the fee switching point, i.e., $w(f^e)\lambda > w(f^l)\lambda$. This is illustrated in Figure~\ref{fig:demandProfileTwoLevelPolicy}b and formalized as follows:
\begin{lemma}[Type of two-level fee] 
    \label{lem:twoLevelTimeDependentShipmentPolicy}
    Let $\pi^{2ti}= (\bar{\tau}, f', \dots, f', f'', \dots, f'')\in\Pi^{2t}$ and $\pi^{2td}= (\bar{\tau},  f'', \dots, f'', f', \dots, f')\in\Pi^{2t}$ where $f'< f''$ and $C_{T-1}^{\pi^{2ti}}=C_{T-1}^{\pi^{2td}}$. \\
    Then,  $\pi^{2ti}$ generates fewer expected backorders and achieves a higher expected profit than $\pi^{2td}$, i.e., 
    \begin{equation*}
         \E[B^{\pi^{2ti}}] \leq \E[B^{\pi^{2td}}] \quad\quad\quad\quad\quad\quad\quad\quad\quad\
         \E[G^{\pi^{2ti}}] \geq \E[G^{\pi^{2td}}].
    \end{equation*}
\end{lemma}
As a consequence, an order fulfillment center that operates a two-level time-dependent shipment policy $\pi^{2t}\in\Pi^{2t}$ increases the express shipment fee at the fee switching point and sets $f^e<f^l$. 

\begin{figure}
    \centering
	\subfloat[Effect of two-level time-dependent fee]{%
		\begin{tikzpicture}[scale = 0.75]\small
            \draw[-stealth] (-1,0) -- (8,0);
            \draw(8,-0.4) node{$\tau$};
            \draw[-stealth] (-1,0) -- (-1,9);
            \draw(-1.4,9) node{$C_\tau$};
            \draw (0, -0.15) -- (0, 0.15);
            \draw (1, -0.15) -- (1, 0.15);
            \draw (2, -0.15) -- (2, 0.15);
            \draw (3, -0.15) -- (3, 0.15);
            \draw (4, -0.15) -- (4, 0.15);
            \draw (5, -0.15) -- (5, 0.15);
            \draw (6, -0.15) -- (6, 0.15);
            \draw (7, -0.15) -- (7, 0.15);
            \draw(0,-0.4) node{\scriptsize0};
            \draw(1,-0.4) node{\scriptsize1}; \draw(1, -0.85) node[red]{$\hat{\tau}$};
            \draw(2,-0.4) node{\scriptsize2};
            \draw(3,-0.4) node{\scriptsize3};
            \draw(4,-0.4) node{\scriptsize4};
            \draw(5,-0.4) node{\scriptsize5}; \draw(5, -0.85) node[black]{$\bar{\tau}$};
            \draw(6,-0.4) node{\scriptsize6};
            \draw(7,-0.4) node{\scriptsize7};

            \filldraw[blue] (-1.1cm,0) -- (-0.9cm,0) -- (-1cm,0.2cm);
            \filldraw[blue] (-0.1cm,4/3) -- (0.1cm,4/3) -- (0,1.533333cm);
            \filldraw[blue] (0.9cm,8/3) -- (1.1cm,8/3) -- (1cm,2.86667cm);
            \filldraw[blue] (1.9cm,4) -- (2.1cm,4) -- (2cm,4.2cm);
            \filldraw[blue] (2.9cm,16/3) -- (3.1cm,16/3) -- (3cm,5.533333cm);
            \filldraw[blue] (3.9cm,20/3) -- (4.1cm,20/3) -- (4cm,6.866667cm);
            \filldraw[blue] (4.9cm,8) -- (5.1cm,8) -- (5cm,8.2cm);
            \filldraw[blue] (5.9cm,8) -- (6.1cm,8) -- (6cm,8.2cm);
            \filldraw[blue] (6.9cm,8) -- (7.1cm,8) -- (7cm,8.2cm);
            \draw[dashed, thick, blue] (-1,0) -- (5,8);
            \draw[dashed, thick, blue] (5,8) -- (7,8);
            
            \filldraw[red] (-1,0) circle (2.5pt) node[anchor=west]{};
            \filldraw[red] (0,2.5) circle (2.5pt) node[anchor=west]{};
            \filldraw[red] (1,5) circle (2.5pt) node[anchor=west]{};
            \filldraw[red] (2,5.75) circle (2.5pt) node[anchor=west]{};
            \filldraw[red] (3,6.5) circle (2.5pt) node[anchor=west]{};
            \filldraw[red] (4,7.25) circle (2.5pt) node[anchor=west]{};
            \filldraw[red] (5,8) circle (2.5pt) node[anchor=west]{};        
            \filldraw[red] (6,8) circle (2.5pt) node[anchor=west]{};        
            \filldraw[red] (7,8) circle (2.5pt) node[anchor=west]{};
            \draw[dashed, thick, red] (-1,0) -- (1,5);
            \draw[dashed, thick, red] (1,5) -- (5,8);
            \draw[dashed, thick, red] (5,8) -- (7,8);

            \filldraw[red] (2.35,1.3) circle (2.5pt); 
            \draw (5.1,1.3) node {\tiny \textcolor{black}{$\pi^{2t}= (\bar{\tau}, f^e, \dots, f^e, f^l, \dots, f^l)$}}; 
            \filldraw[blue] (2.25cm,0.6) -- (2.45cm,0.6) -- (2.35cm,0.8cm); 
            \draw (4.1,0.7) node {\tiny \textcolor{black}{$\pi^{cc} = (\bar{\tau}, f,\dots, f)$}};   
            \draw[gray] (2.0,1.65) rectangle (7.9,0.4);
        \end{tikzpicture}
	}
    \hfill
    	\subfloat[Type of two-level time-dependent fee; $f'<f''$]{%
		    \begin{tikzpicture}[scale = 0.75]\small
        \draw[-stealth] (-1,0) -- (8,0);
        \draw(8,-0.4) node{$\tau$};
        \draw[-stealth] (-1,0) -- (-1,9);
        \draw(-1.4,9) node{$C_\tau$};
        \draw (0, -0.15) -- (0, 0.15);
        \draw (1, -0.15) -- (1, 0.15);
        \draw (2, -0.15) -- (2, 0.15);
        \draw (3, -0.15) -- (3, 0.15);
        \draw (4, -0.15) -- (4, 0.15);
        \draw (5, -0.15) -- (5, 0.15);
        \draw (6, -0.15) -- (6, 0.15);
        \draw (7, -0.15) -- (7, 0.15);
        \draw(0,-0.4) node{\scriptsize0};
        \draw(1,-0.4) node{\scriptsize1}; \draw(1, -0.85) node[red]{$\hat{\tau}^i$};
        \draw(2,-0.4) node{\scriptsize2};
        \draw(3,-0.4) node{\scriptsize3}; \draw(3, -0.85) node[blue]{$\hat{\tau}^d$};
        \draw(4,-0.4) node{\scriptsize4};
        \draw(5,-0.4) node{\scriptsize5}; \draw(5, -0.85) node[black]{$\bar{\tau}$};
        \draw(6,-0.4) node{\scriptsize6};
        \draw(7,-0.4) node{\scriptsize7};

        \draw[<->] (1,2.5) -- (1,5);
        \draw (1.65,3.75) node{\tiny$w(f')\lambda$};
        \draw[dotted] (1,2.5) -- (0,2.5);   
        \draw[<->] (2,5.75) -- (2,5);
        \draw (2.7,5.375) node{\tiny$w(f'')\lambda$};
        \draw[dotted] (2,5) -- (1,5); 
        \draw[<->] (4,5.5) -- (4,3);
        \draw (4.65,4.25) node{\tiny$w(f')\lambda$};
        \draw[dotted] (3,3) -- (4,3);   
        \draw[<->] (3,3) -- (3,2.25);
        \draw (3.7,2.625) node{\tiny$w(f'')\lambda$};
        \draw[dotted] (2,2.25) -- (3,2.25);

         \filldraw[blue] (-1.1cm,0) -- (-0.9cm,0) -- (-1cm,0.2cm);
        \filldraw[blue] (-0.1cm,0.75) -- (0.1cm,0.75) -- (0,0.95cm);
        \filldraw[blue] (0.9cm,1.5) -- (1.1cm,1.5) -- (1cm,1.7cm);
        \filldraw[blue] (1.9cm,2.25) -- (2.1cm,2.25) -- (2cm,2.45cm);
        \filldraw[blue] (2.9cm,3) -- (3.1cm,3) -- (3cm,3.2cm);
        \filldraw[blue] (3.9cm,5.5) -- (4.1cm,5.5) -- (4cm,5.7cm);
        \filldraw[blue] (4.9cm,8) -- (5.1cm,8) -- (5cm,8.2cm);
        \filldraw[blue] (5.9cm,8) -- (6.1cm,8) -- (6cm,8.2cm);
        \filldraw[blue] (6.9cm,8) -- (7.1cm,8) -- (7cm,8.2cm);
        \draw[dashed, thick, blue] (-1,0) -- (3,3);
        \draw[dashed, thick, blue] (3,3) -- (5,8);
        \draw[dashed, thick, blue] (5,8) -- (7,8);
        
        \filldraw[red] (-1,0) circle (2.5pt) node[anchor=west]{};
        \filldraw[red] (0,2.5) circle (2.5pt) node[anchor=west]{};
        \filldraw[red] (1,5) circle (2.5pt) node[anchor=west]{};
        \filldraw[red] (2,5.75) circle (2.5pt) node[anchor=west]{};
        \filldraw[red] (3,6.5) circle (2.5pt) node[anchor=west]{};
        \filldraw[red] (4,7.25) circle (2.5pt) node[anchor=west]{};
        \filldraw[red] (5,8) circle (2.5pt) node[anchor=west]{};        
        \filldraw[red] (6,8) circle (2.5pt) node[anchor=west]{};        
        \filldraw[red] (7,8) circle (2.5pt) node[anchor=west]{};
        \draw[dashed, thick, red] (-1,0) -- (1,5);
        \draw[dashed, thick, red] (1,5) -- (5,8);
        \draw[dashed, thick, red] (5,8) -- (7,8);

        \filldraw[red] (2.35,1.3) circle (2.5pt); 
        \draw (5.1,1.3) node {\tiny \textcolor{black}{$\pi^{2ti}= (\bar{\tau}, f', \dots, f', f'', \dots, f'')$}}; 
        \filldraw[blue] (2.25cm,0.6) -- (2.45cm,0.6) -- (2.35cm,0.8cm); 
        \draw (5.1,0.7) node {\tiny \textcolor{black}{$\pi^{2td}= (\bar{\tau},  f'', \dots, f'', f', \dots, f')$}};   
        \draw[gray] (2.0,1.65) rectangle (7.9,0.4);
        
    \end{tikzpicture}
	}
    \caption{Express shipment demand profiles of two-level time-dependent shipment policies (dashed line for guidance only).}
    \label{fig:demandProfileTwoLevelPolicy}
\end{figure}
To characterize the benefit of a two-level time-dependent fee for express shipment, we exploit the problem-domain transformation to the space of express shipment demand profiles. We compare the express shipment demand profiles of a two-level time-dependent shipment policy $\pi^{2t}=(\bar{\tau}, f^e, \dots, f^e, f^l, \dots, f^l)\in\Pi^{2t}$ and a constant cutoff point shipment policy $\pi^{cc}=(\bar{\tau}, f,\dots, f)\in\Pi^{cc}$ that have the same cutoff point and the same total expected express shipment orders per cycle, i.e., $C_{T-1}^{\pi^{2t}}=C_{T-1}^{\pi^{cc}}$, illustrated in Figure~\ref{fig:demandProfileTwoLevelPolicy}a. Note that $\mathbf{C}^{\pi^{2t}}$ has larger demand increments than $\mathbf{C}^{\pi^{cc}}$ until the fee switching point, and smaller ones thereafter. Thus, $\pi^{2t}$ has a higher expected fraction of express shipment orders than $\pi^{cc}$ at early ages of the operating cycle. This ensures that $\pi^{2t}$ can better utilize the processing capacity at early ages to process express shipment orders and results in fewer expected backorders than $\pi^{cc}$. 

Furthermore, these differences in the express shipment demand profiles imply that $\pi^{2t}$ charges a lower shipment fee than $\pi^{cc}$ until the fee switching point, but a higher one thereafter, i.e., $f^e < f < f^l$. As the expected revenue is concave in the shipment fee (Lemma~\ref{lem:revenueMaximizingPolicy}), $\pi^{2t}$ yields a lower expected revenue than $\pi^{cc}$. Consequently, $\pi^{2t}$ achieves a higher expected profit than $\pi^{cc}$ if the penalty cost rate is sufficiently high. 
These observations are formalized as follows:
\begin{lemma}[Effect of two-level time-dependent fee] \label{lem:twoLevelTimeDependentShipmentPolicyComparedConstant}
    Let $\pi^{cc} = (\bar{\tau}, f,\dots, f)\in\Pi^{cc}$ and $\pi^{2t}= (\bar{\tau}, f^e, \dots, f^e, f^l, \dots, f^l)\in\Pi^{2t}$ where $f^e< f^l$ and $C_{T-1}^{\pi^{cc}}=C_{T-1}^{\pi^{2t}}$. 
    \begin{enumerate}
        \item[(1)] The express shipment fees satisfy $f = \epsilon f^e + \left(1 - \epsilon\right) f^l$ with $\epsilon = (\hat{\tau}+1)/(\bar{\tau}+1)$.
        \item[(2)]  $\pi^{2t}$ generates fewer expected backorders than $\pi^{cc}$, i.e., $    \E[B^{\pi^{2t}}] \leq \E[B^{\pi^{cc}}]$.
        \item[(3)] $\pi^{2t}$ yields a higher expected profit than $\pi^{cc}$, i.e., $\E[G^{\pi^{2t}}] \geq \E[G^{\pi^{cc}}]$, if and only if
        \begin{equation}
            \beta \geq \frac{\E[R^{\pi^{cc}}]-\E[R^{\pi^{2t}}]}{\E[B^{\pi^{cc}}]-\E[B^{\pi^{2t}}]}.
        \end{equation}
    \end{enumerate}
\end{lemma}

The key managerial insight from this lemma is that an e-fulfillment center that operates a constant cutoff point shipment policy can reduce the expected backorders by charging a two-level time-dependent fee for express shipment that increases at the fee switching point. This will yield a higher expected profit if the penalty cost rate is sufficiently high.

\subsection{General time-dependent shipment policy}
\label{sec:timeDependentShipmentPolicy}
A natural next step is to analyze general shipment policies $\pi\in\Pi$. Recall from Definition~\ref{def:setShipmentPolicies} that such policies $\pi=(\bar{\tau},f_0, f_1, \hdots, f_{\bar{\tau}})\in\Pi$ may charge different express shipment fees $f_\tau$ at every age $\tau$ of the operating cycle. Such policies can potentially outperform the restricted policy classes $\Pi^{c}$, $\Pi^{cc}$, $\Pi^{2t}$ analyzed in \S\ref{sec:constantShipmentPolicy}-\ref{sec:twoLevelTimeDependentShipmentPolicy}, but are less intuitive for customers. In this section, we establish structural properties of general time-dependent shipment policies that improve our understanding and facilitate the efficient computation of the optimal shipment policy in \S\ref{sec:optimizingshipmentfees}.

Recall that we seek to maximize the expected profit $\E[G]=\E[R]- \beta \E[B]$, and that we established in Lemma~\ref{lem:revenueMaximizingPolicy} that the constant shipment policy with $f = f^r$ maximizes the expected revenue $\E[R]$. Thus, any deviation from constant pricing hurts the expected revenue. On the other hand, we know from Lemma \ref{lem:twoLevelTimeDependentShipmentPolicyComparedConstant} that deviations from constant pricing may be beneficial, as they can reduce expected backorders $\E[B]$, and therefore penalty costs $ \beta \E[B]$. It turns out that the scope of this insight can be broadened towards general time-dependent shipment policies:
\begin{lemma}[Effect of time-dependent fee] \label{lem:timeDependentShipmentPolicy}
    Let $\pi^{cc} = (\bar{\tau}, \bar{f},\dots, \bar{f})\in\Pi^{cc}$, $\pi = (\bar{\tau}, f_0, f_1, \dots, f_{\bar{\tau}})\in\Pi$, and $\pi' = (\bar{\tau}, f'_0, f'_1, \hdots, f'_{\bar{\tau}})\in\Pi$, and suppose all three shipment policies have the same average express shipment fee ($\bar{f}=\frac{1}{\bar{\tau}+1}\sum_{j=0}^{\bar{\tau}} f_j =\frac{1}{\bar{\tau}+1} \sum_{j=0}^{\bar{\tau}} f'_j$). 
    Moreover, for each $\tau\in\{1,\ldots,\bar{\tau}\}$, the \emph{early-cycle average express shipment fee} is lowest under $\pi$, higher under $\pi'$, and highest under $\pi^{cc}$:
    \begin{equation}
    \frac{1}{\tau} \sum_{j=0}^{\tau-1} f_j \leq \frac{1}{\tau} \sum_{j=0}^{\tau-1} f'_j \le \bar{f}.
    \end{equation}
Then, $\E[B^{\pi}] \le \E[B^{\pi'}] \le \E[B^{\pi^{cc}}]$, i.e., $\pi$ has the fewest expected backorders and $\pi^{cc}$ the most.
\end{lemma}
This lemma generalizes the insights from earlier sections by demonstrating that also when fixing the average express shipment fee of a \emph{general} time-dependent shipment policy, lowering early-cycle average express shipment fees (and consequently increasing express shipment fees later in the cycle) tends to reduce the expected backorders. As a consequence of this insight, we furthermore expect that an optimal time-dependent shipment policy would be monotonic, and this can indeed be established formally:
\begin{lemma}[Monotonicity]\label{lem:monotonicityFee}
There exists a profit-maximizing shipment policy $\pi^*=(\bar{\tau},f_0,f_1,\ldots, f_{\bar{\tau}})\in\Pi$ with monotonous increasing fees, i.e., $f_\tau \geq f_{\tau-1}$ $\forall \tau \in\{1, 2,\dots,\bar{\tau}\}$. 
\end{lemma}

Our analytic results yield a clear understanding of the key trade-offs in designing time-dependent shipment policies. On the one hand, keeping the express shipment fees relatively stable throughout the operating cycle helps to preserve the expected revenue, as shown in Lemma \ref{lem:revenueMaximizingPolicy}. On the other hand, raising late-cycle express shipment fees (and decreasing early-cycle express shipment fees) can reduce expected backorders and thus penalty costs, as highlighted in Lemma \ref{lem:timeDependentShipmentPolicy}. We also know from Lemma \ref{lem:monotonicityFee} that there exists an optimal shipment policy with non-decreasing express shipment fees over time, which restricts the set of candidates to consider. However, even with this structural simplification, the search for the optimal general time-dependent shipment policy remains challenging: the number of shipment policies with monotone express shipment fees grows exponentially in the length of the operating cycle $T$. This naturally raises the question: how can we best balance these opposing forces to identify shipment policies that maximize the expected profit? 

\subsection{Optimizing shipment fees}
\label{sec:optimizingshipmentfees}
In this section, we discuss the identification of optimal policy parameters for any  restricted policy class $\Pi^c,\Pi^{cc},\Pi^{2t}$, as well as the identification of the profit-maximizing shipment policy $\pi^* \in \Pi$. 

For this purpose, it is natural to restrict attention to a discrete set of possible express shipment fees. Real-world pricing decisions rarely involve arbitrary continuous values; instead, companies typically choose from a limited set of price points, for example, in multiples of \$0.10 or \$0.20, to keep pricing simple for customers. To this end, let $\delta_f > 0$ denote the price increment between allowable express shipment fees; let $
n_f = \frac{\umax - \umin}{\delta_f}$, and suppose that $n_f\in \mathbb{N}$. Then, the set of discrete express shipment fees is
\begin{equation}
\mathcal{F} := \{\umin + m\,\delta_f \mid m = 1, 2, \ldots, n_f\}.
\end{equation}

Now, we will restrict attention to the following set of shipment policies: 
\begin{equation}
    \Pi^{\mathcal{F}} = \left\{ (\bar{\tau},f_0, f_1, \dots, f_{\bar{\tau}}) \mid  \bar{\tau} = T-1, \forall \tau = 0,1,\dots,T-1:f_\tau \in \mathcal{F}\right\}.
\end{equation}
The restriction to $\bar{\tau} = T-1$ does not entail a loss of generality: Any shipment policy with an earlier cutoff point $\bar{\tau} < T-1$ can be replicated by a shipment policy with $\bar{\tau} = T-1$ by setting $f_\tau= \umax$ for $\tau > \bar{\tau}$, so that no customer chooses express shipment in these periods (as illustrated and discussed in Figure~\ref{fig:demandProfile}).

We note that for any shipment policy $\pi$, we may determine $\E[G^\pi]$ by solving for the steady-state of a Markov chain. To optimize pricing, one can adopt enumeration techniques; for example, for identifying the best discrete policy in $\Pi^{c}$, we may solve the Markov chain for all discrete constant shipment policies and select the best one. The best discrete policies in $\Pi^{cc}$ and $\Pi^{2t}$ can be identified similarly. However, such enumeration techniques are no longer viable when solving for the optimal general time-dependent shipment policy $\pi \in \Pi^\mathcal{F}$, since the number of feasible policies $|\Pi^{\mathcal{F}}|$ grows exponentially with the length of the operating cycle $T$, even while accounting for the monotonicity guarantee of Lemma~\ref{lem:monotonicityFee}. Thus, to optimize general time-dependent shipment policies, it will be essential to gain further insight into how the penalty cost and profit functions behave as we vary the shipment fees.

It turns out that the problem-domain transformation to the space of express shipment demand profiles is again instrumental to achieving this insight. When shipment policies are compared in terms of these demand profiles, their effect on expected backorders and expected profit exhibits a useful lattice structural property:
\begin{proposition}[Supermodularity]\label{prop:submodularityShipmentPolicy}
Let $\pi,\pi'\in\Pi$ with associated express shipment demand profiles $\mathbf{C}^\pi$ and $\mathbf{C}^{\pi'}$.
\begin{enumerate}
    \item[(1)] There exist unique shipment policies $\pi^{min}$ and $\pi^{max}$ such that $\mathbf{C}^{\pi^{min}} =\min(\mathbf{C}^{\pi},\mathbf{C}^{\pi'})$ and $\mathbf{C}^{\pi^{max}}=\max(\mathbf{C}^{\pi},\mathbf{C}^{\pi'})$, where $\min$ and $\max$ are taken componentwise. 
    \item[(2)]  The expected backorders are submodular in $\mathbf{C}$, i.e., $\E[B^{\pi^{min}}]+\E[B^{\pi^{max}}]\le \E[B^{\pi}]+\E[B^{\pi'}]$.
    \item[(3)] The expected profit is supermodular in $\mathbf{C}$, i.e., $\E[G^{\pi^{min}}]+\E[G^{\pi^{max}}]\ge \E[G^{\pi}]+\E[G^{\pi'}]$. 
\end{enumerate}
\end{proposition}
The submodularity of expected backorders implies diminishing returns when reallocating express shipment demand earlier in the operating cycle: once early-cycle capacity is largely utilized, additional shifts produce progressively smaller reductions in expected backorders. Beyond this useful managerial perspective, the structural insight also provides the theoretical foundation needed to guide the search for high-quality time-dependent shipment policies.

To proceed, we explore the express shipment demand profiles $\mathbf{C}$ that arise for shipment policies $\pi\in\Pi^\mathcal{F}$. 
It is easy to see that $C^{\pi}_\tau\in \mathcal{C}_\tau$ for such policies, where $\mathcal{C}_\tau=\{m\lambda/n_f  | m=0,1,\ldots,(\tau+1) n_f \}$. Recalling $C^\pi_\tau-C^\pi_{\tau-1}=w(f_\tau)\lambda$ from (\ref{eq:demandProfileShipmentPolicy}), the set of demand profiles that may arise for shipment policies $\pi\in\Pi^{\mathcal{F}}$ becomes:
\[
\mathcal{C} := \left\{ (C_0,\dots,C_{T-1}) \ \middle|\ C_\tau \in \mathcal{C}_\tau,  0\le  C_\tau - C_{\tau-1} \le \lambda \ \forall \tau \ge 1 \right\}.
\]
Abusing notation, let $\E[G(\mathbf{C})]$ denote the expected profit in our transformed problem domain. Then, $\max_{\pi\in\Pi^{\mathcal{F}}} \E[G^\pi] = \max_{\mathbf{C}\in\mathcal{C}}\E[G(\mathbf{C})]$, i.e., identifying the profit-maximizing shipment policy is equivalent to finding the express shipment demand profile that maximizes the expected profit. By Proposition~\ref{prop:submodularityShipmentPolicy}, this casts our problem into the domain of supermodular function maximization / submodular function minimization (SFM) \citep{mccormick2005submodular}. Building on this connection, it turns out that we can reformulate the search for a profit-maximizing shipment policy as the maximization of a supermodular function over a ring family defined on a suitable defined finite ground set, which places our problem within a class for which polynomial-time algorithms are known to exist:

\newcommand{\EO}{\mathsf{EO}}

\begin{proposition}[Efficient identification of optimal policy]
\label{prop:polytime}
Let $\Pi^{\mathcal{F}}$ denote the set of discrete shipment policies. Given an algorithm that, for any $\pi\in\Pi^{\mathcal{F}}$, evaluates $\E[G^\pi]$ in time $\EO$ (via the steady-state of the associated Markov chain), there exists an algorithm that returns
\[
\pi^\star \in \arg\max_{\pi\in\Pi^{\mathcal{F}}} \E[G^\pi]
\]
in time $\mathrm{poly}(|\mathcal{U}|)\cdot \EO$, where $|\mathcal{U}|=\frac{n_f\,T(T+1)}{2}$ and $\EO$ is the time required for evaluating a single policy. In particular, by \citet[§3.3.2]{mccormick2005submodular}, a strongly polynomial oracle-time bound of $O\!\big(N^{7}\log N\cdot \EO\big)$ is achievable.
\end{proposition}
The polynomial-time algorithms underlying Proposition~\ref{prop:polytime}, while theoretically appealing, are often regarded as impractical for large problem instances. Motivated by this, we solve the supermodular maximization problem via a column generation framework as suggested in \citet[§4.1]{mccormick2005submodular}. At each iteration, a Carathéodory decomposition routine is applied to the current fractional solution in the master problem, identifying the integral demand profiles that support it. We then evaluate the shipment policies corresponding to these extreme points and add the associated inequality to the master problem. This process is repeated until no improving demand profiles are found, at which point optimality is guaranteed.

\section{Numerical analysis}
\label{sec:numericalAnalysis}
In this section, we demonstrate the performance of the two-level time-dependent shipment policy in a numerical study. In \S\ref{sec:AnalysisDoEImplementation}, we present the experimental design and discuss the model implementation. In \S\ref{sec:AnalysisBenefit}, we compare the optimal two-level time-dependent shipment policy with the optimal general time-dependent shipment policy as well as several benchmark policies prevalent in practice. In \S\ref{sec:AnalysisOptDesign}, we discuss the optimal policy parameters of the two-level time-dependent shipment policy. 

\subsection{Experimental design and model implementation}
\label{sec:AnalysisDoEImplementation}
We consider e-commerce fulfillment centers with an operating cycle of length $T = 8$ and a utilization of $\rho \in \{0.85, 0.9, 0.95\}$. The customer demand $D_t$, $t \in \mathbb{N}_0$, is $D_t\sim \poi(\lambda)$ with $\lambda = 5$. The processing capacity $K_t$, $t \in \mathbb{N}_0$, is $K_t\sim K$ where $K$ follows a discretized Beta-distribution with support $\{0, 1, \hdots, 20\}$, which is specified by the expected value $\E[K]$ and squared coefficient of variation $scv[K]$ \citep[295--297]{Law2007}. Note that $\E[K]=\lambda / \rho$ by (\ref{eq:U}) and $scv[K]=0.5$.
The product price is $\bar{p} = 40$ and the penalty cost rate is $\beta\in\{8,12\}$. The minimum and maximum express shipment valuations are $\umin=0$ and $\umax=0.1\bar{p}=4$: Note that the maximum express shipment fee is smaller than the penalty cost rate. We discretize shipment fees using $\delta_f=0.2$.

To implement and compute our model, we need to derive an upper bound for the state space of the Markov chain. Recalling from \S\ref{sec:ModelMarkovChain} that $(X_t^h, X_t^\Sigma, \tau)$ with $\tau = t \Mod T$ specifies the state of the Markov chain, the state space is $\mathbb{N}_0\times\mathbb{N}_0\times\{0,1,\dots,T-1\}$. As $X_t^h\leq X_t^\Sigma$ for every $t\in\mathbb{N}_0$, it is sufficient to derive an upper bound for $X_t^\Sigma$, denoted by $\bar{X}^\Sigma$. We determine $\bar{X}^\Sigma$ such that the probability of rejecting an order that arrives at any time $t$ is negligible small. In this study, the threshold value is $0.023$. We refer to Appendix~\ref{sec:AppendixImplementation} for further details on the construction of $\bar{X}^\Sigma$ and the resulting Markov chain model with finite state space. 

\subsection{Performance of two-level time-dependent shipment policy}
\label{sec:AnalysisBenefit}
We evaluate the performance of the two-level time-dependent shipment policy. We therefore compare the optimal two-level time-dependent shipment policy $\pi^{2t*}$ with the optimal general time-dependent shipment policy $\pi^*:=\argmax_{\pi\in\Pi} \E[G^\pi]$ and the following benchmark policies:
\begin{enumerate} 
   \item[(1)] profit-maximizing constant cutoff point shipment policy $\pi^{cc*} := \argmax_{\pi\in\Pi^{cc}} \E[G^\pi]$ where $\Pi^{cc} = \left\{\pi \in\Pi \mid \pi = (\bar{\tau}, f, \dots, f), \bar{\tau}\in\{0,1,\dots, 7\}, f\in  \mathcal{F} \right\}$;
    \item[(2)] profit-maximizing constant cutoff point shipment policy with revenue-maximizing fee $\pi^{ccr}:= \argmax_{\pi\in\Pi^{ccr}} \E[G^\pi]$ where $\Pi^{ccr} = \left\{\pi \in\Pi \mid \pi = (\bar{\tau}, f, \dots, f), \bar{\tau}\in\{0,1,\dots, 7\}, f=\umax/2 \right\}$;
     \item[(3)] revenue-maximizing shipment policy $\pi^{r}=(T-1,f^r, \dots, f^r)$ where $f^r=\umax/2$. 
\end{enumerate}

\begin{table}[h!]
    \centering
    \small
\begin{threeparttable}
    
    \caption{Optimal two-level time-dependent shipment policy $\pi^{2t*}$ compared to the optimal general time-dependent shipment policy $\pi^*$ and the benchmark policies $\pi^{cc*}$, $\pi^{ccr}$, $\pi^r$.} 
    \label{tab:ComparisonStrategies}
    \begin{tabular}{ll l llllllll r rrrr}
        \toprule
        && \multicolumn{9}{c}{Shipment policy} &	\multicolumn{1}{c}{Profit} &\multicolumn{4}{c}{Benefit\tnote{1} [\%] compared to} \\
        \cmidrule{3-11} \cmidrule{13-16} 
       Setting & Policy &	$\bar{\tau}$ & $f_0$ & $f_1$ & $f_2$ & $f_3$ & $f_4$ & $f_5$ & $f_6$ & $f_7$	& \multicolumn{1}{l}{$\E[G]$}	& \multicolumn{1}{l}{$\pi^r$} & \multicolumn{1}{l}{$\pi^{ccr}$} & \multicolumn{1}{l}{ $\pi^{cc*}$} & \multicolumn{1}{l}{$\pi^{2t*}$} \\
       \midrule
       $\rho=0.85$,  & $\pi^*$	&	7	&	2.2	&	2.2	&	2.2	&	2.2	&	2.2	&	2.4	&	2.6	&	3.0	&	33.03	&	11.36	&	8.69	&	2.50	&	0.51	\\
       $\beta=8$  & $\pi^{2t*}$	&	7	&	2.4	&	2.4	&	2.4	&	2.4	&	2.4	&	2.4	&	2.4	&	3.0	&	32.86	&	10.80	&	8.14	&	1.98	&		\\
        & $\pi^{cc*}$	&	7	&	2.4	&	2.4	&	2.4	&	2.4	&	2.4	&	2.4	&	2.4	&	2.4	&	32.22	&	8.64	&	6.04	&		&		\\
        & $\pi^{ccr}$	&	6	&	2.0	&	2.0	&	2.0	&	2.0	&	2.0	&	2.0	&	2.0	&		&	30.39	&	2.45	&		&		&		\\
        & $\pi^r$	&	7	&	2.0	&	2.0	&	2.0	&	2.0	&	2.0	&	2.0	&	2.0	&	2.0	&	29.66	&		&		&		&		\\
        \midrule
        $\rho=0.85$, &$\pi^*$	&	7	&	2.2	&	2.4	&	2.4	&	2.4	&	2.4	&	2.4	&	2.6	&	3.2	&	30.88	&	26.09	&	9.96	&	3.97	&	0.34	\\
        $\beta=12$ &$\pi^{2t*}$	&	7	&	2.4	&	2.4	&	2.4	&	2.4	&	2.4	&	2.4	&	2.4	&	3.2	&	30.78	&	25.67	&	9.59	&	3.63	&		\\
        &$\pi^{cc*}$	&	6	&	2.4	&	2.4	&	2.4	&	2.4	&	2.4	&	2.4	&	2.4	&		&	29.70	&	21.27	&	5.76	&		&		\\
        &$\pi^{ccr}$	&	6	&	2.0	&	2.0	&	2.0	&	2.0	&	2.0	&	2.0	&	2.0	&		&	28.08	&	14.67	&		&		&		\\
        &$\pi^r$	&	7	&	2.0	&	2.0	&	2.0	&	2.0	&	2.0	&	2.0	&	2.0	&	2.0	&	24.49	&		&		&		&		\\
        \midrule
        $\rho=0.9$, &$\pi^*$	&	7	&	2.4	&	2.4	&	2.4	&	2.6	&	2.6	&	2.6	&	2.8	&	3.2	&	25.62	&	38.84	&	22.26	&	2.17	&	0.09	\\
        $\beta=8$ &$\pi^{2t*}$	&	7	&	2.6	&	2.6	&	2.6	&	2.6	&	2.6	&	2.6	&	2.6	&	3.2	&	25.60	&	38.72	&	22.16	&	2.08	&		\\
        &$\pi^{cc*}$	&	7	&	2.6	&	2.6	&	2.6	&	2.6	&	2.6	&	2.6	&	2.6	&	2.6	&	25.08	&	35.89	&	19.67	&		&		\\
        &$\pi^{ccr}$	&	6	&	2.0	&	2.0	&	2.0	&	2.0	&	2.0	&	2.0	&	2.0	&		&	20.95	&	13.56	&		&		&		\\
        &$\pi^r$	&	7	&	2.0	&	2.0	&	2.0	&	2.0	&	2.0	&	2.0	&	2.0	&	2.0	&	18.45	&		&		&		&		\\
        \midrule
        $\rho=0.9$, &$\pi^*$	&	7	&	2.6	&	2.6	&	2.6	&	2.6	&	2.6	&	2.8	&	3.0	&	3.4	&	21.22	&	176.29	&	43.35	&	3.86	&	0.71	\\
        $\beta=12$ &$\pi^{2t*}$	&	7	&	2.8	&	2.8	&	2.8	&	2.8	&	2.8	&	2.8	&	2.8	&	3.4	&	21.07	&	174.35	&	42.34	&	3.13	&		\\
        &$\pi^{cc*}$	&	6	&	2.6	&	2.6	&	2.6	&	2.6	&	2.6	&	2.6	&	2.6	&		&	20.43	&	166.02	&	38.02	&		&		\\
        &$\pi^{ccr}$	&	5	&	2.0	&	2.0	&	2.0	&	2.0	&	2.0	&	2.0	&		&		&	14.80	&	92.74	&		&		&		\\
        &$\pi^r$	&	7	&	2.0	&	2.0	&	2.0	&	2.0	&	2.0	&	2.0	&	2.0	&	2.0	&	7.68	&		&		&		&		\\
        \midrule
        $\rho=0.95$, &$\pi^*$	&	7	&	3.0	&	3.0	&	3.0	&	3.0	&	3.0	&	3.0	&	3.0	&	3.4	&	6.91	&	163.34	&	439.89	&	3.67	&	0.00	\\
        $\beta=8$ &$\pi^{2t*}$	&	7	&	3.0	&	3.0	&	3.0	&	3.0	&	3.0	&	3.0	&	3.0	&	3.4	&	6.91	&	163.34	&	439.89	&	3.67	&		\\
        &$\pi^{cc*}$	&	7	&	3.0	&	3.0	&	3.0	&	3.0	&	3.0	&	3.0	&	3.0	&	3.0	&	6.67	&	161.10	&	427.87	&		&		\\
        &$\pi^{ccr}$	&	2	&	2.0	&	2.0	&	2.0	&		&		&		&		&		&	-2.03	&	81.36	&		&		&		\\
        &$\pi^r$	&	7	&	2.0	&	2.0	&	2.0	&	2.0	&	2.0	&	2.0	&	2.0	&	2.0	&	-10.91	&		&		&		&		\\
        \midrule
        $\rho=0.95$, &$\pi^*$	&	7	&	3.2	&	3.2	&	3.2	&	3.2	&	3.2	&	3.2	&	3.4	&	3.6	&	-2.98	&	91.81	&	72.28	&	6.43	&	0.69	\\
        $\beta=12$ &$\pi^{2t*}$	&	7	&	3.2	&	3.2	&	3.2	&	3.2	&	3.2	&	3.2	&	3.2	&	3.6	&	-3.00	&	91.75	&	72.08	&	5.78	&		\\
        &$\pi^{cc*}$	&	6	&	3.2	&	3.2	&	3.2	&	3.2	&	3.2	&	3.2	&	3.2	&		&	-3.18	&	91.24	&	70.37	&		&		\\
        &$\pi^{ccr}$	&	1	&	2.0	&	2.0	&		&		&		&		&		&		&	-10.75	&	70.45	&		&		&		\\
        &$\pi^r$	&	7	&	2.0	&	2.0	&	2.0	&	2.0	&	2.0	&	2.0	&	2.0	&	2.0	&	-36.37	&		&		&		&		\\
        \midrule
        \midrule
        Median & $\pi^*$ &&&&&&&&&&& 65.32	&	32.80	&	3.77	&	0.42 \\ 
        & $\pi^{2t*}$	&&&&&&&&&&& 65.24	&	32.25	&	3.38	& \\ 
        &  $\pi^{cc*}$	&&&&&&&&&&& 63.57	&	28.85	&		& \\
        & $\pi^{ccr}$	&&&&&&&&&&& 42.56	&		&		& \\
        \bottomrule
    \end{tabular}
    \begin{tablenotes}
        \item[1] {\footnotesize Benefit is calculated as the relative deviation in the expected profit with respect to a given policy.}
    \end{tablenotes}
\end{threeparttable} 
\end{table}

Table~\ref{tab:ComparisonStrategies} gives the policy parameters and expected profit of $\pi^{2t*}$, $\pi^*$, and the benchmark policies $\pi^{cc*}$, $\pi^{ccr}$, and $\pi^r$. 
Comparing $\pi^{2t*}$ and $\pi^*$, we note that the optimality gap of $\pi^{2t*}$ is 0.4\% in the median and 0.7\% at maximum. In system setting $(\rho = 0.95, c=8)$, $\pi^{2t*}$ is optimal. These results indicate that e-commerce companies can achieve near-optimal expected profits by operating the two-level time-dependent shipment policy that is much easier to understand for customers than $\pi^*$. 

Next, we compare $\pi^{2t*}$ with the benchmark policies $\pi^{cc*}$, $\pi^{ccr}$, and $\pi^r$. For instance, in system setting $(\rho = 0.95,c=12)$, $\pi^{2t*}$ achieves a 5.8\% higher expected profit per cycle than $\pi^{cc*}$. As expected, $\pi^{2t*}$ has an even higher benefit of 72.1\% and 91.8\% compared to $\pi^{ccr}$ and $\pi^r$, respectively. Similar results apply to the other system settings. In the median, $\pi^{2t*}$ outperforms $\pi^{cc*}$ by 3.4\%, $\pi^{ccr}$ by 32.3\%, and $\pi^r$ by 65.2\%. 
Furthermore, there is a clear ranking of the benchmark policies, i.e., $\pi^r$ achieves the lowest expected profit, $\pi^{ccr}$ improves the expected profit by 42.6\% (compared to $\pi^r$), and  $\pi^{cc*}$ further improves the expected profit by 28.9\% (compared to $\pi^{ccr}$). 
These results indicate that 
e-commerce companies can improve the expected profit by considerable margins when evolving the shipment policy as follows: 
(1) offer time-dependent instead of constant shipment options by introducing a cutoff point for express shipment; 
(2) charge the profit-maximizing instead of the revenue-maximizing express shipment fee; and 
(3) switch from a constant to a two-level time-dependent express shipment fee. 

\subsection{Optimal policy parameters of two-level time-dependent shipment policy}
\label{sec:AnalysisOptDesign}
We analyze the optimal policy parameters of the two-level time-dependent shipment policy $\pi^{2t*}$ given in Table~\ref{tab:ComparisonStrategies}. In each system setting,  the optimal cutoff point and fee switching point are $\bar{\tau}^* = 7$ and $\hat{\tau}^*=6$. Noting that $\bar{\tau}^* = T-1$ and $\hat{\tau}^*=\bar{\tau}^*-1$, these results indicate that e-commerce companies should offer express shipment until the last time period before the deadline and switch from the express to the last-minute express fee in the last time period before the cutoff point. In this way, all customers have the choice between express and regular shipment and the vast majority of them pays the express fee when choosing express shipment; only those customers that arrive shortly before the deadline and still want express shipment pay the last-minute express fee. 
Further, in each system setting, the optimal last-minute express fee is higher than the optimal express fee, i.e., $f^{l*}>f^{e*}$, which is consistent with our analytical results in Lemmas~\ref{lem:twoLevelTimeDependentShipmentPolicy} and \ref{lem:monotonicityFee}. A higher last-minute express fee reduces the risk of receiving a high proportion of express shipment orders shortly before the deadline while increasing the revenue from customers that really want a last-minute order. 
Moreover, note that $f^{e*}$ and $f^{l*}$ are higher than the revenue-maximizing fee, which is $f^r =2$, and increase in the penalty cost rate $\beta$ and utilization $\rho$. Again, these are clear actionable insights for practitioners. 
When the system becomes more congested or a backorder becomes more expensive, the e-commerce company should increase the express shipment fees. 
\begin{figure}
    \centering
    \subfloat[System setting $(\rho=0.85,c=12)$]{%
    \begin{tikzpicture}[scale = 0.75]\small
            \draw[-stealth] (-1,0) -- (8,0);
            \draw(8,-0.4) node{$\tau$};
            \draw[-stealth] (-1,0) -- (-1,8.5);
            \draw(-1.4,8.5) node{$C_\tau$};
            \draw (0, -0.15) -- (0, 0.15);
            \draw (1, -0.15) -- (1, 0.15);
            \draw (2, -0.15) -- (2, 0.15);
            \draw (3, -0.15) -- (3, 0.15);
            \draw (4, -0.15) -- (4, 0.15);
            \draw (5, -0.15) -- (5, 0.15);
            \draw (6, -0.15) -- (6, 0.15);
            \draw (7, -0.15) -- (7, 0.15);
            \draw(0,-0.4) node{\scriptsize0};
            \draw(1,-0.4) node{\scriptsize1}; 
            \draw(2,-0.4) node{\scriptsize2};
            \draw(3,-0.4) node{\scriptsize3};
            \draw(4,-0.4) node{\scriptsize4};
            \draw(5,-0.4) node{\scriptsize5}; 
            \draw(6,-0.4) node{\scriptsize6};
            \draw(7,-0.4) node{\scriptsize7};
             \draw (-0.85, 5/2) -- (-1.15, 5/2);
             \draw(-1.5,5/2) node{\scriptsize\phantom{1}5};
             \draw (-0.85, 10/2) -- (-1.15, 10/2);
             \draw(-1.5,10/2) node{\scriptsize10};
             \draw (-0.85, 15/2) -- (-1.15, 15/2);
             \draw(-1.5,15/2) node{\scriptsize15};
            
            \filldraw[blue] (-1.1cm,0) -- (-0.9cm,0) -- (-1cm,0.2cm);
            \filldraw[blue] (-0.1cm,1.125) -- (0.1cm,1.125) -- (0,1.325cm);
            \filldraw[blue] (0.9cm,2.125) -- (1.1cm,2.125) -- (1cm,2.325cm);
            \filldraw[blue] (1.9cm,3.125) -- (2.1cm,3.125) -- (2cm,3.325cm);
            \filldraw[blue] (2.9cm,4.125) -- (3.1cm,4.125) -- (3cm,4.325cm);
            \filldraw[blue] (3.9cm,5.125) -- (4.1cm,5.125) -- (4cm,5.325cm);
            \filldraw[blue] (4.9cm,6.125) -- (5.1cm,6.125) -- (5cm,6.325cm);
            \filldraw[blue] (5.9cm,7) -- (6.1cm,7) -- (6cm,7.2cm);
            \filldraw[blue] (6.9cm,7.5) -- (7.1cm,7.5) -- (7cm,7.7cm);
            \draw[dashed, thick, blue] (-1,0) -- (0,2.25/2);
            \draw[dashed, thick, blue] (0,2.25/2) -- (5,12.25/2);
            \draw[dashed, thick, blue] (5,12.25/2) -- (6,14/2);
            \draw[dashed, thick, blue] (6,14/2) -- (7,15/2);
            
            \filldraw[red] (-1,0) circle (2.5pt) node[anchor=west]{};
            \filldraw[red] (0,2/2) circle (2.5pt) node[anchor=west]{};
            \filldraw[red] (1,4/2) circle (2.5pt) node[anchor=west]{};
            \filldraw[red] (2,6/2) circle (2.5pt) node[anchor=west]{};
            \filldraw[red] (3,8/2) circle (2.5pt) node[anchor=west]{};
            \filldraw[red] (4,10/2) circle (2.5pt) node[anchor=west]{};
            \filldraw[red] (5,12/2) circle (2.5pt) node[anchor=west]{};        
            \filldraw[red] (6,14/2) circle (2.5pt) node[anchor=west]{};        
            \filldraw[red] (7,15/2) circle (2.5pt) node[anchor=west]{};
            \draw[dashed, thick, red] (-1,0) -- (6,14/2);
            \draw[dashed, thick, red] (6,14/2) -- (7,15/2);

            \filldraw[red] (6.35,1.3) circle (2.5pt); 
            \draw (7.15,1.3) node {\tiny \textcolor{black}{$\pi^{2t*}$}}; 
            \filldraw[blue] (6.25cm,0.6) -- (6.45cm,0.6) -- (6.35cm,0.8cm); 
            \draw (7,0.7) node {\tiny \textcolor{black}{$\pi^{*}$}};   
            \draw[gray] (5.90,1.65) rectangle (7.8,0.4);

        \end{tikzpicture}
    }
    \hfill
    \subfloat[System setting $(\rho=0.9,c=8)$]{%
    \begin{tikzpicture}[scale = 0.75]\small
            \draw[-stealth] (-1,0) -- (8,0);
            \draw(8,-0.4) node{$\tau$};
            \draw[-stealth] (-1,0) -- (-1,8.5);
            \draw(-1.4,8.5) node{$C_\tau$};
            \draw (0, -0.15) -- (0, 0.15);
            \draw (1, -0.15) -- (1, 0.15);
            \draw (2, -0.15) -- (2, 0.15);
            \draw (3, -0.15) -- (3, 0.15);
            \draw (4, -0.15) -- (4, 0.15);
            \draw (5, -0.15) -- (5, 0.15);
            \draw (6, -0.15) -- (6, 0.15);
            \draw (7, -0.15) -- (7, 0.15);
            \draw(0,-0.4) node{\scriptsize0};
            \draw(1,-0.4) node{\scriptsize1}; 
            \draw(2,-0.4) node{\scriptsize2};
            \draw(3,-0.4) node{\scriptsize3};
            \draw(4,-0.4) node{\scriptsize4};
            \draw(5,-0.4) node{\scriptsize5}; 
            \draw(6,-0.4) node{\scriptsize6};
            \draw(7,-0.4) node{\scriptsize7};
             \draw (-0.85, 5/2) -- (-1.15, 5/2);
             \draw(-1.5,5/2) node{\scriptsize\phantom{1}5};
             \draw (-0.85, 10/2) -- (-1.15, 10/2);
             \draw(-1.5,10/2) node{\scriptsize10};
             \draw (-0.85, 15/2) -- (-1.15, 15/2);
             \draw(-1.5,15/2) node{\scriptsize15};
            
          \filldraw[blue] (-1.1cm,0) -- (-0.9cm,0) -- (-1cm,0.2cm);
            \filldraw[blue] (-0.1cm,1) -- (0.1cm,1) -- (0,1.2cm);
            \filldraw[blue] (0.9cm,2) -- (1.1cm,2) -- (1cm,2.2cm);
            \filldraw[blue] (1.9cm,3) -- (2.1cm,3) -- (2cm,3.2cm);
            \filldraw[blue] (2.9cm,3.85) -- (3.1cm,3.85) -- (3cm,4.05cm);
            \filldraw[blue] (3.9cm,4.75) -- (4.1cm,4.75) -- (4cm,4.95cm);
            \filldraw[blue] (4.9cm,5.625) -- (5.1cm,5.625) -- (5cm,5.825cm);
            \filldraw[blue] (5.9cm,6.375) -- (6.1cm,6.375) -- (6cm,6.575cm);
            \filldraw[blue] (6.9cm,6.85) -- (7.1cm,6.85) -- (7cm,7.05cm);
            \draw[dashed, thick, blue] (-1,0) -- (2,6/2);
            \draw[dashed, thick, blue] (2,6/2) -- (6,12.75/2);
            \draw[dashed, thick, blue] (6,12.75/2) -- (7,13.75/2);
            
            \filldraw[red] (-1,0) circle (2.5pt) node[anchor=west]{};
            \filldraw[red] (0,1.75/2) circle (2.5pt) node[anchor=west]{};
            \filldraw[red] (1,3.5/2) circle (2.5pt) node[anchor=west]{};
            \filldraw[red] (2,5.25/2) circle (2.5pt) node[anchor=west]{};
            \filldraw[red] (3,7/2) circle (2.5pt) node[anchor=west]{};
            \filldraw[red] (4,8.75/2) circle (2.5pt) node[anchor=west]{};
            \filldraw[red] (5,10.5/2) circle (2.5pt) node[anchor=west]{};        
            \filldraw[red] (6,12.25/2) circle (2.5pt) node[anchor=west]{};        
            \filldraw[red] (7,13.25/2) circle (2.5pt) node[anchor=west]{};
            \draw[dashed, thick, red] (-1,0) -- (6,12.25/2);
            \draw[dashed, thick, red] (6,12.25/2) -- (7,13.25/2);

            \filldraw[red] (6.35,1.3) circle (2.5pt); 
            \draw (7.15,1.3) node {\tiny \textcolor{black}{$\pi^{2t*}$}}; 
            \filldraw[blue] (6.25cm,0.6) -- (6.45cm,0.6) -- (6.35cm,0.8cm); 
            \draw (7,0.7) node {\tiny \textcolor{black}{$\pi^{*}$}};   
            \draw[gray] (5.90,1.65) rectangle (7.8,0.4);
    
    \end{tikzpicture}
    }
    \caption{Express shipment demand profiles of the optimal two-level time-dependent shipment policy and the optimal general time-dependent shipment policy (dashed line for guidance only).}
    \label{fig:comparisonDemandProfiles}
\end{figure}

Comparing the optimal policy parameters of $\pi^{2t*}$ to those of the optimal general time-dependent shipment policy $\pi^*$, we note that $\pi^{2t*}$ and $\pi^*$ set the same cutoff point $\bar{\tau} = 7$, charge the same express shipment fee  at the cutoff point, and adjust the express shipment fee in the time period before the cutoff point, i.e., $\tau = 6$ is a fee switching point. While this is the unique fee switching point of $\pi^{2t*}$, $\pi^*$ additionally adjusts the express shipment fee at earlier ages of the operating cycle. Except for system setting $(\rho=0.95, c=8)$, where $\pi^*=\pi^{2t*}$, $\pi^*$ has two or three fee switching points per cycle. 
As a consequence, the express shipment demand profiles of $\pi^{2t*}$ and $\pi^*$ are very similar. 
For instance, in system setting $(\rho=0.85,c=12)$ illustrated in Figure \ref{fig:comparisonDemandProfiles}a, $\pi^{2t*}$ and $\pi^*$ have the same total expected express shipment orders per cycle, i.e., $C^{\pi^{2t*}}_{T-1}=C^{\pi^*}_{T-1}$, and equivalently the same average express shipment fee of $\bar{f} = 2.5$.  However, as at every age $\tau\in\{0,1,...,\bar{\tau}-1\}$, $\pi^*$ has a lower early-cycle average express shipment fee than $\pi^{2t*}$, $\pi^*$ has a higher fraction of express shipment orders at early ages of the operating cycle, resulting in a 0.34\% higher expected profit than $\pi^{2t*}$ (Table~\ref{tab:ComparisonStrategies}).
Conversely, in system setting $(\rho=0.9,c=8)$ illustrated in Figure~\ref{fig:comparisonDemandProfiles}b, $\pi^{*}$ has more total expected express shipment orders per cycle, i.e., $C^{\pi^{*}}_{T-1}>C^{\pi^{2t*}}_{T-1}$, which coincides with a lower average express shipment fee than $\pi^{2t*}$ ($\bar{f}^{\pi^{*}} = 2.625 < 2.675 = \bar{f}^{\pi^{2t*}}$). Nevertheless, $\pi^*$ achieves a 0.09\% higher expected profit than $\pi^{2t*}$ (Table~\ref{tab:ComparisonStrategies}) by charging lower express shipment fees than $\pi^{2t*}$ at early ages and weakly higher fees at late ages of the operating cycle.

\section{Concluding remarks} 
\label{sec:Conclusions}
Fulfilling online orders quickly is increasingly vital in e-commerce, with many companies offering same-day shipment services to meet rising customer expectations. However, companies risk overpromising their shipment services when the offered shipment options and corresponding fees are not aligned with the ability to fulfill these orders on time. To mitigate this risk, we propose a transparent time-dependent shipment policy where the shipment options and corresponding fees depend on the remaining time until the next shipment deadline. These shipment options and fees are consistent each day and can be easily communicated on the ordering platform, so that any change in the shipment policy is transparent and predictable to customers.  

We considered multiple operating cycles of an e-commerce fulfillment center, where each cycle consists of a fixed number of time periods and ends with a deadline upon which orders which are due for shipment by this deadline need to be handed over to the parcel delivery company. Throughout each operating cycle, utility-maximizing customers arrive according to a Poisson process and processing capacity is stochastic. If the available capacity is insufficient to process  all orders due for shipment in a given operating cycle, these orders are said to be late and carried over for shipment in the subsequent operating cycle. 
We developed a discrete-time period Markov chain model for exact steady-state performance analysis of transparent time-dependent shipment policies. 
As direct analysis of transparent time-dependent shipment policies in the price domain is analytically and computationally intractable, we introduced a novel problem-domain transformation that maps each policy to its induced cumulative express demand profile. This transformation proved central to our analysis, as it both enabled the derivation of several structural properties—showing that optimal fees increase over time and that introducing a cutoff point and time-dependent fees enhances profit—and revealed a supermodular structure of the profit function, which allows the optimal transparent shipment policy to be computed in polynomial time.
Beyond these analytical results, we proposed the two-level time-dependent shipment policy that differentiates two fees for express shipment and switches between these once per cycle. As demonstrated numerically, this policy achieves near-optimal profits while remaining easy to communicate and market in practice.

Our results demonstrate that transparent time-dependent shipment policies are a customer-friendly alternative to opaque dynamic pricing, allowing e-commerce companies to align express shipment demand with fulfillment capacity while maintaining fairness and simplicity. We provide practitioners with clear, actionable guidelines on how to design the shipment options and fees of the transparent time-dependent shipment policy.


\bibliographystyle{informs2014}
\bibliography{bibliography}

\newpage
\section*{Online companion paper}
\begin{APPENDICES}

\section{Time-dependent customer demand} \label{sec:modelExtension}
In our model, the customer demand $D_t$, $t\in\mathbb{N}_0$, is i.i.d. drawn from a stationary Poisson distribution with expected customer demand $\lambda$. However, in practice,  customer demand is likely to follow some time-dependent pattern. In this section, we discuss how our model and analytical results extend to different types of time-dependent customer demand. 

Customer demand may depend on the time of the day. We can extend our model to such demand patterns by modeling an age-dependent expected demand rate, i.e., $\lambda_t = \lambda_\tau$ where $\tau = t \Mod T$. Except for the monotonicity of the express shipment fee (Lemma~\ref{lem:monotonicityFee}, Lemma~\ref{lem:twoLevelTimeDependentShipmentPolicy}), all analytical results in \S\ref{sec:analyticalResults} continue to hold with the following minor modifications.
In Lemma~\ref{lem:twoLevelTimeDependentShipmentPolicyComparedConstant}, the express shipment fee $f$ of the constant cutoff point shipment policy is a convex combination of the fees $f^e$ and $f^l$ of the two-level time-dependent shipment policy, i.e., $f= \epsilon f^e + (1-\epsilon)f^l$, with $\epsilon$ being demand-dependent now, i.e., $\epsilon = \left(\sum_{\tau=0}^{\hat{\tau}}\lambda_\tau\right) / \left(\sum_{\tau=0}^{\bar{\tau}}\lambda_\tau\right)$. 
In Lemma~\ref{lem:timeDependentShipmentPolicy}, the assumptions refer to \emph{weighted} average shipment fees now. More specifically, suppose that the shipment policies $\pi^{cc} = (\bar{\tau}, \bar{f},\dots, \bar{f})\in\Pi^{cc}$, $\pi = (\bar{\tau}, f_0, f_1, \dots, f_{\bar{\tau}})\in\Pi$, and $\pi' = (\bar{\tau}, f'_0, f'_1, \hdots, f'_{\bar{\tau}})\in\Pi$ have the same \emph{weighted} average express shipment fee and for each $\tau \in\{ 1,\dots, \bar{\tau}\}$, the early-cycle \emph{weighted} average express shipment fee is lowest under $\pi$, higher under $\pi'$, and highest under $\pi^{cc}$, i.e., 
\begin{equation*}
    \bar{f} = \frac{\sum_{j=0}^{\bar{\tau}} \lambda_jf_j}{\sum_{j=0}^{\bar{\tau}} \lambda_j} = \frac{\sum_{j=0}^{\bar{\tau}} \lambda_jf'_j}{\sum_{j=0}^{\bar{\tau}} \lambda_j} 
    \quad \quad \quad \quad \text{and} \quad \quad \quad \quad
    \frac{\sum_{j=0}^{\tau-1} \lambda_jf_j}{\sum_{j=0}^{\tau-1} \lambda_j} \leq \frac{\sum_{j=0}^{\tau-1} \lambda_jf'_j}{\sum_{j=0}^{\tau-1} \lambda_j} \leq \bar{f} \quad \quad \forall \tau \in\{ 1,\dots, \bar{\tau}\}.
    \end{equation*}
Then, $\E[B^{\pi}] \le \E[B^{\pi'}] \le \E[B^{\pi^{cc}}]$, i.e., $\pi$ has the fewest expected backorders and $\pi^{cc}$ the most. 

Moreover, the customer demand during a certain period of the year may be higher than during the rest of the year. A typical example for such peak periods is Christmas sales in December. Peak periods are usually known in advance and sufficiently long so that the e-commerce company can adjust their processing capacity (e.g., by hiring temporary workers) to the expected demand increase and maintain a stable system. As the duration of a peak period (e.g., multiple weeks) is long relative to the duration of a single operating cycle (e.g., one day), the system reaches steady state. Thus, computing our model with the adjusted parameters $\lambda$ and $K$ provides us with an exact performance analysis of the fulfillment center during the peak period and allows us to identify a suitable shipment policy for the peak period. The peak and off-peak shipment policies can be transparently communicated on the ordering platform. 
 
\section{Proofs}\label{AppendixProofs}
In this section, we provide the proofs of the analytical properties established in \S\ref{sec:analyticalResults}. We start with some preliminary analysis in \ref{sec:interplayBacklogShipmentPolicy} and \ref{sec:feeStructures} that allow us to present the proofs more conveniently. In \ref{sec:proofsShipmentPolicies}, we provide the proofs of the analytical properties of the shipment policies (\S\ref{sec:constantShipmentPolicy}-\ref{sec:timeDependentShipmentPolicy}). In \ref{sec:proofsOptimization}, we provide the proofs related to optimizing the shipment fees (\S\ref{sec:optimizingshipmentfees}).  

\subsection{Evolution of order backlog and backorders and the interplay with the shipment policy}
\label{sec:interplayBacklogShipmentPolicy}

We characterize the effects of the shipment policy $\pi=(\bar{\tau},f_0, f_1, \dots, f_{\bar{\tau}})$ on the evolution of the order backlog and the backorders. 
We first note from (\ref{eq:transitionTotalBacklog}) that the total order backlog $X^\Sigma_t$ is independent of $\pi$ since $X^\Sigma_t$ only depends on the customer demand $D_t$ and the processing capacity $K_t$, which are independent of $\pi$. The high-urgency order backlog $X^h_{t}$ at the beginning of each operating cycle $l$, i.e., at times $t = lT$, is also independent of $\pi$ since $X^h_{lT} = X^\Sigma_{lT}$; see (\ref{eq:transitionHighUrgencyBacklog}).  
We can therefore restrict the analysis to the effects of $\pi$ on the high-urgency order backlog $X^h_t$ in time periods $t=lT+\tau$, $\tau\in\{1,2\dots,T-1\}$, of any arbitrarily chosen operating cycle $l$ in the following. 

To ease the burden of notation, we refer to these time periods by $\tau$ only. In particular, the deadline at the end of operating cycle $l$ is denoted by $T$, and we refer to the quantities $X^h_{lT+\tau}$, $K_{lT+\tau}$, and $E_{lT+\tau}$ by $X^h_\tau$, $K_\tau$, and $E_\tau$, respectively. 
Let $E[\tau_1,\tau_2]:=\sum_{j=\tau_1}^{\tau_2}E_j$ and $K[\tau_1,\tau_2]:=\sum_{j=\tau_1}^{\tau_2}K_j$ denote the \emph{cumulative express shipment orders} and \emph{cumulative processing capacity} in operating cycle $l$ between $\tau_1\in\{0,\ldots,T-1\}$ and $\tau_2\in \{\tau_1,\ldots,T-1\}$, respectively. 

Consider the e-fulfillment center at any age $\tau\in\{1,2,\dots,T-1\}$. 
The processing capacity $K_\tau$ breaks down into three categories: (a) capacity used to process orders that are due by deadline $T$ or already too late, (b) capacity used to process orders that are due by deadline $2T$, (c) idle capacity. We call the processing capacity that falls into categories (b) and (c) \emph{low-urgency processing capacity}, denoted by $L_\tau$. Note that 
\begin{equation} \label{eq:lowUrgencyProcessingCapacity}
    L_\tau=(K_\tau -X^h_{\tau-1}-E_\tau)^+.
\end{equation}
Further, $L[\tau_1,\tau_2]:=\sum_{j=\tau_1}^{\tau_2}L_j$ denotes the \emph{cumulative low-urgency processing capacity} in operating cycle $l$ between $\tau_1\in\{0,\ldots,T-1\}$ and $\tau_2\in \{\tau_1,\ldots,T-1\}$. 

The cumulative processing capacity $K[0,\tau]$ up to age $\tau$ is either used to process orders that are due by deadline $T$ or already too late, or it is low-urgency processing capacity. If $L_\tau>0$ for some $\tau\in \{1,\ldots,T-1\}$, we know that the cumulative express shipment orders $E[0,\tau]$ up to age $\tau$ have been processed since high-urgency order backlog and express shipment orders are prioritized over regular shipment orders. Thus, we have $K[0,\tau]=X^h_0+E[0,\tau]+L[0,\tau]$ and equivalently $L[0,\tau] =  K[0,\tau]-E[0,\tau]-X^h_0$ if $L_{\tau}>0$. Further, note that $L[0,\tau]=L[0,\tau-1]$ if $L_{\tau}=0$ and $L[0,\tau]$ is non-decreasing in $\tau$. Then, the cumulative low-urgency processing capacity $L[0,\tau]$ up to age $\tau$ satisfies 
\begin{equation}\label{eq:lowUrgencyProcessingCapacityCumulative}
  L[0,\tau] = \max_{\tau'\in\{0,\ldots,\tau\}} \left( K[0,\tau']-E[0,\tau']-X^h_0\right)^+.
\end{equation}

The high-urgency order backlog $X_\tau^h$ at age $\tau$ can be rewritten as follows:
\begin{equation} \label{eq:highUrgencyBacklogTau}
    X^h_\tau=(X^h_{\tau-1}+E_\tau-K_\tau)^+ = X^h_{\tau-1}+E_\tau-K_\tau + L_\tau = X_0^h + E[0,\tau]-K[0,\tau] + L[0,\tau].
\end{equation}
Recalling that $X^h_0$ and $K[0,\tau]$ are independent of $\pi$, it follows together with (\ref{eq:lowUrgencyProcessingCapacityCumulative}) that the high-urgency order backlog $X_\tau^h$ at age $\tau$ depends on $\pi$ by the express shipment orders in operating cycle $l$ up to age $\tau$, i.e., $\{E_0, E_1, \dots, E_\tau\}$. 

Similarly, the backorders $B_l$ at the end of operating cycle $l$ are given by
\begin{align} 
    B_l 
    = X^h_{T-1}+E_{T-1}-K_{T-1} + L_{T-1}
    = X_0^h + E[0,T-1] - K[0,T-1] + L[0,T-1],\label{eq:backorders_2}
\end{align}
and depend on $\pi$ by the express shipment orders in operating cycle $l$, i.e., $\{E_0, E_1, \dots, E_{\bar{\tau}}\}$.

\subsection{Fee structure and equivalent shipment policy}
\label{sec:feeStructures}
We introduce the notion of fee structures (Definition~\ref{def:feeStructure}) and establish some analytical results on the fee structures (Lemma~\ref{lem:feeStructure}, Proposition~\ref{prop:dominance}). These results help us in \ref{sec:proofsShipmentPolicies} to more conveniently present the proofs of the analytical properties of the shipment policies. 

\begin{definition}[Fee structure] \label{def:feeStructure}
\phantom{texttexttexttext}
\begin{enumerate}
    \item[(a)] Let $f_\tau\in[\umin,\umax]$ $\forall \tau\in\{0,1,\dots,T-1\}$ and $f_T = \umax$. We denote by $\phi:=(f_0,f_1,\dots,f_{T-1},f_T)$ a \emph{fee structure}. The set of fee structures is given by $\Phi: = \{(f_0,f_1,\dots,f_{T-1},\umax) \mid f_\tau\in[\umin,\umax] \forall \tau\in\{0,\dots,T-1\}\}$.
    \item[(b)]  Let $\phi\in\Phi$. We denote by $\mathbf{C}^\phi:= (C_{-1}^\phi, C_{0}^\phi, C_{1}^\phi,\dots, C_{T-1}^\phi)$ where $C_{-1}^\phi = 0$ and $C_\tau^\phi = \E[E^\phi[0,\tau]]$ $\forall \tau\in\{0,1,\dots,T-1\}$ the \textit{express shipment demand profile} of fee structure $\phi$.
    \item[(c)] Let $\pi\in\Pi$, $\phi\in\Phi$, and $\mathbf{C}^\pi$ and $\mathbf{C}^\phi$ the associated express shipment demand profiles. We call $\pi$ and $\phi$ \emph{equivalent} if $C_\tau^\pi = C_\tau^\phi$ $\forall \tau=0,1,\dots,T-1$.\hfill \Halmos \endproof
    \end{enumerate}
\end{definition}

Note that the express shipment demand profile $\mathbf{C}^\phi$ of any $\phi\in\Phi$ satisfies
\begin{align}
    C_\tau^\phi - C_{\tau-1}^\phi = \E[E_\tau^\phi] &= w(f_\tau)\lambda \quad \quad \quad \quad \forall \tau \in \{0,1,\dots,T-1\}. 
     \label{eq:demandProfileFeeStructure}
\end{align}

\begin{lemma}[Equivalent fee structure and shipment policy] \label{lem:feeStructure}
\phantom{texttexttexttext}
\begin{enumerate}
    \item[(1)] Let $\pi=(\bar{\tau},f_0,f_1,\dots,f_{\bar{\tau}})\in\Pi$. The fee structure $\phi=(f'_0,f'_1,\dots,f'_{T-1},\umax)$ with $f'_\tau = f_\tau$ $\forall \tau=0,1,\dots,\bar{\tau}$ and $f'_\tau = \umax$ $\forall \tau=\bar{\tau}+1,\dots,T-1$ is the unique equivalent fee structure of shipment policy $\pi$. 
    \item[(2)] Let $\phi=(f'_0,f'_1,\dots,f'_{T-1},\umax)\in\Phi$. The shipment policy $\pi=(\bar{\tau},f_0,f_1,\dots,f_{\bar{\tau}})\in\Pi$ with $\bar{\tau} = \min\{\tau\in\{0,1,\dots,T-1\} \mid f_{\tau+k} = \umax \forall k = 1,2,\dots,T-\tau\}$ and $f_\tau = f'_\tau$ $\forall \tau=0,1,\dots,\bar{\tau}$ is the unique equivalent shipment policy of fee structure $\phi$. 
\end{enumerate}
\end{lemma}

\proof{Proof of Lemma~\ref{lem:feeStructure}.} 
Let $\pi=(\bar{\tau},f_0,f_1,\dots,f_{\bar{\tau}})\in\Pi$ and $\phi=(f'_0,f'_1,\dots,f'_{T-1},\umax)\in\Phi$. We know from Definition~\ref{def:feeStructure} that $\pi$ and $\phi$ are equivalent if $C_\tau^\pi = C_\tau^\phi$ $\forall \tau=0,1,\dots,T-1$. 
As $C_{-1}^\pi = C_{-1}^\phi = 0$ by definition, it suffices to show that 
\begin{align}
    C_\tau^\pi-C_{\tau-1}^\pi = C_\tau^\phi-C_{\tau-1}^\phi \quad \quad \quad \quad \forall \tau\in\{0,1,\dots,T-1\}. \label{eq:conditionEquivalentPiPhi}
\end{align}
Recall from (\ref{eq:demandProfileShipmentPolicy}) and (\ref{eq:demandProfileFeeStructure}) the express shipment demand profiles $\mathbf{C}^\pi$ and $\mathbf{C}^\phi$ associated with any $\pi\in\Pi$ and any $\phi\in\Phi$, respectively:  
\begin{align}
    C_\tau^\pi - C_{\tau-1}^\pi &= \E[E_\tau^\pi] = \begin{cases}
        w(f_\tau)\lambda &\text{if }\tau \leq \bar{\tau}, \\
        0 & \text{otherwise} \\
    \end{cases}\quad \quad \quad \quad &\forall &\tau \in \{0,1,\dots,T-1\}. \label{eq:demandProfilePi} \\
  C_\tau^\phi - C_{\tau-1}^\phi &= \E[E_\tau^\phi] = w(f'_\tau)\lambda\quad \quad \quad \quad &\forall &\tau \in \{0,1,\dots,T-1\}. \label{eq:demandProfilePhi}
\end{align}

\textit{Part (1)}: Suppose that $\pi$ is given. We build on (\ref{eq:conditionEquivalentPiPhi}) and (\ref{eq:demandProfilePi})-(\ref{eq:demandProfilePhi}) to derive the associated equivalent fee structure $\phi$. 
For every $\tau \in \{0,1,\dots,\bar{\tau}\}$, we find:
\begin{align*}
    C_\tau^\phi-C_{\tau-1}^\phi = C_\tau^\pi-C_{\tau-1}^\pi 
    \quad \overset{}{\Leftrightarrow} \quad
   w(f'_\tau)\lambda = w(f_\tau)\lambda 
    \quad \overset{}{\Leftrightarrow} \quad 
   &f'_\tau = f_\tau, 
\end{align*}
where the last equivalence follows as $w$ is strictly decreasing in $f$ for $f\in[\umin,\umax]$ by (\ref{eq:willingnessToPay}), and $f_\tau,f'_\tau\in [\umin,\umax]$ $\forall \tau\in\{0,1,\dots,\bar{\tau}\}$ by $\pi\in\Pi$ and $\phi\in\Phi$.
For every $\tau \in \{\bar{\tau}+1,\dots,T-1\}$, we have:
\begin{align*}
    C_\tau^\phi-C_{\tau-1}^\phi = C_\tau^\pi-C_{\tau-1}^\pi 
    \quad \overset{}{\Leftrightarrow} \quad
   w(f'_\tau)\lambda = 0 
   \quad \overset{}{\Leftrightarrow} \quad 
   f'_\tau = \umax, 
\end{align*}
where the last equivalence follows as $w = 0$ for $f\geq \umax$ by (\ref{eq:willingnessToPay}), and $f'_\tau\in [\umin,\umax]$ $\forall \tau\in\{0,1,\dots,T-1\}$ by $\phi\in\Phi$.
This establishes claim (1).

\textit{Part (2)}:
Suppose that $\phi$ is given. We build on (\ref{eq:conditionEquivalentPiPhi}) and (\ref{eq:demandProfilePi})-(\ref{eq:demandProfilePhi}) to derive the associated equivalent shipment policy $\pi$. 
Let $\tilde{\tau}:=\min\{\tau\in\{0,1,\dots,T-1\} \mid f'_{\tau+k} = \umax \forall k = 1,2,\dots,T-\tau\}$. Then, we have for every $\tau \in \{0,1,\dots,T-1\}$,
\begin{align}
    C_\tau^\pi-C_{\tau-1}^\pi = C_\tau^\phi-C_{\tau-1}^\phi 
    \quad \overset{}{\Leftrightarrow} \quad
    C_\tau^\pi-C_{\tau-1}^\pi = w(f'_\tau)\lambda 
    \quad \overset{}{\Leftrightarrow} \quad
     C_\tau^\pi-C_{\tau-1}^\pi = \begin{cases}
         w(f'_\tau)\lambda &\text{if }\tau\leq \tilde{\tau},\\
         0 & \text{otherwise,}
    \end{cases}
    \label{eq:piGivenPhi}
\end{align}
where the last equivalence follows as $f'_\tau = \umax$ $\forall\tau\in\{\tilde{\tau}+1,\dots,T-1\}$ by construction of $\tilde{\tau}$ and $w(\umax) = 0$ by (\ref{eq:willingnessToPay}). Note that $\pi$ satisfies (\ref{eq:piGivenPhi}) if and only if $\bar{\tau} = \tilde{\tau}$ and $f_\tau = f'_\tau$  $\forall \tau\in\{0,\dots,\bar{\tau}\}$. This follows together with (\ref{eq:demandProfilePi}) as $w$ is strictly decreasing in $f$ for $f\in[\umin,\umax]$ by (\ref{eq:willingnessToPay}), and $f_\tau,f'_\tau\in [\umin,\umax]$ $\forall \tau\in\{0,1,\dots,\bar{\tau}\}$ by $\pi\in\Pi$ and $\phi\in\Phi$.
This establishes claim (2). 
\hfill \Halmos \endproof

\begin{proposition}[Backorders of fee structures]\label{prop:dominance}
    Let $\phi, \phi'\in \Phi$ and $\mathbf{C}^\phi$, $\mathbf{C}^{\phi'}$ the associated express shipment demand profiles. Suppose $C_\tau^\phi\geq C_\tau^{\phi'}$ $\forall \tau\in\{0,1,\dots,T-2\}$ and $C_{T-1}^\phi=C_{T-1}^{\phi'}$. \\
    Then, $\phi$ generates fewer expected backorders than $\phi'$, i.e., $\E[B^{\phi}]\leq \E[B^{\phi'}]$.    
\end{proposition}

\proof{Proof of Proposition~\ref{prop:dominance}.}
Let $\phi, \phi'\in \Phi$ with associated express shipment demand profiles $\mathbf{C}^\phi$ and $\mathbf{C}^{\phi'}$ that satisfy $C_\tau^\phi\geq C_\tau^{\phi'}$ $\forall \tau\in\{0,1,\dots,T-2\}$ and $C_{T-1}^\phi=C_{T-1}^{\phi'}$. Together with  Definition~\ref{def:feeStructure}, we have 
\begin{align}
    C_\tau^\phi\geq C_\tau^{\phi'} \quad &\Leftrightarrow \quad \E[E^\phi[0,\tau]] \geq \E[E^{\phi'}[0,\tau]] &\forall &\tau\in\{0,1,\dots,T-2\} \label{eq:demandProfilesProperty1} \\
    C_{T-1}^\phi=C_{T-1}^{\phi'} \quad&\Leftrightarrow \quad \E[E^\phi[0,T-1]] = \E[E^{\phi'}[0,T-1]].\label{eq:demandProfilesProperty2}
\end{align}
As $E^\phi[0,\tau]$ and $E^{\phi'}[0,\tau]$ are Poisson distributed and (\ref{eq:demandProfilesProperty1}) holds, we may assume the existence of Poisson random variables ${\chi}_\tau$ such that $E^{\phi'}[0,\tau]+\chi_\tau \sim E^\phi[0,\tau]$ for every $\tau\in\{0,1,\dots,T-2\}$. 

Then, recalling from (\ref{eq:lowUrgencyProcessingCapacityCumulative}) that $L[0,\tau] = \max_{\tau'\in\{0,\ldots,\tau\}} \left( K[0,\tau']-E[0,\tau']-X^h_0\right)^+$, we find for every $\tau\in\{0,1,\dots,T-2\}$, 
    \begin{align}
      L^{\phi}[0,\tau] = \max_{\tau'\in\{0,\ldots,\tau\}} \left( K[0,\tau']-E^{\phi}[0,\tau']-X^h_0\right)^+\nonumber
      &= \max_{\tau'\in\{0,\ldots,\tau\}} \left( K[0,\tau']-E^{\phi'}[0,\tau']-{\chi}_{\tau'}-X^h_0\right)^+\nonumber\\
      &\le \max_{\tau'\in\{0,\ldots,\tau\}} \left( K[0,\tau']-E^{\phi'}[0,\tau']-X^h_0\right)^+ = L^{\phi'}[0,\tau].\label{eq:inequalityLowUrgencyCapacityRealization}
    \end{align}
As (\ref{eq:inequalityLowUrgencyCapacityRealization}) holds for every realization of random variables, it also holds in expectation, i.e., 
\begin{align}\label{eq:inequalityLowUrgencyCapacityExpectation}
    \E[L^{\phi}[0,\tau]] \le \E[L^{\phi'}[0,\tau]] \quad \quad  \quad \forall \tau\in\{0,1,\dots,T-2\}.
\end{align}

Next, consider the backorders $B_l^\phi$ and $B_l^{\phi'}$ associated with $\phi$ and $\phi'$. Recall from (\ref{eq:backorders_2}) that $B_l^\phi = X^h_0 + E^\phi[0,T-1] - K[0,T-1] + L^\phi[0,T-1]$ and $B_l^{\phi'} = X^h_0 + E^{\phi'}[0,T-1] - K[0,T-1] + L^{\phi'}[0,T-1]$. As $X^h_0$ and $K[0,T-1]$ are independent of the fee structure, we have for every $l\in\mathbb{N}_0$,  
\begin{equation*}
    B_l^{\phi'}-B_l^\phi = E^{\phi'}[0,T-1] + L^{\phi'}[0,T-1]-E^{\phi}[0,T-1] - L^{\phi}[0,T-1].
\end{equation*}
Taking expectations on both sides and applying (\ref{eq:demandProfilesProperty2}) and (\ref{eq:inequalityLowUrgencyCapacityExpectation}), we find  
\begin{align*}
    &\E[B_l^{\phi'}]-\E[B_l^\phi] = \E[E^{\phi'}[0,T-1]] + \E[L^{\phi'}[0,T-1]] - \E[E^{\phi}[0,T-1]] - \E[L^{\phi}[0,T-1]]\\
    &\quad \Leftrightarrow \quad 
    \E[B_l^{\phi'}]-\E[B_l^\phi] \geq 0.
\end{align*} 
\hfill \Halmos \endproof

\subsection{Proofs of analytical properties of shipment policies (Lemmas~\ref{lem:revenueMaximizingPolicy}-\ref{lem:monotonicityFee})}
\label{sec:proofsShipmentPolicies}

\proof{Proof of Lemma~\ref{lem:revenueMaximizingPolicy}.}
This proof has two parts. Part (1) establishes claim (1) of Lemma~\ref{lem:revenueMaximizingPolicy} etc. 

\textit{Part (1):} We show that the expected revenue is separable and concave in the express shipment fees $f_\tau$, $\tau=0,1,\dots,\bar{\tau}$.
Recalling from (\ref{eq:revenue}) the definition of the expected revenue and applying (\ref{eq:expressRegularOrdersRealization}), we have
\begin{equation} \label{eq:revenue_proof}
    \E[R] = \sum_{\tau=0}^{\bar{\tau}}f_\tau \E[E_\tau] = \lambda \sum_{\tau=0}^{\bar{\tau}}f_\tau w(f_\tau), 
\end{equation}
which is separable in $f_\tau$.
Further, we have
\begin{align}
   \frac{\partial \E[R]}{\partial f_\tau} &= \lambda \left( w(f_\tau) + f_\tau w'(f_\tau) \right) &\forall &\tau \in\{0,1,\dots,\bar{\tau}\} \label{eq:revenueDerivative}\\
    \frac{\partial^2 \E[R]}{\partial^2 f_\tau} &= \lambda  \left( 2w'(f_\tau) + f w''(f_\tau) \right)&\forall &\tau \in\{0,1,\dots,\bar{\tau}\}. \nonumber
\end{align}
As $w'(f_\tau) = - (\umax-\umin)^{-1}< 0$ and $w''(f_\tau) \equiv 0$ by (\ref{eq:willingnessToPay}), $\partial^2 \E[R]/\partial^2 f_\tau\leq0$ and $\E[R]$ is concave in $f_\tau$ for each $\tau \in\{0,1,\dots,\bar{\tau}\}$.

\textit{Part (2):} We derive the revenue-maximizing shipment policy. 
Recalling (\ref{eq:revenue_proof}), it follows that $\bar{\tau}^r=T-1$. Then, as $\E[R]$ is concave in $f_\tau$ by part (1) and together with (\ref{eq:revenueDerivative}), the revenue-maximizing fees $f^r_\tau$, $\tau \in\{0,1,\dots,\bar{\tau}\}$, are unique and satisfy
\begin{align*}
    \frac{\partial \E[R]}{\partial f_\tau} = 0 
    \quad \quad \quad \Rightarrow \quad \quad \quad
    f^r_\tau  = - \frac{w(f_\tau)}{w'(f_\tau)}.
\end{align*}
Applying (\ref{eq:willingnessToPay}), we have $f^r_\tau = f^r = \max\{\umax/2; \umin\}$ $\forall \tau \in\{0,1,\dots,\bar{\tau}\}$ and $\pi^r=(T-1,f^r,\dots,f^r)$ is the unique revenue-maximizing shipment policy. 
\hfill $\Box$

\proof{Proof of Lemma~\ref{lem:constantShipmentPolicy}.}
This proof has two parts. Part (1) proves claim (1) of Lemma~\ref{lem:constantShipmentPolicy} etc. 

\textit{Part (1):} 
We show that the expected backorders of the constant shipment policy are convex in the express shipment fee. Let $f_1,f_2\in[\umin,\umax]$, $f_1\neq f_2$, and $f_3:=\epsilon f_1+(1-\epsilon)f_2$, $\epsilon\in (0,1)$. We denote by $\pi^{c1},\pi^{c2},\pi^{c3}\in\Pi^c$ the associated constant shipment policies, i.e., $\pi^{c1}: = (T-1, f_1, \hdots, f_1)$,  $\pi^{c2}: = (T-1, f_2, \hdots, f_2)$, and  $\pi^{c3}: = (T-1, f_3, \hdots, f_3)$. 
As for every $\pi\in\Pi^c$, $C_\tau^\pi = (\tau+1)w(f)\lambda$, $\tau\in\{0,1,\dots,T-1\}$,  by (\ref{eq:demandProfileShipmentPolicy}) and $w(f)$ linear in $f$ by (\ref{eq:willingnessToPay}), we find
\begin{align*}
    f_3=\epsilon f_1+(1-\epsilon)f_2
    \quad\quad &\Leftrightarrow \quad\quad
    w(f_3)=\epsilon w(f_1) + (1-\epsilon)w(f_2) \\
    \quad\quad &\Leftrightarrow \quad\quad
    C_\tau^{\pi^{c3}} = \epsilon C_\tau^{\pi^{c1}} + (1-\epsilon) C_\tau^{\pi^{c2}} \quad \forall \tau\in\{0,1,\dots,T-1\}.
\end{align*}
Thus, the expected backorders are convex in $f$ if and only if they are convex in $\mathbf{C}$.

In the following, we show that the expected backorders are convex in $\mathbf{C}$. Recalling from Definition~\ref{def:demandProfile} that $C_\tau^\pi = \E[E^\pi[0,\tau]]$, $\tau\in\{0,1,\dots,T-1\}$, we have for every $\tau\in\{0,1,\dots,T-1\}$,  
\begin{equation*}
    E^{\pi^{c1}}[0,\tau] \sim \poi(C_\tau^{\pi^{c1}}) \quad \quad \quad \quad
    E^{\pi^{c2}}[0,\tau] \sim \poi(C_\tau^{\pi^{c2}}) \quad \quad \quad \quad 
    E^{\pi^{c3}}[0,\tau] \sim \poi(C_\tau^{\pi^{c3}}).
\end{equation*}
In particular, together with $ C_\tau^{\pi^{c3}} = \epsilon C_\tau^{\pi^{c1}} + (1-\epsilon) C_\tau^{\pi^{c2}}$,
\begin{equation*}
     E^{\pi^{c3}}[0,\tau] \sim \poi(\epsilon C_\tau^{\pi^{c1}} + (1-\epsilon) C_\tau^{\pi^{c2}}) \sim \epsilon E^{\pi^{c1}}[0,\tau] + (1-\epsilon) E^{\pi^{c2}}[0,\tau].
\end{equation*}
For every $\tau\in\{0,1,\dots,T-1\}$, $E^{\pi^{c3}}[0,\tau]$ is coupled to $E^{\pi^{c1}}[0,\tau]$ and $E^{\pi^{c2}}[0,\tau]$ through
\begin{equation}
    E^{\pi^{c3}}[0,\tau] = \epsilon E^{\pi^{c1}}[0,\tau] + (1-\epsilon) E^{\pi^{c2}}[0,\tau].
    \label{eq:cumulativeExpressOrdersConstantPolicyConvexCombination}
\end{equation}

Next, consider the cumulative low-urgency processing capacity $L[0,\tau]$ up to age $\tau$. Recall from (\ref{eq:lowUrgencyProcessingCapacityCumulative}) that $L[0,\tau] = \max_{\tau'\in\{0,\dots,\tau\}}\left(K[0,\tau']-E[0,\tau'] -X_0^h \right)^+$. Then, together with (\ref{eq:cumulativeExpressOrdersConstantPolicyConvexCombination}), we find that for every $\tau\in\{0,1,\dots,T-1\}$, 
\begin{align}
   L^{\pi^{c3}}[0,\tau] &\overset{}{=} \max_{\tau'\in\{0,\dots,\tau\}}\left(K[0,\tau']-E^{\pi^{c3}}[0,\tau'] -X_0^h \right)^+ \nonumber \\ 
   &\overset{}{=} \max_{\tau'\in\{0,\dots,\tau\}}\left(K[0,\tau']-\epsilon E^{\pi^{c1}}[0,\tau] - (1-\epsilon) E^{\pi^{c2}}[0,\tau] -X^h_0 \right)^+ \nonumber \\
    &\leq \epsilon \max_{\tau'\in\{0,\dots,\tau\}}\left(K[0,\tau']-E^{\pi^{c1}}[0,\tau] -X^h_0 \right)^+ + (1-\epsilon) \max_{\tau'\in\{0,\dots,\tau\}}\left(K[0,\tau']- E^{\pi^{c2}}[0,\tau] -X^h_0 \right)^+ \nonumber\\
    &\overset{}{=} \epsilon L^{\pi^{c1}}[0,\tau] + (1-\epsilon) L^{\pi^{c2}}[0,\tau].
    \label{eq:cumulativeLowUrgencyConstantPolicyConvexCombination}
\end{align}

Recalling from (\ref{eq:backorders_2}) that $B_l = X^h_0 + E[0,T-1] - K[0,T-1] + L[0,T-1]$, it follows together with (\ref{eq:cumulativeExpressOrdersConstantPolicyConvexCombination})-(\ref{eq:cumulativeLowUrgencyConstantPolicyConvexCombination}) 
\begin{align}
    B_l^{\pi^{c3}} 
    &\leq \epsilon \left(X^h_0 + E^{\pi^{c1}}[0,T-1] - K[0,T-1] + L^{\pi^{c1}}[0,T-1]\right) \nonumber \\
    &\quad \quad \quad\quad \quad + (1-\epsilon) \left(X^h_0 + E^{\pi^{c2}}[0,T-1] - K[0,T-1] + L^{\pi^{c2}}[0,T-1]\right) 
    \overset{}{=} \epsilon B_l^{\pi^{c1}} + (1-\epsilon) B_l^{\pi^{c2}}. \label{eq:backordersConstantPolicyConvexCombination}
\end{align}
As (\ref{eq:backordersConstantPolicyConvexCombination}) holds for every $l\in\mathbb{N}_0$ and every realization of random variables, it also holds in expectation, i.e., 
\begin{equation*}
     \E[B^{\pi^{c3}}]\leq \epsilon \E[B^{\pi^{c1}}] + (1-\epsilon) \E[B^{\pi^{c2}}].
\end{equation*}
Thus, $\E[B]$ is convex in $\mathbf{C}$ and therefore convex in $f$.

\textit{Part (2):}  We characterize the profit-maximizing shipment fee of the constant shipment policy. 
Let $\pi^c = (T-1, f, \hdots, f)\in\Pi^c$ and recall from (\ref{eq:ProfitVariable}) that $\E[G^{\pi^c}] = \E[R^{\pi^c}] - \beta\E[B^{\pi^c}]$. As $\E[R^{\pi^c}]$ is concave in $f$ by Lemma~\ref{lem:revenueMaximizingPolicy} and $\E[B^{\pi^c}]$ is convex in $f$ by part (1), it follows that $\E[G^{\pi^c}]$ is concave in $f$ and the profit-maximizing fee $f^*$ is the unique solution to
\begin{equation}
    \frac{\partial\E[R^{\pi^c}]}{\partial f} (f^*)= \beta \frac{\partial\E[B^{\pi^c}]}{\partial f}(f^*).
    \label{eq:optConditionConstantPolicy}
\end{equation}

Let $f^r$ and $f^b$ be the revenue-maximizing fee and the backorder-minimzing fee, respectively, i.e., $\frac{\partial\E[R^{\pi^c}]}{\partial f} (f^r) =0$ and $\frac{\partial\E[B^{\pi^c}]}{\partial f} (f^b) =0$. 
Then, together with (\ref{eq:optConditionConstantPolicy}),  $f^*\in[f^r,f^b]$ and $f^*$ increases in $\beta$.
\hfill \Halmos \endproof

\proof{Proof of Lemma~\ref{lem:constantShipmentPolicyCutoffPoint}.}
Let $\pi^c = (T-1, f,\dots, f)\in \Pi^c$ and $\pi^{cc}= (\bar{\tau}, f',\hdots, f')\in\Pi^{cc}$ as in Lemma~\ref{lem:constantShipmentPolicyCutoffPoint}. Together with (\ref{eq:demandProfileShipmentPolicy}), we find that the associated express shipment demand profiles $\mathbf{C}^{\pi^c}$ and $\mathbf{C}^{\pi^{cc}}$ satisfy
\begin{align}
    C_\tau^{\pi^c}-C_{\tau-1}^{\pi^c} &= \E[E_\tau^{\pi^c}] = w(f)\lambda &\forall &\tau\in\{0,1,\dots,T-1\}, \label{eq:demandProfileConstantShipmentPolicy_2}\\
    C_\tau^{\pi^{cc}}-C_{\tau-1}^{\pi^{cc}} &=  \E[E_\tau^{\pi^{cc}}] 
    = \begin{cases}
        w(f')\lambda & \text{if }\tau\leq \bar{\tau}, \\
        0 &\text{otherwise} \\
    \end{cases} &\forall &\tau\in\{0,1,\dots,T-1\}. \label{eq:demandProfileConstantShipmentPolicyCutoff_2} 
\end{align}
This proof has three parts. Part (1) establishes claim (1) of Lemma~\ref{lem:constantShipmentPolicyCutoffPoint} etc. 

\textit{Part (1):}  We show that $f>f'$. Recalling that $C_{T-1}^{\pi^{c}}=C_{T-1}^{\pi^{cc}}$ by assumption of Lemma~\ref{lem:constantShipmentPolicyCutoffPoint} and applying (\ref{eq:demandProfileConstantShipmentPolicy_2})-(\ref{eq:demandProfileConstantShipmentPolicyCutoff_2}), we find   
\begin{align}
    C_{T-1}^{\pi^{c}}=C_{T-1}^{\pi^{cc}}
    \quad \quad \overset{}{\Leftrightarrow} \quad \quad
    w(f)\lambda T = w(f')\lambda (\bar{\tau}+1) \label{eq:cumulativeDemandPerCycle} 
    \quad \quad \overset{}{\Leftrightarrow} \quad \quad
    w(f) = w(f')\frac{\bar{\tau}+1}{T}. 
\end{align}
Note that $(\bar{\tau}+1)/T < 1$ as $\bar{\tau} \leq T-2$ by assumption of Lemma~\ref{lem:constantShipmentPolicyCutoffPoint}. Then, it follows that $f>f'$ as $w$ is strictly decreasing in $f$ for $f\in[\umin,\umax]$ by (\ref{eq:willingnessToPay}) and $f,f'\in[\umin,\umax]$ by $\pi^c\in\Pi^c$ and $\pi^{cc}\in\Pi^{cc}$.

\textit{Part (2)}: We establish $\E[B^{\pi^{cc}}] \leq \E[B^{\pi^{c}}]$. 
Let $\phi^{cc}\in\Phi$ and $\phi^c\in\Phi$ denote the equivalent fee structures of shipment policies $\pi^{cc}$ and $\pi^c$, respectively. We know from Lemma~\ref{lem:feeStructure} that $\phi^{cc}$ and $\phi^c$ are the unique fee structures that satisfy $C_\tau^{\phi^{cc}} = C_\tau^{\pi^{cc}}$ and $C_\tau^{\phi^{c}} = C_\tau^{\pi^{c}}$ for every $\tau=0,1,\dots,T-1$, respectively. Thus, $\E[B^{\pi^{cc}}] \leq \E[B^{\pi^{c}}]\Leftrightarrow \E[B^{\phi^{cc}}] \leq \E[B^{\phi^{c}}]$. Further, we know from Proposition~\ref{prop:dominance} that $\E[B^{\phi^{cc}}] \leq \E[B^{\phi^{c}}]$ if $\mathbf{C}^{\phi^{cc}}$ and $\mathbf{C}^{\phi^{c}}$ satisfy $C_\tau^{\phi^{cc}}\geq C_\tau^{\phi^{c}}$ $\forall \tau\in\{0,1,\dots,T-2\}$ and $C_{T-1}^{\phi^{cc}}=C_{T-1}^{\phi^{c}}$.

Thus, to establish $\E[B^{\pi^{cc}}] \leq \E[B^{\pi^{c}}]$, it suffices to show that $\mathbf{C}^{\pi^{cc}}$ and $\mathbf{C}^{\pi^{c}}$ satisfy $C_\tau^{\pi^{cc}}\geq C_\tau^{\pi^{c}}$ $\forall \tau\in\{0,1,\dots,T-2\}$ and $C_{T-1}^{\pi^{cc}}=C_{T-1}^{\pi^{c}}$. 
Note that the latter condition holds by assumption of Lemma~\ref{lem:constantShipmentPolicyCutoffPoint}. 
To show that $C_\tau^{\pi^{cc}}\geq C_\tau^{\pi^{c}}$ $\forall \tau\in\{0,1,\dots,T-2\}$, we differentiate two cases:
\begin{enumerate}
    \item[(a)] For every $\tau\in\{0,1,\dots,\bar{\tau}\}$, we recall from  (\ref{eq:demandProfileConstantShipmentPolicy_2})-(\ref{eq:demandProfileConstantShipmentPolicyCutoff_2}) that $C_\tau^{\pi^{cc}}-C_{\tau-1}^{\pi^{cc}} = w(f')\lambda$ and $C_\tau^{\pi^{c}}-C_{\tau-1}^{\pi^{c}} = w(f)\lambda$. Note that $w(f')\lambda > w(f)\lambda$ as $w$ is strictly decreasing in $f$ for $f\in[\umin,\umax]$ by (\ref{eq:willingnessToPay}) and $f'<f$ by part (1) of this proof. Then, together with $C_{-1}^{\pi^{cc}} = C_{-1}^{\pi^{c}}=0$, it follows that $C_\tau^{\pi^{cc}} > C_\tau^{\pi^{c}}$ $\forall\tau\in\{0,1,\dots,\bar{\tau}\}$.
    \item[(b)] For every $\tau\in\{\bar{\tau}+1,\dots,T-2\}$, we recall from  (\ref{eq:demandProfileConstantShipmentPolicy_2})-(\ref{eq:demandProfileConstantShipmentPolicyCutoff_2}) that $C_\tau^{\pi^{cc}}-C_{\tau-1}^{\pi^{cc}} = 0$ and $C_\tau^{\pi^{c}}-C_{\tau-1}^{\pi^{c}} = w(f)\lambda$. It follows that $C_{\bar{\tau}}^{\pi^{cc}} = C_{\bar{\tau}+1}^{\pi^{cc}} = \dots = C_{T-1}^{\pi^{cc}}$ while $C_{\bar{\tau}}^{\pi^{c}} \leq C_{\bar{\tau}+1}^{\pi^{c}} \leq \dots \leq C_{T-1}^{\pi^{c}}$. As $C_{T-1}^{\pi^{cc}}=C_{T-1}^{\pi^{c}}$ by assumption of Lemma~\ref{lem:constantShipmentPolicyCutoffPoint}, we have  $C_\tau^{\pi^{cc}} \geq C_\tau^{\pi^{c}}$ $\forall\tau\in\{\bar{\tau}+1,\dots,T-1\}$.
\end{enumerate}

\textit{Part (3)}: We show that $\E[G^{\pi^{cc}}] \geq \E[G^{\pi^{c}}]$ if and only if $\beta \geq (f-f') C_\tau^{\pi^{cc}}(\E[B^{\pi^{c}}] -\E[B^{\pi^{cc}}])^{-1}$. 
Recalling from (\ref{eq:revenue})-(\ref{eq:ProfitVariable}) that $\E[G] = \E[R] -  \beta \E[B]$ where $\E[R] = \sum_{\tau=0}^{T-1}f_\tau\E[E_\tau]$ and applying (\ref{eq:demandProfileConstantShipmentPolicy_2})-(\ref{eq:demandProfileConstantShipmentPolicyCutoff_2}), we have  
\begin{align*}
   \E[G^{\pi^{cc}}] = f'w(f')\lambda(\bar{\tau}+1) - \beta \E[B^{\pi^{cc}}] 
   \quad \quad \quad \quad \text{and} \quad \quad \quad \quad 
   \E[G^{\pi^{c}}] = fw(f)\lambda T - \beta \E[B^{\pi^{c}}].
\end{align*}
Then, we find 
\begin{align*}
   \E[G^{\pi^{cc}}] \geq \E[G^{\pi^{c}}]
   \quad &\Leftrightarrow \quad
   \beta (\E[B^{\pi^{c}}] -\E[B^{\pi^{cc}}]) \geq fw(f)\lambda T - f'w(f')\lambda(\bar{\tau}+1) \\
   \quad &\overset{}{\Leftrightarrow} \quad
   \beta  \geq (f - f')C_{T-1}^{\pi^{cc}}(\E[B^{\pi^{c}}] -\E[B^{\pi^{cc}}])^{-1},
\end{align*}
where the last equivalence applies $w(f)\lambda T = w(f')\lambda(\bar{\tau}+1) = C_{T-1}^{\pi^{cc}}$ by (\ref{eq:cumulativeDemandPerCycle}).
\hfill \Halmos \endproof

\proof{Proof of Lemma~\ref{lem:twoLevelTimeDependentShipmentPolicy}.} 
Let $\pi^{2ti}= (\bar{\tau}, f', \dots, f', f'', \dots, f'')\in\Pi^{2t}$ and $\pi^{2td}= (\bar{\tau},  f'', \dots, f'', f', \dots, f')\in\Pi^{2t}$ as in Lemma~\ref{lem:twoLevelTimeDependentShipmentPolicy}. We denote by $\hat{\tau}^{i}$ and $\hat{\tau}^{d}$ the fee switching points of $\pi^{2ti}$ and $\pi^{2td}$, respectively. Together with (\ref{eq:demandProfileShipmentPolicy}), we find that the associated express shipment demand profiles $\mathbf{C}^{\pi^{2ti}}$ and $\mathbf{C}^{\pi^{2td}}$ satisfy:
\begin{align}
  C_\tau^{\pi^{2ti}}-C_{\tau-1}^{\pi^{2ti}} &=  \E[E_\tau^{\pi^{2ti}}] 
    = \begin{cases}
        w(f')\lambda & \text{if }\tau\leq \hat{\tau}^i, \\
        w(f'')\lambda & \text{if }\hat{\tau}^i<\tau\leq \bar{\tau}, \\
        0 &\text{otherwise} \\
    \end{cases} &\forall &\tau\in\{0,1,\dots,T-1\}, \label{eq:demandProfile2ti} \\
   C_\tau^{\pi^{2td}}-C_{\tau-1}^{\pi^{2td}} &=  \E[E_\tau^{\pi^{2td}}] 
    = \begin{cases}
        w(f'')\lambda & \text{if }\tau\leq \hat{\tau}^d, \\
        w(f')\lambda & \text{if }\hat{\tau}^d<\tau\leq \bar{\tau}, \\
        0 &\text{otherwise} \\
    \end{cases} &\forall &\tau\in\{0,1,\dots,T-1\}. \label{eq:demandProfile2td} 
\end{align}
Recalling that $C_{T-1}^{\pi^{2ti}}=C_{T-1}^{\pi^{2td}}$ by assumption of Lemma~\ref{lem:twoLevelTimeDependentShipmentPolicy} and applying (\ref{eq:demandProfile2ti})-(\ref{eq:demandProfile2td}), we note that
\begin{align} \label{eq:feeSwitchingPointRelation}
    C_{T-1}^{\pi^{2ti}}=C_{T-1}^{\pi^{2td}}
    \quad \quad \Leftrightarrow \quad \quad
    (\hat{\tau}^i+\hat{\tau}^d-\bar{\tau}+1)\left(w(f'')-w(f')\right) = 0 
    \quad \quad \Leftrightarrow \quad \quad
    \hat{\tau}^d = \bar{\tau}-\hat{\tau}^i-1. 
\end{align}
The last equivalence follows since $w(f'')\neq w(f')$ as $w$ is strictly decreasing in $f$ for $f\in[\umin,\umax]$, $f'\neq f''$, and $f',f''\in[\umin,\umax]$.

The proof proceeds as follows. In part (1), we show that $\E[B^{\pi^{2ti}}] \leq \E[B^{\pi^{2td}}]$. In part (2), we show that $\E[G^{\pi^{2ti}}] \geq \E[G^{\pi^{2td}}]$. 

\textit{Part (1)}: We establish $\E[B^{\pi^{2ti}}] \leq \E[B^{\pi^{2td}}]$. 
Let $\phi^{2ti}$ and $\phi^{2td}$ denote the equivalent fee structures of shipment policies $\pi^{2ti}$ and $\pi^{2td}$, respectively. We know from Lemma~\ref{lem:feeStructure} that $\phi^{2ti}$ and $\phi^{2td}$ are the unique fee structures that satisfy $C_\tau^{\phi^{2ti}} = C_\tau^{\pi^{2ti}}$ and $C_\tau^{\phi^{2td}} = C_\tau^{\pi^{2td}}$ for every $\tau=0,1,\dots,T-1$, respectively. Thus, $\E[B^{\pi^{2ti}}] \leq \E[B^{\pi^{2td}}]\Leftrightarrow\E[B^{\phi^{2ti}}] \leq \E[B^{\phi^{2td}}]$. 
Further, we know from Proposition~\ref{prop:dominance} that $\E[B^{\phi^{2ti}}] \leq \E[B^{\phi^{2td}}]$ if $\mathbf{C}^{\phi^{2ti}}$ and $\mathbf{C}^{\phi^{2td}}$ satisfy $C_\tau^{\phi^{2ti}}\geq C_\tau^{\phi^{2td}}$ $\forall \tau\in\{0,1,\dots,T-2\}$ and $C_{T-1}^{\phi^{2ti}}=C_{T-1}^{\phi^{2td}}$.

Thus, to establish $\E[B^{\pi^{2ti}}] \leq \E[B^{\pi^{2td}}]$, it suffices to show that $\mathbf{C}^{\pi^{2ti}}$ and $\mathbf{C}^{\pi^{2td}}$  satisfy $C_\tau^{\pi^{2ti}}\geq C_\tau^{\pi^{2td}}$ $\forall \tau\in\{0,1,\dots,T-2\}$ and $C_{T-1}^{\pi^{2ti}}=C_{T-1}^{\pi^{2td}}$. 
Note that the latter condition holds by assumption of Lemma~\ref{lem:twoLevelTimeDependentShipmentPolicy}. To show that $C_\tau^{\pi^{2ti}}\geq C_\tau^{\pi^{2td}}$ $\forall \tau\in\{0,1,\dots,T-2\}$, we suppose that $\hat{\tau}^i\leq\hat{\tau}^d$ and differentiate four cases:
\begin{enumerate}
\item[(a)] For every $\tau\in\{0,1,\dots,\hat{\tau}^i\}$, we have $C_\tau^{\pi^{2ti}} = C_{\tau-1}^{\pi^{2ti}} + w(f')\lambda$ and $C_\tau^{\pi^{2td}} = C_{\tau-1}^{\pi^{2td}} + w(f'')\lambda$  by (\ref{eq:demandProfile2ti})-(\ref{eq:demandProfile2td}). Note that $w(f') > w(f'')$ since $w$ is strictly decreasing in $f$ for $f\in[\umin,\umax]$ by (\ref{eq:willingnessToPay}) and $f' < f''$ by assumption of Lemma~\ref{lem:timeDependentShipmentPolicy}. Then, together with $C_{-1}^{\pi^{2ti}} = C_{-1}^{\pi^{2td}}=0$, it follows that $C_\tau^{\pi^{2ti}} > C_\tau^{\pi^{2td}}$ $\forall\tau\in\{0,1,\dots,\hat{\tau}^i\}$.
\item[(b)] For every $\tau\in\{\hat{\tau}^i+1,\dots,\hat{\tau}^d\}$, we have $C_\tau^{\pi^{2ti}} = C_{\tau-1}^{\pi^{2ti}} + w(f'')\lambda$ and $C_\tau^{\pi^{2td}} = C_{\tau-1}^{\pi^{2td}} + w(f'')\lambda$  by (\ref{eq:demandProfile2ti})-(\ref{eq:demandProfile2td}). Together with $C_{\hat{\tau}^i}^{\pi^{2ti}} > C_{\hat{\tau}^i}^{\pi^{2td}}$ by case (a), it follows that $C_\tau^{\pi^{2ti}} > C_\tau^{\pi^{2td}}$ $\forall\tau\in\{\hat{\tau}^i+1,\dots,\hat{\tau}^d\}$.
\item[(c)] For every $\tau\in\{\bar{\tau},\dots,T-2\}$, we have $C_\tau^{\pi^{2ti}} = C_{\tau+1}^{\pi^{2ti}}$ and $C_\tau^{\pi^{2td}} = C_{\tau+1}^{\pi^{2td}}$ by (\ref{eq:demandProfile2ti})-(\ref{eq:demandProfile2td}). As $C_{T-1}^{\pi^{2ti}} = C_{T-1}^{\pi^{2td}}$ by assumption of Lemma~\ref{lem:twoLevelTimeDependentShipmentPolicy}, it follows that $C_\tau^{\pi^{2ti}} = C_\tau^{\pi^{2td}}$ $\forall\tau\in\{\bar{\tau},\dots,T-2\}$.
\item[(d)] For every $\tau\in\{\hat{\tau}^d+1,\dots,\bar{\tau}-1\}$, we have $C_{\tau}^{\pi^{2ti}} = C_{\tau+1}^{\pi^{2ti}} - w(f'')\lambda$ and $C_{\tau}^{\pi^{2td}} = C_{\tau+1}^{\pi^{2td}} - w(f')\lambda$ by (\ref{eq:demandProfile2ti})-(\ref{eq:demandProfile2td}). Note that $w(f'') < w(f')$ since $w$ is strictly decreasing in $f$ for $f\in[\umin,\umax]$ by (\ref{eq:willingnessToPay}) and $f'' > f'$ by assumption of Lemma~\ref{lem:twoLevelTimeDependentShipmentPolicy}. Then, together with case (c), i.e., $C_{\bar{\tau}}^{\pi^{2ti}} = C_{\bar{\tau}}^{\pi^{2td}}$, it follows that $C_\tau^{\pi^{2ti}} > C_\tau^{\pi^{2td}}$ $\forall\tau\in\{\hat{\tau}^d+1,\dots,\bar{\tau}-1\}$.
\end{enumerate}

The converse case, i.e., $\hat{\tau}^i>\hat{\tau}^d$, is handled similarly. This completes the proof of part (1). 

\textit{Part (2)}: We establish $\E[G^{\pi^{2ti}}] \geq \E[G^{\pi^{2td}}]$. Recall from (\ref{eq:ProfitVariable}) that $\E[G] = \E[R]- \beta \E[B]$. 
We first show that $\E[R^{\pi^{2ti}}] = \E[R^{\pi^{2td}}]$. Recalling from (\ref{eq:revenue}) that $\E[R] = \sum_{\tau=0}^{T-1}f_\tau\E[E_\tau]$ and applying (\ref{eq:demandProfile2ti})-(\ref{eq:demandProfile2td}), we have
\begin{align*}
    \E[R^{\pi^{2ti}}] 
    &= (\hat{\tau}^i+1)f'w(f')\lambda + (\bar{\tau}-\hat{\tau}^i)f''w(f'')\lambda \\
    \E[R^{\pi^{2td}}] 
    &= (\hat{\tau}^d+1)f''w(f'')\lambda + (\bar{\tau}-\hat{\tau}^d)f'w(f')\lambda.   
\end{align*}
Applying $\hat{\tau}^d = \bar{\tau}-\hat{\tau}^i-1$ by (\ref{eq:feeSwitchingPointRelation}), we find 
\begin{align*}
    \E[R^{\pi^{2td}}] = (\hat{\tau}^i+1)f'w(f')\lambda + (\bar{\tau}-\hat{\tau}^i)f''w(f'')\lambda = \E[R^{\pi^{2ti}}].
\end{align*}

Then, together with $\E[B^{\pi^{2ti}}] \leq \E[B^{\pi^{2td}}]$ by part (2) of this proof, it follows that $\E[G^{\pi^{2ti}}] \geq \E[G^{\pi^{2td}}]$.
\hfill \Halmos \endproof

\proof{Proof of Lemma~\ref{lem:twoLevelTimeDependentShipmentPolicyComparedConstant}.} 
Let $\pi^{cc}= (\bar{\tau}, f,\dots, f)\in\Pi^{cc}$ and $\pi^{2t}= (\bar{\tau}, f^e, \dots, f^e, f^l, \dots, f^l)\in\Pi^{2t}$ as in Lemma~\ref{lem:twoLevelTimeDependentShipmentPolicyComparedConstant}. Together with (\ref{eq:demandProfileShipmentPolicy}), we find that the associated express shipment demand profiles $\mathbf{C}^{\pi^{cc}}$ and $\mathbf{C}^{\pi^{2t}}$ satisfy:
\begin{align}
   C_\tau^{\pi^{cc}}-C_{\tau-1}^{\pi^{cc}} &=  \E[E_\tau^{\pi^{cc}}] 
    = \begin{cases}
        w(f)\lambda & \text{if }\tau\leq \bar{\tau}, \\
        0 &\text{otherwise} \\
    \end{cases} &\forall &\tau\in\{0,1,\dots,T-1\}, \label{eq:demandProfilecc} \\
   C_\tau^{\pi^{2t}}-C_{\tau-1}^{\pi^{2t}} &=  \E[E_\tau^{\pi^{2t}}] 
    = \begin{cases}
        w(f^e)\lambda & \text{if }\tau\leq \hat{\tau}, \\
        w(f^l)\lambda & \text{if }\hat{\tau}<\tau\leq \bar{\tau}, \\
        0 &\text{otherwise} \\
    \end{cases} &\forall &\tau\in\{0,1,\dots,T-1\}. \label{eq:demandProfile2ti_2}
\end{align}
This proof has three parts. Part (1) establishes claim (1) of Lemma~\ref{lem:twoLevelTimeDependentShipmentPolicyComparedConstant} etc.

\textit{Part (1)}: We characterize the relationship of the express shipment fees $f$, $f^e$, and $f^l$. Recalling that $C_{T-1}^{\pi^{2t}}=C_{T-1}^{\pi^{cc}}$ by assumption of Lemma~\ref{lem:twoLevelTimeDependentShipmentPolicyComparedConstant} and applying (\ref{eq:demandProfilecc})-(\ref{eq:demandProfile2ti_2}), we find  
\begin{align}
    C_{T-1}^{\pi^{2t}}=C_{T-1}^{\pi^{cc}}
    \quad &\overset{}{\Leftrightarrow} \quad
    (\hat{\tau}+1)w(f^e)\lambda + (\bar{\tau}-\hat{\tau})w(f^l)\lambda = (\bar{\tau}+1)w(f)\lambda \nonumber \\
    \quad &\Leftrightarrow \quad
    \frac{\hat{\tau}+1}{\bar{\tau}+1}w(f^e) + \frac{\bar{\tau}-\hat{\tau}}{\bar{\tau}+1}w(f^l) = w(f) \nonumber \\
    \quad &\overset{}{\Leftrightarrow} \quad
    \frac{\hat{\tau}+1}{\bar{\tau}+1} f^e + \frac{\bar{\tau}-\hat{\tau}}{\bar{\tau}+1} f^l =f. \label{eq:feeComparison}
\end{align}
The last equivalence follows as $w$ is linear and strictly decreasing in $f$ for $f\in[\umin,\umax]$ by (\ref{eq:willingnessToPay}) and $f,f^e,f^l\in[\umin,\umax]$ by $\pi^{cc}\in\Pi^{cc}$ and $\pi^{2t}\in\Pi^{2t}$. 
As $f^e<f^l$ by assumption of Lemma~\ref{lem:twoLevelTimeDependentShipmentPolicyComparedConstant}, (\ref{eq:feeComparison}) implies that $f^e<f<f^l$. 

\textit{Part (2)}: We establish $\E[B^{\pi^{2t}}] \leq \E[B^{\pi^{cc}}]$. 
Let $\phi^{2ti}$ and $\phi^{cc}$ denote the equivalent fee structures of shipment policies $\pi^{2t}$ and $\pi^{cc}$, respectively. We know from Lemma~\ref{lem:feeStructure} that $\phi^{2ti}$ and $\phi^{cc}$ are the unique fee structures that satisfy $C_\tau^{\phi^{2ti}} = C_\tau^{\pi^{2t}}$ and $C_\tau^{\phi^{cc}} = C_\tau^{\pi^{cc}}$ for every $\tau=0,1,\dots,T-1$, respectively. Thus, $\E[B^{\pi^{2t}}] \leq \E[B^{\pi^{cc}}] \Leftrightarrow\E[B^{\phi^{2ti}}] \leq \E[B^{\phi^{cc}}]$. 
Further, we know from Proposition~\ref{prop:dominance} that $\E[B^{\phi^{2ti}}] \leq \E[B^{\phi^{cc}}]$ if $\mathbf{C}^{\phi^{2ti}}$ and $\mathbf{C}^{\phi^{cc}}$ satisfy $C_\tau^{\phi^{2ti}}\geq C_\tau^{\phi^{cc}}$ $\forall \tau\in\{0,1,\dots,T-2\}$ and $C_{T-1}^{\phi^{2ti}}=C_{T-1}^{\phi^{cc}}$.

Thus, to establish $\E[B^{\pi^{2t}}] \leq \E[B^{\pi^{cc}}]$, it suffices to show that $\mathbf{C}^{\pi^{2t}}$ and $\mathbf{C}^{\pi^{cc}}$ satisfy $C_\tau^{\pi^{2t}}\geq C_\tau^{\pi^{cc}}$ $\forall \tau\in\{0,1,\dots,T-2\}$ and $C_{T-1}^{\pi^{2t}}=C_{T-1}^{\pi^{cc}}$. 
Note that the latter condition holds by assumption of Lemma~\ref{lem:twoLevelTimeDependentShipmentPolicyComparedConstant}. 
To show that $C_\tau^{\pi^{2t}}\geq C_\tau^{\pi^{cc}}$ $\forall \tau\in\{0,1,\dots,T-2\}$, we differentiate three cases:
\begin{enumerate}
\item[(a)] For every $\tau\in\{0,1,\dots,\hat{\tau}\}$, we have $C_\tau^{\pi^{2t}} = C_{\tau-1}^{\pi^{2t}} + w(f^e)\lambda$ and $C_\tau^{\pi^{cc}} = C_{\tau-1}^{\pi^{cc}} + w(f)\lambda$  by (\ref{eq:demandProfilecc})-(\ref{eq:demandProfile2ti_2}). Note that $w(f^e) > w(f)$ since $w$ is strictly decreasing in $f$ for $f\in[\umin,\umax]$ by (\ref{eq:willingnessToPay}) and $f^e < f$ by part (1) of this proof. Then, together with $C_{-1}^{\pi^{2t}} = C_{-1}^{\pi^{cc}}=0$, it follows that $C_\tau^{\pi^{2t}} > C_\tau^{\pi^{cc}}$ $\forall\tau\in\{0,1,\dots,\hat{\tau}\}$.
\item[(b)] For every $\tau\in\{\bar{\tau},\dots,T-2\}$, we have $C_\tau^{\pi^{2t}} = C_{\tau+1}^{\pi^{2t}}$ and $C_\tau^{\pi^{cc}} = C_{\tau+1}^{\pi^{cc}}$ by (\ref{eq:demandProfilecc})-(\ref{eq:demandProfile2ti_2}). As $C_{T-1}^{\pi^{2t}} = C_{T-1}^{\pi^{cc}}$ by assumption of Lemma~\ref{lem:twoLevelTimeDependentShipmentPolicyComparedConstant}, it follows that $C_\tau^{\pi^{2t}} = C_\tau^{\pi^{cc}}$ $\forall\tau\in\{\bar{\tau},\dots,T-2\}$.
\item[(c)] For every $\tau\in\{\hat{\tau}+1,\dots,\bar{\tau}-1\}$, we have $C_{\tau}^{\pi^{2t}} = C_{\tau+1}^{\pi^{2t}} - w(f^l)\lambda$ and $C_{\tau}^{\pi^{cc}} = C_{\tau+1}^{\pi^{cc}} - w(f)\lambda$ by (\ref{eq:demandProfilecc})-(\ref{eq:demandProfile2ti_2}). Note that $w(f^l) < w(f)$ since $w$ is strictly decreasing in $f$ for $f\in[\umin,\umax]$ by (\ref{eq:willingnessToPay}) and $f^l > f$ by part (1) of this proof. Then, together with (b), i.e., $C_{\bar{\tau}}^{\pi^{2t}} = C_{\bar{\tau}}^{\pi^{cc}}$, it follows that $C_\tau^{\pi^{2t}} > C_\tau^{\pi^{cc}}$ $\forall\tau\in\{\hat{\tau}+1,\dots,\bar{\tau}-1\}$.
\end{enumerate}

\textit{Part (3)}: We show that $\E[G^{\pi^{2t}}]\geq \E[G^{\pi^{cc}}]$ if and only if $\beta \geq (\E[R^{\pi^{cc}}]-\E[R^{\pi^{2t}}])(\E[B^{\pi^{cc}}]-\E[B^{\pi^{2t}}])^{-1}$. 
This follows when recalling from (\ref{eq:ProfitVariable}) that $\E[G] = \E[R] -  \beta \E[B]$ since 
\begin{align*}
   \E[G^{\pi^{2t}}] \geq \E[G^{\pi^{cc}}]
   \quad &\Leftrightarrow \quad
   \beta \geq (\E[R^{\pi^{cc}}]-\E[R^{\pi^{2t}}])(\E[B^{\pi^{cc}}]-\E[B^{\pi^{2t}}])^{-1}.
\end{align*}
Note that $\E[R^{\pi^{cc}}]\geq\E[R^{\pi^{2t}}]$ as $\E[R]$ is concave in $f$ by Lemma~\ref{lem:revenueMaximizingPolicy} and $f = \epsilon f^e+(1-\epsilon)f^l$ with $\epsilon=(\hat{\tau}+1)/(\bar{\tau}+1)$ by part (1). Further, $\E[B^{\pi^{cc}}]\geq\E[B^{\pi^{2t}}]$ by part (2). 
\hfill \Halmos \endproof

\proof{Proof of Lemma~\ref{lem:timeDependentShipmentPolicy}.} 
Let $\pi = (\bar{\tau}, f_0, f_1, \dots, f_{\bar{\tau}})\in\Pi$, $\pi' = (\bar{\tau}, f'_0, f'_1, \dots, f'_{\bar{\tau}})\in\Pi$, and $\pi^{cc} = (\bar{\tau}, \bar{f},\dots, \bar{f})\in\Pi^{cc}$ as in Lemma~\ref{lem:timeDependentShipmentPolicy}, i.e., 
\begin{equation} \label{eq:conditionsFeesGeneralPolicies_1}
    \frac{1}{\tau} \sum_{j=0}^{\tau-1} f_j \leq \frac{1}{\tau} \sum_{j=0}^{\tau-1} f'_j \le \bar{f} \quad \quad \forall \tau \in\{1,2,\dots,\bar{\tau}\}
    \quad \quad \quad \quad \text{and} \quad \quad \quad \quad 
    \frac{1}{\bar{\tau}+1} \sum_{j=0}^{\bar{\tau}} f_j = \frac{1}{\bar{\tau}+1} \sum_{j=0}^{\bar{\tau}} f'_j = \bar{f}.
\end{equation}
We denote by $\phi$, $\phi'$, and  $\phi^{cc}$ the equivalent fee structures of shipment policies $\pi$, $\pi'$, and $\pi^{cc}$, respectively. 
We know from Lemma~\ref{lem:feeStructure} that $\phi$, $\phi'$, and $\phi^{cc}$ are the unique fee structures that satisfy for every $\tau=0,1,\dots,T-1$, $C_\tau^{\phi} = C_\tau^{\pi}$, $C_\tau^{\phi'} = C_\tau^{\pi'}$, and $C_\tau^{\phi^{cc}} = C_\tau^{\pi^{cc}}$, respectively. 
Thus, $\E[B^{\pi}]\leq \E[B^{\pi'}]\leq \E[B^{\pi^{cc}}] \Leftrightarrow\E[B^{\phi}]\leq \E[B^{\phi'}]\leq \E[B^{\phi^{cc}}]$. Further, we know from Proposition~\ref{prop:dominance} that $\E[B^{\phi}]\leq \E[B^{\phi'}]\leq \E[B^{\phi^{cc}}]$ if $\mathbf{C}^{\phi}$, $\mathbf{C}^{\phi'}$, and $\mathbf{C}^{\phi^{cc}}$ satisfy $C_\tau^{\phi}\geq C_\tau^{\phi'}\geq C_\tau^{\phi^{cc}}$ $\forall \tau\in\{0,1,\dots,T-2\}$ and $C_{T-1}^{\phi}=C_{T-1}^{\phi'}=C_{T-1}^{\phi^{cc}}$.
Thus, it suffices to show that 
\begin{equation}
   C_{T-1}^{\pi}=C_{T-1}^{\pi'}=C_{T-1}^{\pi^{cc}} \quad \quad \quad \quad \quad \text{and} \quad \quad \quad \quad \quad
    C_\tau^{\pi}\geq C_\tau^{\pi'}\geq C_{\tau}^{\pi^{cc}} \quad \quad \forall \tau\in\{0,1,\dots,T-2\} \label{eq:conditionsFeesGeneralPolicies_2}
\end{equation}

As the express shipment fees satisfy (\ref{eq:conditionsFeesGeneralPolicies_1}) and recalling from (\ref{eq:willingnessToPay}) that $w$ is linear and strictly decreasing in $f$ for $f\in[\umin,\umax]$, we find for every $\tau = 1, \dots, \bar{\tau}$:
\begin{align*}
     \frac{1}{\tau}\sum_{j=0}^{\tau-1} f_j \leq \frac{1}{\tau}\sum_{j=0}^{\tau-1} f'_j \leq \bar{f}
     \quad \quad \overset{}{\Leftrightarrow} \quad \quad 
     \lambda \sum_{j=0}^{\tau-1} w(f_j) \geq \lambda \sum_{j=0}^{\tau-1} w(f'_j) \geq \tau \lambda w(\bar{f})
    \quad \quad \overset{}{\Leftrightarrow} \quad \quad 
    C_{\tau-1}^{\pi}\geq C_{\tau-1}^{\pi'} \geq C_{\tau-1}^{\pi^{cc}}
\end{align*}
where the last equivalence follows by construction of $\mathbf{C}^{\pi}$, $\mathbf{C}^{\pi'}$, and $\mathbf{C}^{\pi^{cc}}$ (Definition~\ref{def:demandProfile}) together with $\E[E_\tau^{\pi^{cc}}] = \tau w(\bar{f})\lambda$, $\tau=0,\dots, \bar{\tau}$, 
Similarly, 
\begin{equation*}
    \frac{1}{\bar{\tau}+1}\sum_{j=0}^{\bar{\tau}} f_j = \frac{1}{\bar{\tau}+1}\sum_{j=0}^{\bar{\tau}} f'_j = \bar{f}
     \quad  \Leftrightarrow \quad 
     \lambda \sum_{j=0}^{\bar{\tau}} w(f_j) = \lambda \sum_{j=0}^{\bar{\tau}} w(f'_j) = (\bar{\tau}+1)\lambda w(\bar{f})
     \quad \Leftrightarrow \quad 
     C_{\bar{\tau}}^{\pi}= C_{\bar{\tau}}^{\pi'} = C_{\bar{\tau}}^{\pi^{cc}}. 
\end{equation*}
Then, recalling from (\ref{eq:demandProfileShipmentPolicy}) that $C_\tau = C_{\tau+1}$ for $\tau=\bar{\tau}, \dots, T-2$, we have $C_{\tau}^{\pi}= C_{\tau}^{\pi'} = C_{\tau}^{\pi^{cc}}$ $\forall\tau\in\{\bar{\tau}, \dots, T-1\}$. Thus, $\pi$, $\pi'$, and $\pi^{cc}$ satisfy (\ref{eq:conditionsFeesGeneralPolicies_2}). 
\hfill \Halmos \endproof

\proof{Proof of Lemma~\ref{lem:monotonicityFee}.}
Let $\pi^* = (\bar{\tau}, f_0, f_1, \dots, f_{\bar{\tau}})$ be any profit-maximizing shipment policy. We construct shipment policy $\pi'= (\bar{\tau}, f'_0, f'_1, \dots, f'_{\bar{\tau}})$ such that the fees $f'_\tau$, $\tau=0,1,\dots,\bar{\tau}$, correspond to the fees $f_\tau$, $\tau=0,1,\dots,\bar{\tau}$, but $f'_\tau$, $\tau=0,1,\dots,\bar{\tau}$, are ordered from low to high, i.e., $f'_\tau\geq f'_{\tau-1}$ $\forall\tau=1,\dots,\bar{\tau}$. 
As 
\begin{equation*}
    \frac{1}{\bar{\tau}+1}\sum_{j=0}^{\bar{\tau}}f'_j = \frac{1}{\bar{\tau}+1}\sum_{j=0}^{\bar{\tau}}f_j
    \quad \quad \quad \quad \text{and} \quad \quad \quad \quad 
    \frac{1}{\tau}\sum_{j=0}^{\tau-1}f'_j \leq \frac{1}{\tau}\sum_{j=0}^{\tau-1}f_j \quad \quad \forall \tau\in\{1,2,\dots,\bar{\tau}\},
\end{equation*}
we know from Lemma~\ref{lem:timeDependentShipmentPolicy} that $\E[B^{\pi'}] \leq \E[B^{\pi^*}]$. 
Further, together with (\ref{eq:revenue}), we have
\begin{equation*}
    \E[R^{\pi'}] = \sum_{j=0}^{\bar{\tau}} f'_\tau \E[E_\tau^{\pi'}] = \sum_{j=0}^{\bar{\tau}} f_\tau \E[E_\tau^{\pi^*}] = \E[R^{\pi^*}].
\end{equation*}
Recalling from (\ref{eq:ProfitVariable}) that $\E[G] = \E[R]- \beta \E[B]$, it follows the existence of a profit-maximizing policy with monotonous increasing fees.
\hfill \Halmos \endproof

\subsection{Proofs related to optimizing shipment fees (Propositions~\ref{prop:submodularityShipmentPolicy}-\ref{prop:polytime})}
\label{sec:proofsOptimization}

We first establish a supermodularity result for the express shipment demand profiles of fee structures (Proposition~\ref{prop:submodularity}), which helps us to prove Proposition \ref{prop:submodularityShipmentPolicy}.

\begin{proposition}[Supermodularity of fee structures]\label{prop:submodularity}
Let $\phi,\phi'\in \Phi$ with associated express shipment demand profiles $\mathbf{C}^\phi$ and $\mathbf{C}^{\phi'}$.
\begin{enumerate}
    \item[(1)] There exist unique fee structures $\phi'',\phi'''\in\Phi$ such that $\mathbf{C}^{\phi''} =\min(\mathbf{C}^{\phi},\mathbf{C}^{\phi'})$ and $\mathbf{C}^{\phi'''}=\max(\mathbf{C}^{\phi},\mathbf{C}^{\phi'})$, where $\min$ and $\max$ are taken componentwise. 
    \item[(2)]  The expected backorders are submodular in $\mathbf{C}$, i.e., $\E[B^{\phi''}]+\E[B^{\phi'''}]\le \E[B^{\phi}]+\E[B^{\phi'}]$.
    \item[(3)] The expected profit is supermodular in $\mathbf{C}$, i.e., $\E[G^{\phi''}]+\E[G^{\phi'''}]\ge \E[G^{\phi}]+\E[G^{\phi'}]$. 
\end{enumerate}
\end{proposition}

\proof{Proof of Proposition~\ref{prop:submodularity}.} 
This proof has three parts. Part (1) establishes claim (1) of Proposition~\ref{prop:submodularity} etc.

\textit{Part (1)}:
Recall from (\ref{eq:demandProfileFeeStructure}) that the express shipment demand profile $\mathbf{C}^\phi$ of any $\phi=(f_0,f_1,\dots,f_{T-1},\umax)\in \Phi$ satisfies $C^{\phi}_\tau-C^{\phi}_{\tau-1}=\E[E^{\phi}_\tau]= w({f}_\tau) \lambda \in [0,\lambda]$  $\forall \tau\in\{0,1,\dots,T-1\}$. Then, it is easy to see that $C_\tau-C_{\tau-1}\in [0,\lambda]$ $\forall \tau\in\{0,1,\dots,T-1\}$ is also a sufficient condition for the existence of a fee structure $\phi\in\Phi$ that generates demand profile $\mathbf{C}$. 

We build on this observation to prove claim (1).  Let $\phi,\phi' \in \Phi$ with associated express shipment demand profiles $\mathbf{C}^\phi$ and $\mathbf{C}^{\phi'}$. Further, $\mathbf{C}^{\phi''} =\min(\mathbf{C}^{\phi},\mathbf{C}^{\phi'})$ and $\mathbf{C}^{\phi'''}=\max(\mathbf{C}^{\phi},\mathbf{C}^{\phi'})$ as in Proposition~\ref{prop:submodularity}. 
Together with the previous observation, we know that there exist unique $\phi'',\phi'''\in\Phi$ if
\begin{align}
    C_\tau^{\phi''} - C_{\tau-1}^{\phi''} &= \min(C^{\phi}_{\tau},C^{\phi'}_{\tau})-\min(C^{\phi}_{\tau-1},C^{\phi'}_{\tau-1}) \in [0,\lambda] \quad \quad \quad \quad &\forall &\tau\in\{0,1,\dots,T-1\} \label{eq:conditionMin}\\
    C_\tau^{\phi'''} - C_{\tau-1}^{\phi'''} &= \max(C^{\phi}_{\tau},C^{\phi'}_{\tau})-\max(C^{\phi}_{\tau-1},C^{\phi'}_{\tau-1}) \in [0,\lambda]\quad \quad \quad \quad &\forall &\tau\in\{0,1,\dots,T-1\}.\label{eq:conditionMax}
\end{align}

To show that (\ref{eq:conditionMin})-(\ref{eq:conditionMax}) hold, we fix some $\tau\in\{0,1,\dots,T-1\}$ and differentiate the following four cases:
\begin{enumerate}
    \item[(a)] Suppose that $C^{\phi'}_{\tau} \ge C^{\phi}_{\tau}$ and $C^{\phi'}_{\tau-1} \ge C^{\phi}_{\tau-1}$. Then, we have $\min(C^{\phi}_{\tau},C^{\phi'}_{\tau})-\min(C^{\phi}_{\tau-1},C^{\phi'}_{\tau-1})=C^{\phi}_{\tau}-C^{\phi}_{\tau-1}\in[0,\lambda]$, which follows since $\phi\in\Phi$. Similarly, $\max(C^{\phi}_{\tau},C^{\phi'}_{\tau})-\max(C^{\phi}_{\tau-1},C^{\phi'}_{\tau-1})=C^{\phi'}_{\tau}-C^{\phi'}_{\tau-1}\in [0,\lambda]$ since $\phi'\in\Phi$.
    \item[(b)] Suppose that $C^{\phi'}_{\tau} \ge C^{\phi}_{\tau}$ and $C^{\phi'}_{\tau-1} \le C^{\phi}_{\tau-1}$. We have $\min(C^{\phi}_{\tau},C^{\phi'}_{\tau})-\min(C^{\phi}_{\tau-1},C^{\phi'}_{\tau-1})=C^{\phi}_{\tau} - C^{\phi'}_{\tau-1}$. 
    Note that $C^{\phi}_{\tau} - C^{\phi'}_{\tau-1} \overset{}{\le} C^{\phi'}_{\tau} - C^{\phi'}_{\tau-1} \overset{}{\le} \lambda$, where the first inequality applies $C^{\phi'}_{\tau} \ge C^{\phi}_{\tau}$ and the second follows as $\phi'\in\Phi$. 
    Further, $C^{\phi}_{\tau} - C^{\phi'}_{\tau-1} \overset{}{\ge} C^{\phi}_{\tau} - C^{\phi}_{\tau-1} \overset{}{\ge} 0$, where the first inequality applies $C^{\phi'}_{\tau-1} \le C^{\phi}_{\tau-1}$ and the second follows as $\phi\in\Phi$. Thus,  $\min(C^{\phi}_{\tau},C^{\phi'}_{\tau})-\min(C^{\phi}_{\tau-1},C^{\phi'}_{\tau-1})\in[0,\lambda]$.     
     We also have $\max(C^{\phi}_{\tau},C^{\phi'}_{\tau})-\max(C^{\phi}_{\tau-1},C^{\phi'}_{\tau-1})=C^{\phi'}_{\tau} - C^{\phi}_{\tau-1}\in[0,\lambda]$ since $C^{\phi'}_{\tau} - C^{\phi}_{\tau-1}\le C^{\phi'}_{\tau} - C^{\phi'}_{\tau-1}\le \lambda$ and $C^{\phi'}_{\tau} - C^{\phi}_{\tau-1}\ge C^{\phi}_{\tau} - C^{\phi}_{\tau-1}\ge 0$ by case assumptions and $\phi,\phi'\in \Phi$.
     \item[(c)] Suppose that $C^{\phi'}_{\tau} \le C^{\phi}_{\tau}$ and $C^{\phi'}_{\tau-1} \le C^{\phi}_{\tau-1}$. Then, (\ref{eq:conditionMin})-(\ref{eq:conditionMax}) hold by a similar argument as in case (a).
     \item[(d)] Suppose that $C^{\phi'}_{\tau} \le C^{\phi}_{\tau}$ and $C^{\phi'}_{\tau-1} \ge C^{\phi}_{\tau-1}$. Then, (\ref{eq:conditionMin})-(\ref{eq:conditionMax}) hold by a similar argument as in case (b).
\end{enumerate}

\textit{Part (2)}: 
Let $\phi,\phi',\phi'',\phi''' \in \Phi$ with associated express shipment demand profiles $\mathbf{C}^\phi$, $\mathbf{C}^{\phi'}$, $\mathbf{C}^{\phi''}$, $\mathbf{C}^{\phi'''}$ as in Proposition~\ref{prop:submodularity}. We establish that the expected backorders are submodular in the express shipment demand profile $\mathbf{C}$, i.e., 
\begin{equation}\label{eq:submodularityBackorders}
    \E[B^{\phi''}_l]+\E[B^{\phi'''}_l]\le \E[B^{\phi}_l]+\E[B^{\phi'}_l].
\end{equation}

We start with the cumulative express shipment orders $E[0,\tau]$ associated with $\phi,\phi',\phi'',\phi'''$. Recall from Definition~\ref{def:demandProfile} that for any $\phi\in \Phi$,  $C_\tau^\phi = \E[E^\phi[0,\tau]]$ $\forall \tau\in\{0,1,\dots,T-1\}$. Thus, for every $\tau\in\{0,1,\dots,T-1\}$,
\begin{align}
    E^{\phi}[0,\tau]&\sim \poi(C^\phi_\tau) \quad \quad 
    &E^{\phi'}[0,\tau]&\sim \poi(C^{\phi'}_\tau) \quad \quad \label{eq:distributionE_1}\\
    E^{\phi''}[0,\tau]&\sim \poi(\min(C^{\phi}_\tau,C^{\phi'}_\tau)) \quad \quad
    &E^{\phi'''}[0,\tau]&\sim \poi(\max(C^{\phi}_\tau,C^{\phi'}_\tau)). \label{eq:distributionE_2}
\end{align}
Let $E[0,\tau]\sim \poi(\min(C^{\phi}_{\tau},C^{\phi'}_{\tau}))$, $\chi_\tau^{\phi}\sim\poi((C^{\phi}_{\tau}-C^{\phi'}_{\tau})^+)$, and $\chi^{\phi'}_\tau\sim\poi((C^{\phi'}_{\tau}-C^{\phi}_{\tau})^+)$. Note that we understand $E[0,\tau]$ to be independent of both $\chi^{\phi'}_\tau$ and $\chi^{\phi}_\tau$, and $\chi\sim \poi(0)$ to mean $\chi:=0$.
Then, it follows together with (\ref{eq:distributionE_1})-(\ref{eq:distributionE_2}) that for every $\tau\in\{0,1,\dots,T-1\}$,
\begin{align}
    E^{\phi}[0,\tau]&\sim  E[0,\tau] +\chi^\phi_\tau \quad \quad 
    &E^{\phi'}[0,\tau]&\sim E[0,\tau] +\chi^{\phi'}_\tau \quad \quad \label{eq:distributionE_3} \\
    E^{\phi''}[0,\tau]&\sim E[0,\tau] \quad \quad
    &E^{\phi'''}[0,\tau]&\sim E[0,\tau]+\chi_\tau^\phi+\chi_\tau^{\phi'}. \label{eq:distributionE_4}
\end{align}
Further, for every $\tau\in\{0,1,\dots,T-1\}$, 
\begin{equation} \label{eq:cumulativeExpressDemandFeesStructures}
    E^{\phi}[0,\tau] + E^{\phi'}[0,\tau] = E[0,\tau] +\chi^\phi_\tau + E[0,\tau] +\chi^{\phi'}_\tau = E^{\phi''}[0,\tau] + E^{\phi'''}[0,\tau].
\end{equation}

Next, consider the cumulative low-urgency processing capacity $L[0,T-1]$ associated with $\phi,\phi',\phi'',\phi'''$. Let
\begin{equation*}
    \tau_1 :=\argmax_{\tau\in\{0,\ldots,T-1\}} \left( K[0,\tau]-E[0,\tau]-\chi^\phi_\tau-\chi^{\phi'}_\tau-X^h_0\right)^+ \quad \quad \quad 
    \tau_2:= \argmax_{\tau\in\{0,\ldots,T-1\}} \left( K[0,\tau]-E[0,\tau]-X^h_0\right)^+.
\end{equation*}
Recalling from (\ref{eq:lowUrgencyProcessingCapacityCumulative}) that $L[0,T-1]=\max_{\tau\in\{0,\ldots,T-1\}} \left( K[0,\tau]-E[0,\tau]-X^h_0\right)^+$, we have 
\begin{align*}
L^{\phi'''}[0,T-1]&+L^{\phi''}[0,T-1]\\
&\overset{}{=}\max_{\tau\in\{0,\ldots,T-1\}} \left( K[0,\tau]-E^{\phi'''}[0,\tau]-X^h_0\right)^+ + \max_{\tau\in\{0,\ldots,T-1\}} \left( K[0,\tau]-E^{\phi''}[0,\tau]-X^h_0\right)^+\\ 
&\overset{[2]}{=} \max_{\tau\in\{0,\ldots,T-1\}} \left( K[0,\tau]-E[0,\tau]-\chi^\phi_\tau-\chi^{\phi'}_\tau-X^h_0\right)^+ + \max_{\tau\in\{0,\ldots,T-1\}} \left( K[0,\tau]-E[0,\tau]-X^h_0\right)^+\\
&\overset{[3]}{=} \left( K[0,\tau_1]-E[0,\tau_1]-\chi^\phi_{\tau_1}-\chi^{\phi'}_{\tau_1}-X^h_0\right)^+ + \left( K[0,\tau_2]-E[0,\tau_2]-X^h_0\right)^+.
\end{align*}
Equality [2] applies (\ref{eq:distributionE_4}); equality [3] follows by construction of $\tau_1$ and $\tau_2$. 
Now, consider the case where $E^\phi[0,\tau_2]\ge E^{\phi'}[0,\tau_2] \Leftrightarrow C_{\tau_2}^{\phi}\ge C_{\tau_2}^{\phi'}$. Then, recalling that $\chi^{\phi'}_{\tau_2}\sim\poi((C^{\phi'}_{\tau_2}-C^{\phi}_{\tau_2})^+)$, we have $\chi^{\phi'}_{\tau_2}=0$ and find:
\begin{align*}
&\left( K[0,\tau_1]-E[0,\tau_1]-\chi^\phi_{\tau_1}-\chi^{\phi'}_{\tau_1}-X^h_0\right)^+ + \left( K[0,\tau_2]-E[0,\tau_2]-X^h_0\right)^+\\
&\quad \quad = \left( K[0,\tau_1]-E[0,\tau_1]-\chi^{\phi}_{\tau_1}-\chi^{\phi'}_{\tau_1}-X^h_0\right)^+ + \left( K[0,\tau_2]-E[0,\tau_2]-\chi^{\phi'}_{\tau_2}-X^h_0\right)^+\\
&\quad \quad \le\max_{\tau\in\{0,\ldots,T-1\}} \left( K[0,\tau]-E[0,\tau]-\chi^\phi_\tau-X^h_0\right)^+ + \max_{\tau\in\{0,\ldots,T-1\}} \left( K[0,\tau]-E[0,\tau]-\chi^{\phi'}_\tau-X^h_0\right)^+\\
&\quad \quad =\max_{\tau\in\{0,\ldots,T-1\}} \left( K[0,\tau]-E^{\phi}[0,\tau]-X^h_0\right)^+ + \max_{\tau\in\{0,\ldots,T-1\}} \left( K[0,\tau]-E^{\phi'}[0,\tau]-X^h_0\right)^+\\
&\quad \quad =L^\phi[0,T-1]+L^{\phi'}[0,T-1]. 
\end{align*}
The converse case, i.e., $E^\phi[0,\tau_2]\le E^{\phi'}[0,\tau_2]\Leftrightarrow C_{\tau_2}^{\phi}\le C_{\tau_2}^{\phi'}$, implies $\chi^{\phi}_{\tau_2}=0$ and is handled similarly. Thus, 
\begin{equation} \label{eq:cumulativeLowUrgencyCapacityFeeStructures}
    L^{\phi'''}[0,T-1]+L^{\phi''}[0,T-1] \leq L^\phi[0,T-1]+L^{\phi'}[0,T-1].
\end{equation}

Then, recalling from (\ref{eq:backorders_2}) the backorders $B_l=X^h_0 + E[0,T-1] - K[0,T-1] + L[0,T-1]$, we have 
\begin{align}
    B^{\phi''}_l+B^{\phi'''}_l-B^{\phi}_l-B^{\phi'}_l 
    &\overset{[1]}{=} E^{\phi''}[0,T-1]+E^{\phi'''}[0,T-1]-E^{\phi}[0,T-1]-E^{\phi'}[0,T-1]\nonumber\\
    &\quad \quad \quad \quad +L^{\phi''}[0,T-1]+L^{\phi'''}[0,T-1]-L^{\phi}[0,T-1]-L^{\phi'}[0,T-1] \nonumber\\
    &\overset{[2]}{=} L^{\phi''}[0,T-1]+L^{\phi'''}[0,T-1]-L^{\phi}[0,T-1]-L^{\phi'}[0,T-1] \nonumber\\
    &\overset{[3]}{\le} 0. \label{eq:backordersFeeStructures}
\end{align}
Equality [1] follows as $X^h_0$ and $K[0,T-1]$ are independent of $\phi$; equality [2] applies (\ref{eq:cumulativeExpressDemandFeesStructures}); inequality [3] follows together with (\ref{eq:cumulativeLowUrgencyCapacityFeeStructures}). As (\ref{eq:backordersFeeStructures}) holds for every realization of random variables, it also holds in expectation, which establishes (\ref{eq:submodularityBackorders}).

\textit{Part (3)}: Let $\phi,\phi',\phi'',\phi''' \in \Phi$ with associated express shipment demand profiles $\mathbf{C}^\phi$, $\mathbf{C}^{\phi'}$, $\mathbf{C}^{\phi''}$, $\mathbf{C}^{\phi'''}$ as in Proposition~\ref{prop:submodularity}. We show that the expected profit is supermodular in the express shipment demand profile $\mathbf{C}$, i.e., $\E[G^{\phi''}]+\E[G^{\phi'''}]\ge \E[G^{\phi}]+\E[G^{\phi'}]$. 

Recall from (\ref{eq:ProfitVariable}) that $\E[G] = \E[R]  - \beta \E[B]$. We further know from part (2) of this proof that $\E[B]$ is submodular in $\mathbf{C}$, i.e., $- \beta \E[B]$ is supermodular in $\mathbf{C}$. Thus, it remains to show that $\E[R]$ is supermodular in $\mathbf{C}$, i.e., $\E[R^{\phi''}]+\E[R^{\phi'''}]\ge \E[R^{\phi}]+\E[R^{\phi'}]$. 
As for any $\phi=(f_0, f_1, \dots, f_{T-1}, \umax) \in\Phi$, $\E[R^\phi] = \sum_{\tau=0}^{{T-1}} f_{\tau} \E[E^\phi_\tau] = \sum_{\tau=0}^{{T-1}} f_{\tau} (C^\phi_\tau-C^\phi_{\tau-1})$, it suffices to show that 
\begin{equation}
 f''_{\tau} (C^{\phi''}_\tau -C^{\phi''}_{\tau-1})+f'''_{\tau} (C^{\phi'''}_\tau -C^{\phi'''}_{\tau-1}) \geq f_{\tau} (C^\phi_\tau -C^\phi_{\tau-1})+f'_{\tau} (C^{\phi'}_{\tau} -C^{\phi'}_{\tau-1}) \quad \quad \forall \tau\in\{0,1,\dots,T-1\}.  \label{eq:tobeprovedforsubmodofprices}  
\end{equation}

We therefore fix some $\tau\in\{0,1,\dots,T-1\}$ and differentiate the following four cases:
\begin{enumerate}
    \item[(a)] Suppose that $C^{\phi'}_\tau\ge C^{\phi}_\tau$ and $C^{\phi'}_{\tau-1}\ge C^{\phi}_{\tau-1}$. Then, we have $C^{\phi'''}_{\tau}-C^{\phi'''}_{\tau-1}=C^{\phi'}_{\tau}-C^{\phi'}_{\tau-1}$. Note that $C^{\phi'''}_{\tau}-C^{\phi'''}_{\tau-1}=C^{\phi'}_{\tau}-C^{\phi'}_{\tau-1} \overset{}{\Leftrightarrow}  w(f'''_\tau) \lambda = w(f'_\tau)\lambda\overset{}{\Leftrightarrow} f'''_\tau = f'_\tau$, where the first equivalence applies (\ref{eq:demandProfileFeeStructure}) and the second equivalence follows as $w$ is strictly decreasing in $f$ for $f\in[\umin,\umax]$ and $f'_\tau,f'''_\tau\in[\umin,\umax]$. Similarly, we have $C^{\phi''}_{\tau}-C^{\phi''}_{\tau-1}=C^{\phi}_{\tau}-C^{\phi}_{\tau-1}$ and thus $f''_{\tau}=f_{\tau}$. Together, this implies that \eqref{eq:tobeprovedforsubmodofprices} holds with equality in case (a). 
    \item[(b)] Suppose that $C^{\phi'}_\tau\ge C^{\phi}_\tau$ and $C^{\phi'}_{\tau-1}< C^{\phi}_{\tau-1}$. Then, we have $C^{\phi'''}_\tau=C^{\phi'}_\tau$, $C^{\phi'''}_{\tau-1}=C^{\phi}_{\tau-1}$, $C^{\phi''}_\tau=C^{\phi}_\tau$, and $C^{\phi''}_{\tau-1}=C^{\phi'}_{\tau-1}$. It follows that
    \begin{align}
        \E[E^\phi_\tau]+\E[E^{\phi'}_\tau]=C_\tau^{\phi}-C_{\tau-1}^{\phi}+C_\tau^{\phi'}-C_{\tau-1}^{\phi'} &= \E[E^{\phi''}_\tau]+\E[E^{\phi'''}_\tau] \label{eq:orderIncomeExpressFeesStructures_1}\\
        C_\tau^{\phi}-C_{\tau-1}^{\phi} \le C_\tau^{\phi'}-C_{\tau-1}^{\phi} \le C_\tau^{\phi'}-C_{\tau-1}^{\phi'} \quad \quad &\Leftrightarrow \quad \quad \E[E^{\phi}_\tau] \le \E[E^{\phi''}_\tau] \le \E[E^{\phi'}_\tau] \label{eq:orderIncomeExpressFeesStructures_2}\\
        C_\tau^{\phi}-C_{\tau-1}^{\phi} \le C_\tau^{\phi}-C_{\tau-1}^{\phi'} \le C_\tau^{\phi'}-C_{\tau-1}^{\phi'} \quad \quad &\Leftrightarrow \quad \quad \E[E^{\phi}_\tau] \le \E[E^{\phi'''}_\tau] \le \E[E^{\phi'}_\tau]. \label{eq:orderIncomeExpressFeesStructures_3}
    \end{align}
    Further, note that $\E[E^\phi_\tau] = w(f_\tau)\lambda \Leftrightarrow f_\tau=w^{-1}(\E[E^\phi_{\tau}]/\lambda)$, where $w^{-1}:[0,1] \rightarrow [\umin,\umax]$ is the inverse of $w$. As $w^{-1}$ is linear and strictly decreasing in $\E[E^\phi_{\tau}]$, the term $f_{\tau} \E[E^\phi_\tau]=(w^{-1}(\E[E^\phi_{\tau}]/\lambda))\E[E^\phi_\tau]$ is concave in $\E[E^\phi_{\tau}]$. Together with (\ref{eq:orderIncomeExpressFeesStructures_1})-(\ref{eq:orderIncomeExpressFeesStructures_3}), this implies that (\ref{eq:tobeprovedforsubmodofprices}) holds in case (b).
    \item[(c)] Suppose that $C^{\phi'}_{\tau} \le C^{\phi}_{\tau}$ and $C^{\phi'}_{\tau-1} \le C^{\phi}_{\tau-1}$. Then, (\ref{eq:tobeprovedforsubmodofprices}) holds by a similar argument as in case (a).
     \item[(d)] Suppose that $C^{\phi'}_{\tau} \le C^{\phi}_{\tau}$ and $C^{\phi'}_{\tau-1} > C^{\phi}_{\tau-1}$. Then, (\ref{eq:tobeprovedforsubmodofprices}) holds by a similar argument as in case (b).
\end{enumerate}
\hfill \Halmos \endproof

\proof{Proof of Proposition~\ref{prop:submodularityShipmentPolicy}.}
Proposition~\ref{prop:submodularityShipmentPolicy} follows by Proposition~\ref{prop:submodularity} when recalling from Lemma~\ref{lem:feeStructure} that each shipment policy $\pi\in\Pi$ has a unique equivalent fee structure $\phi\in\Phi$, i.e., $\pi$ and $\phi$ yield the same express shipment demand profile. 
\hfill \Halmos \endproof

\proof{Proof of Proposition \ref{prop:polytime}.} To utilize established theory for submodular function minimization (SFM) \cite[see e.g.][]{mccormick2005submodular}, we connect to the standard standard setting which involves optimization over subsets of a finite ground set. To connect our formulation to this standard setting, we represent an express shipment demand profile $\mathbf{C}$ as a subset of the finite ground set $\mathcal{U} := \{ ( \tau, m ) \mid \tau \in \{0,1,\dots,T-1\}, \ m \in \{1,\dots,(\tau+1) n_f\} \},$ where the element $(\tau,m)$ indicates that the cumulative expected express shipment orders at age~$\tau$ are at least $m\lambda/n_f$.  
Given a profile $\mathbf{C}\in\mathcal{C}$, the corresponding subset is
\[
S(\mathbf{C}) := \{ (\tau,m) \in \mathcal{U} \mid C_\tau \ge m\lambda/n_f \}.
\]

Now, define the subsets of $\mathcal{U}$ that correspond to express shipment demand profiles $\mathbf{C}\in\mathcal{C}$ as follows: $\mathcal{R} \;:=\; \big\{ S(\mathbf{C}) \ \big|\ \mathbf{C}\in\mathcal{C} \big\} \;\subseteq\; 2^{\mathcal{U}}.$ A key observation is then that $\mathcal{R}$ is a \emph{ring family}, i.e., it is closed under union and intersection. Indeed, for any $\mathbf{C},\mathbf{C}'\in\mathcal{C}$, let
\[
(\mathbf{C}\vee\mathbf{C}')_\tau := \max\{C_\tau,C'_\tau\},\qquad
(\mathbf{C}\wedge\mathbf{C}')_\tau := \min\{C_\tau,C'_\tau\}, \quad \tau=0,\dots,T-1.
\]
Since $\mathcal{C}$ enforces non-decreasing express shipment demand profiles with increments in $[0,\lambda]$ and grid points in $\mathcal{C}_\tau$, both $\mathbf{C}\vee\mathbf{C}'$ and $\mathbf{C}\wedge\mathbf{C}'$ remain in $\mathcal{C}$. Moreover,
\[
S(\mathbf{C}) \cup S(\mathbf{C}') \;=\; S(\mathbf{C}\vee\mathbf{C}'), 
\qquad
S(\mathbf{C}) \cap S(\mathbf{C}') \;=\; S(\mathbf{C}\wedge\mathbf{C}'),
\]
so $\mathcal{R}$ is closed under $\cup$ and $\cap$.

On this domain, we define the set function
\[
F:\mathcal{R}\to\mathbb{R},\qquad F(S) := \E\big[ G(\mathbf{C}(S)) \big],
\]
where $\mathbf{C}(S)$ denotes the unique express shipment demand profile with $S(\mathbf{C}(S))=S$.

By these observations, and by construction in \S\ref{sec:optimizingshipmentfees}, any policy $\pi\in\Pi^{\mathcal{F}}$ corresponds to a unique demand profile $\mathbf{C}^\pi\in\mathcal{C}$, and hence to a unique set $S(\mathbf{C}^\pi)\in\mathcal{R}$, where $\mathcal{R}$ is a ring family over the ground set $\mathcal{U}$.  
The size of the ground set is 
\[
n_{\mathcal{U}} = |\mathcal{U}| = \frac{n_f\,T(T+1)}{2}.
\]
By Proposition~\ref{prop:submodularityShipmentPolicy}, the profit function $\pi \mapsto \E[G^\pi]$ is supermodular in $\mathbf{C}$ and hence in $S(\mathbf{C})$.  
Maximizing a supermodular function over a ring family is equivalent to minimizing a submodular function over the same ring family.  
By \citet[§3.3.2]{mccormick2005submodular}, this can be done in strongly polynomial time using only evaluations of $F(S) = \E[G(\mathbf{C}(S))]$ as an oracle.  
Therefore, the optimal policy $\pi^\star$ can be found in time $O\!\big(n_{\mathcal{U}}^{7} \log n_{\mathcal{U}} \cdot \EO\big)$, where $\EO$ is the time to evaluate $\E[G^\pi]$ for a given $\pi$.
\hfill \Halmos \endproof

\section{Implementation details}
\label{sec:AppendixImplementation}
In this section, we derive an upper bound for the state space of the Markov chain, which is needed to implement and compute the model. 
Recalling from \S\ref{sec:ModelMarkovChain} that $(X_t^h, X_t^\Sigma, \tau)$ with $\tau = t \Mod T$ specifies the state of the Markov chain, the state space is given by $\mathbb{N}_0\times\mathbb{N}_0\times\{0,1,\dots,T-1\}$. As $X_t^h\leq X_t^\Sigma$ for every $t\in\mathbb{N}_0$, it is sufficient to derive an upper bound for $X_t^\Sigma$, which we denote by $\bar{X}^\Sigma$.  

We determine $\bar{X}^\Sigma$ using a binary search algorithm such that the probability of rejecting a customer order is below a predefined threshold. For any given $\bar{X}^\Sigma\in \mathbb{N}$, a customer order arriving at any time $t$ is rejected if by accepting this order, the total order backlog $X_{t+1}^\Sigma$ at time $t+1$ would exceed the upper bound $\bar{X}^\Sigma$. The rejection probability is calculated as
\begin{equation*}
         \lim_{L\rightarrow \infty} \frac{1}{L} \sum_{t=0}^L \p \left(X^\Sigma_t + D_t - K_t > \bar{X}^\Sigma \right).
\end{equation*}

For the Markov chain with finite state space $\bar{X}^\Sigma\times\bar{X}^\Sigma\times\{0,1,\dots,T-1\}$, the state transition (\ref{eq:transitionHighUrgencyBacklog})-(\ref{eq:transitionTotalBacklog}) is adjusted as follows. Define $O_t := (X^\Sigma_t + D_t - K_t - \bar{X}^\Sigma)^+$ as the \emph{overflow of customer orders} at time $t$. If $O_t>0$, some customer orders are rejected at time $t$. We assume that regular shipment orders are rejected first. Then, the adjusted state transition from $(X^h_t,X^\Sigma_t)$ at time $t$ to $(X^h_{t+1},X^\Sigma_{t+1})$ at time $t+1$ is given by
\begin{align*}
     X_{t+1}^h &= \begin{cases}
            \left(X^\Sigma_t + D_t - K_t \right)^+    & \text{if }t\Mod T = T-1 \text{ and } O_t = 0, \\
            \bar{X}^\Sigma    & \text{if }t\Mod T = T-1 \text{ and } O_t > 0, \\
            \left(X^h_t + E_t - K_t \right)^+     &\text{if } t\Mod T < T-1 \text{ and } O_t = 0,   \\
            \left(X^h_t +\left( E_t - \left(O_t - N_t\right)^+ \right)^+- K_t \right)^+     &\text{otherwise},   \\
    \end{cases} \\
    X_{t+1}^\Sigma &= 
    \begin{cases}
        \left(X^\Sigma_t + D_t - K_t \right)^+  & \text{if }O_t = 0,\\
        \bar{X}^\Sigma                          & \text{otherwise}.
    \end{cases}
\end{align*}

\end{APPENDICES}

\end{document}